\documentclass[opre,nonblindrev]{informs3} 
\usepackage{longtable, comment, booktabs}
\usepackage{setspace}
\usepackage{pdfpages}
\usepackage{bbm}
\usepackage{algorithmicx}
\usepackage[]{algpseudocode}
\usepackage{multibib}
\usepackage{amsmath}
\usepackage{arydshln}
\usepackage{makecell}
\usepackage{textcomp}
\usepackage{subfigure}
\usepackage{wrapfig}
\usepackage{threeparttable}
\usepackage[ruled,vlined]{algorithm2e}
\newcites{annex}{References}
\OneAndAHalfSpacedXI 

\usepackage{amsmath, bm}
\usepackage{amssymb}
\usepackage{mathrsfs}
\usepackage{graphicx}
\usepackage{verbatim}
\usepackage{setspace}
\usepackage{float}
\usepackage{enumerate}
\usepackage{enumitem}
\usepackage{multirow}
\usepackage{array}
\usepackage{epstopdf}
\usepackage{hyperref}
\hypersetup{
	colorlinks=true,
	linkcolor=black,
	citecolor=black,
	filecolor=black,
	urlcolor=black,
}

\usepackage{booktabs} 
\usepackage{lineno}

\renewcommand{\arraystretch}{1.5} %

\newtheorem{myTheo}{Theorem}

\newtheorem{myFind}{\textbf{Finding}}
\newtheorem{myremark}{\textbf{Remark}}
\newtheorem{myPro}{Proposition}
\usepackage{natbib}
 \bibpunct[, ]{(}{)}{,}{a}{}{,}%
 \def\bibfont{\small}%
\TheoremsNumberedThrough

\EquationsNumberedThrough
{\theoremstyle{EXkey}}

\begin{document}

\TITLE{On Coordinated Drone-Courier Logistics for Intra-city Express Services}

\ARTICLEAUTHORS{%
\AUTHOR{Shuiwang Chen}
\AFF{Department of Aeronautical and Aviation, The Hong Kong Polytechnic University, Hong Kong, 999077, China} 
\AUTHOR{Kai Wang}
\AFF{School of Vehicle and Mobility, Tsinghua University, Beijing 100084, China}
\AUTHOR{Lingxiao Wu}
\AFF{Department of Aeronautical and Aviation, The Hong Kong Polytechnic University, Hong Kong, 999077, China, \EMAIL{lingxiao-leo.wu@polyu.edu.hk}} 
\AUTHOR{Wei Qi}
\AFF{Department of Industrial Engineering, Tsinghua University, Beijing 100084, China}
\AFF{Desautels Faculty of Management, McGill University, Montreal, Canada} 
} 

\ABSTRACT{%
\textbf{\textit{Problem definition}:} 
Drones, despite being acknowledged as a transformative force in the city logistics sector, are unable to execute the \textit{last-meter delivery} (unloading goods directly to customers' doorsteps) due to airspace restrictions and safety concerns.
To leverage advancements and overcome the limitations of drones in providing intra-city express services, we introduce a coordinated drone-courier logistics system where drones operate within a closed network among vertiports, while couriers connect customers to the drone delivery system.
This paper aims to shed light on this coordinated system in terms of system feasibility, network interactivity, and long-term sustainability.
\textbf{\textit{Methodology/Results}:} 
We develop an integrated optimization model to optimize the network planning of the coordinated logistics system. 
The interplay between network planning and tactical operations is mirrored by a queueing network model, resulting in the nonlinear and nonconvex (partially convex and partially concave) feasible region of the optimization model.
An iterative exact algorithm that tightens lower and upper bounds by adaptively refining the linear approximations of nonlinear constraints is developed to provide optimality-guaranteed solutions with finite convergence.
The computational experiments demonstrate the scalability and robustness of our algorithm across various network configurations and scenarios.
\textbf{\textit{Managerial implications:}}
The case study, based on a real-world dataset from SF Express, a logistics giant in China, validates that the coordinated logistics system efficiently attains cost and time savings by leveraging the effective turnover of drones and the coordination between drones and couriers.
The optimal network design features a concentrated structure, streamlining demand consolidation and reducing deadhead repositioning.
Vertiport selection is based not only on attractiveness to demands but also on geographical competitiveness.
Moving forward, maintaining enduring cost and time savings necessitates a dual focus on advancing drone technologies and reducing courier delivery expenses.
}%


\maketitle

\section{Introduction}

\subsection{Background}

Urbanization and the e-commerce boom have necessitated agile and responsive intra-city express services to efficiently handle the door-to-door delivery of items such as parcels, medicines, documents, and meals in urban areas.
The intra-city express service market size is anticipated to reach USD 26.2 billion by 2026, exhibiting a compound annual growth rate of over 9.2\% \citep{Intra2023}.
To stay ahead of the efficiency competition, online retailers partnering with courier providers have launched on-demand logistics options to offer their customers rapid courier services \citep{savelsbergh201650th}. 
However, the on-demand city logistics industry still faces challenges.
On one hand, the current city logistics system faces inefficiencies, partly attributed to traffic congestion, especially prevalent in metropolitan areas
\citep{fatnassi2015planning}. 
On the other hand, the costs for on-demand city logistics are notably higher than traditional scheduled deliveries due to the expenses linked to dedicated couriers and vehicles for door-to-door services \citep{guo2021auction}.
Breakthrough transportation technologies and practices that can simultaneously reduce costs and improve efficiency are actively pursued in the on-demand city logistics industry.

Over the past decade, drones, or unmanned aerial vehicles, have been regarded as a game changer poised to redefine city logistics \citep{savelsbergh201650th, he2022smart}. 
Since the initiation of the world’s first drone delivery project, Amazon Prime Air, in 2013 \citep{Rose2013}, drones have been deployed to provide timely delivery services in various city logistics scenarios, including package delivery \citep{baloch2020strategic}, online retailing \citep{perera2020retail}, and healthcare services \citep{enayati2023multimodal}.
Drones' benefits in enhancing delivery efficiency are evident, as they can travel at high speeds without being hindered by traffic congestion. 
Advancements in battery technologies also enable drones to deliver goods on a citywide scale.
Nevertheless, the shortcomings are equally apparent.
As per a McKinsey report \citep{Mckinsey2023}, a single-package drone delivery costs nearly \$13.50, which is not competitive with the traditional truck-routing delivery cost of \$1.90 per package.
Drones also face limitations in achieving door-to-door delivery due to airspace restrictions in densely populated areas and privacy concerns \citep{baloch2020strategic}.
Furthermore, the current ``last-meter'' delivery---unloading parcels from drones to customers---like the parachute and drop \citep{Mckinsey2023-2}, may lack precision and safety, requiring customers to be present at pick-up points to deter losses.
In response, we propose a mode in which drones are separated from the customers and confined their operations to a closed network. 
Human couriers are employed in the mode to bridge the drone-delivery system and customers.

\begin{figure}
    \centering
    \caption{Sample of Coordination Logistics Systems. (a) Meituan. (b) Hive Box}
    \includegraphics[width=\linewidth]{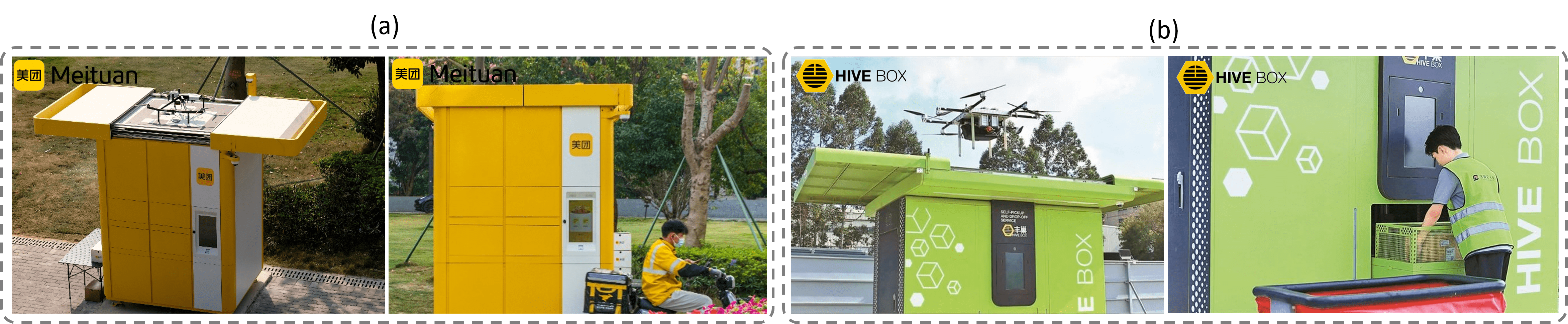}
    \label{fig.sample}
\end{figure}

We therefore propose a novel \textit{coordinated drone-courier logistics} mode that not only leverages the efficient advantages of drones but also tackles the challenges drones face in door-to-door logistics settings.
Drones transport multiple packages between dedicated take-off and landing sites, known as \textit{vertiports}, over citywide areas, while couriers manage the first-mile (collecting) and last-mile (distributing) delivery tasks in local areas.
Drones can be free from customer interaction.
Couriers can excel in collecting and distributing without incurring high costs in local areas.
This coordinated mode has been adopted by some city logistics players.
Meituan, a food delivery giant in China, for instance, has set up pickup lockers that double as vertiports for taking off and landing drones (Figure \ref{fig.sample}).
Drones do not deliver parcels directly to customers' doorsteps \citep{Meituan2023}, and the last-mile delivery is still carried out by couriers. 
Similar stories also unfolded in the DHL \citep{DHL20232} and Hive Box \citep{HiveBox2024}.
Despite the strides in practices, there is currently no practical tool available to assess the feasibility, guide network planning, and identify the drivers of sustainable development, which motivates this study.

\subsection{Contributions and Outline}

We introduce a hub location model that offers strategic guidance for long-term vertiport planning in a coordinated logistics network. 
Considering the coupling between strategic network design and tactical operations (including drone fleet sizing, dispatching, charging, and parking) in such an on-demand delivery system, we tailor a unified hub location-queueing network model that integrates these decision aspects.
Leveraging real-world data, we address three practical questions: (i) Can cost and time savings be realized through the coordination? (ii) What are the characteristics of the coordinated logistics network? (iii) What are the key drivers that enhance economical and operational efficiency for future intra-city logistics operators?

Our first contribution is the introduction of a novel coordinated logistics mode and a unified optimization model for its service network design (Section \ref{sec.model}).
From a practical standpoint, to our knowledge, this paper proposes the first coordinated drone-courier logistics mode for intra-city express services.
Existing research has only considered employing drones to perform the entire process for door-to-door delivery, which can be challenging in city logistics scenarios.
Our coordinated mode operates as a hub-and-spoke network, where drones facilitate transport between hubs (vertiports) while couriers manage customer-side pickup and delivery between vertiports and customers.
From a modeling perspective, the unified optimization model combines a hub location model (to optimize vertiport locations and capacities) with a queueing network model (to replicate the tactical system dynamics).
This unified model can account for the interplay between strategic network planning and tactical operations.
However, integrating the convex queue length function $f(\rho) = \frac{\rho}{1-\rho}$ and its negative $-f(\rho)$, where $\rho$ represents the traffic intensity and service level of the vertiport, leads to a mixed-integer nonlinear programming (MINLP) model with a feasible region that is partially convex and partially concave, motivating our algorithm design.

Our second contribution includes an exact algorithm, the adaptive discretization algorithm, for tackling the proposed MINLP model (Section \ref{sec.algorithm}).
We develop the algorithm based on the piecewise linear approximations for $f(\rho)$ and $-f(\rho)$. 
The algorithm alternates between a \textit{conservative model}---constructed from pessimistic linear approximations---to yield feasible solutions and a \textit{relaxed model}---built from optimistic linear approximations---to yield valid lower bounds.
In the iterative framework, the discretization dynamically expands, gradually relaxing the pessimistic approximations of the conservative model to improve the quality of feasible solutions, while simultaneously tightening the optimistic approximations by cutting out infeasible relaxations of the relaxed model to lift bounds.
The algorithm can provide optimality-guaranteed solutions with finite convergence.
In general, our algorithm is versatile in dealing with nonlinearity in a broad class of location-queueing models that integrate nonlinear queue length or waiting time functions. 


Our third contribution is to verify the algorithm efficacy and derive managerial insights into the coordinated drone-courier logistics (Section \ref{sec.results}).
We use an intra-city logistics order dataset provided by a prominent industry player, SF Express, in China, to study the coordinated system.
Across various network configurations and demand scenarios, our algorithm showcases scalability and robustness in delivering optimality-guaranteed solutions within reasonable computational times and fewer iterations.
The results also suggest insights into the aforementioned three questions:
(i) The proposed coordinated logistics yields a win-win outcome, achieving both cost savings and efficiency enhancements, with the benefits from coordination. 
A small fleet, empowered by the rapid turnover efficiency of drones, can mitigate the high ownership costs, while the consolidation facilitated by the coordinated network helps reduce turnaround times and improve efficiency. 
(ii) The optimal network represent a ``concentrated structure'' that enhances synergies for parcel consolidation and reduces deadheads through centralized traffic management.
The selection of vertiports is not solely based on their attractiveness to demands but also on their geographic competitive advantages.
(iii) Achieving sustainable benefits in cost and time savings depends not only on technological advancements but also on effectively managing expenses in courier-involved deliveries.

Notation tables, technical proofs, and additional analyses are available in the online appendices.

\section{Literature Review}

We review related work on the coordinated city logistics, drone delivery network planning, and hub location-queueing models and algorithms.

\subsection{Coordinated City Logistics}

In general, this research falls within the realm of smart-city operations, with a specific focus on smart-city logistics characterized by the coupled operations of multiple transport modes \citep{qi2019smart, mak2022enabling}.
With the emergence of novel mobility and technologies, new solutions for coordinated city logistics have surfaced.
\cite{qi2018shared} propose an integrated logistics system in which trucks depart from a depot to unload parcels at local terminals, while the last-mile deliveries from terminals to customers are handled by nearby ride-sharing passenger cars with available mobility.
\cite{reed2022impact} explore the possibility for cost and time savings through the use of autonomous vehicle-assisted delivery.
The truck-drone tandem last-mile delivery systems, in which truck-mounted drones replace the trucks to delivery parcels to hard-to-reach or expensive-to-reach destinations, are also being deliberated \citep{carlsson2018coordinated, dayarian2020same}.

Owing to airspace restrictions \citep{FAA2024} and security concerns \citep{Mckinsey2023-2} associated with pure drone delivery, this paper introduces a coordinated courier-drone logistics mode. 
This system features a hub-and-spoke network.
Likewise, \cite{ma2022game} suggest a coordinated metro-courier logistics mode that shares similarities in network topology with our mode. 
However, the \textit{scheduled} operations of trains stand in contrast to the \textit{on-demand} operations of drones in our system. 
Another scenario where the vehicle cannot complete the service independently is the drone delivery of automated external defibrillators to respond to out-of-hospital cardiac arrests. 
Although drones can transport these devices to patients, the emergency service can only be completed with the assistance of willing bystanders \citep{boutilier2022drone,gao2024shortening}. 
This system shares similarities in the vehicle-human coordination mode with our system. 
In contrast, it provides a \textit{depot-to-door} service instead of the \textit{door-to-door} service typical of city logistics.

\subsection{Drone Delivery Network Planning}

To enable drone delivery services, efforts have focused on developing methodologies for strategic network planning and operational drone routing (refer to \citealp{chung2020optimization} for a recent review). 
This review focuses on network planning problems from standpoints of stockholders in different scenarios.
From a delivery company standpoint, 
\cite{chauhan2019maximum} propose a model to optimize the vertiport locations for online retailers under the range constraint of drones. 
From a retailer standpoint, 
\cite{baloch2020strategic} consider the competitive facility location problem where there exists service-based competition between existing courier-based and emerging drone-based delivery services.
From an emergency response standpoint, 
\cite{lejeune2024drone} develop a network planning model aimed at maximizing survival probability in combating opioid overdoses.

The approaches in previous work, while offering valuable insights, cannot be directly applied to our problem. 
First, the intra-city express service is a \textit{door-to-door} service on a city scale, whereas the scenarios mentioned above relate to the \textit{last-mile} delivery service.
Second, the service network planning in our settings refers to a \textit{hub location problem} where hubs are interconnected by vehicles (in this case, drones) \citep{ma2022game}, while the existing research pertains to a \textit{facility location problem} without inter-facility traffic.  
Methodologically, this paper is, to our knowledge, the first attempt to address the hub (in this case, vertiport) location for the intra-city express service.

\subsection{Hub Location-Queueing Problem}

Queueing models have been well integrated into hub location problems to evaluate service metrics, such as queue length and waiting time, which directly influence the selection of hubs.
Hubs are formulated as various queueing systems, including the $M/D/1$ system \citep{marianov2003location} and $M/M/1$ system \citep{bayram2023hub}. 
In these research, the queueing networks are open. 
This means that the queueing model for each hub is separate and independent, overlooking the interconnection between hubs facilitated by vehicles.
In contrast, the closed queueing network models can formulate the on-demand transport system with traffic flowing between facilities \citep{adelman2007price}. 
We leverage the closed queueing network model proposed by \cite{pavone2015autonomous} and \cite{he2017service}.
Our paper customizes it in three ways: (1) by allowing traffic intensity $\rho$ as a continuous endogenous decision variable;
(2) by focusing on the single-allocation $P$-hub version, commonly applied in logistics settings; (3) by considering specific restrictions in drone delivery.

Differing from the typical nonlinear yet convex hub location-queueing problems \citep{elhedhli2010lagrangean, bayram2023hub}, our model deals with both nonlinearity and nonconvexity arising from $f(\rho)$ and $-f(\rho)$.
Some reformulation techniques, such as second-order cone programming (SOCP) reformulation \citep{ahmadi2022convexification,bayram2023hub} and (integer) linear reformulation \citep{lejeune2024drone}, are devised to transform the model into tractable equivalents.
These methods, however, are not viable for our nonconvex model.
\cite{he2021charging} propose a SOCP reformulation for $f(\rho)$, while it is not applicable for $-f(\rho)$ due to its nonconvexity.
Linear approximation techniques are applicable to numerically linearize the nonlinear terms \citep{de2011hybrid}.
Nevertheless, the outer approximation can only produce relaxations.
The feasible solutions must be established based on heuristics \citep{elhedhli2006service, elhedhli2010lagrangean}.

The approach proposed in \cite{elhedhli2006service} is somewhat similar to a partial component of our algorithm.
The results in Section \ref{sec.computionalResults} demonstrate the scalability of our approach over the one presented in \cite{elhedhli2006service}.
Our adaptive discretization algorithm is a generalization of the approach outlined in \cite{kai2022vertiport}.
First, our method tackles nonconvex models with a partially convex and partially concave feasible region. 
Second, we introduce a general linear approximation scheme for convex and concave terms instead of stepwise constant approximation or local piecewise approximation, which can yield tighter bounds with less computational effort. 
Third, we propose some general acceleration strategies to improve the efficiency of the algorithm.

\section{Problem Statement and Model Formulation}\label{sec.model}

\subsection{Problem Description and Assumptions}

\begin{figure}[h]
    \centering
    \caption{Traditional Dedicated Courier-based and Coordinated Drone-Courier Logistics Networks}
    \includegraphics[width=0.6\linewidth]{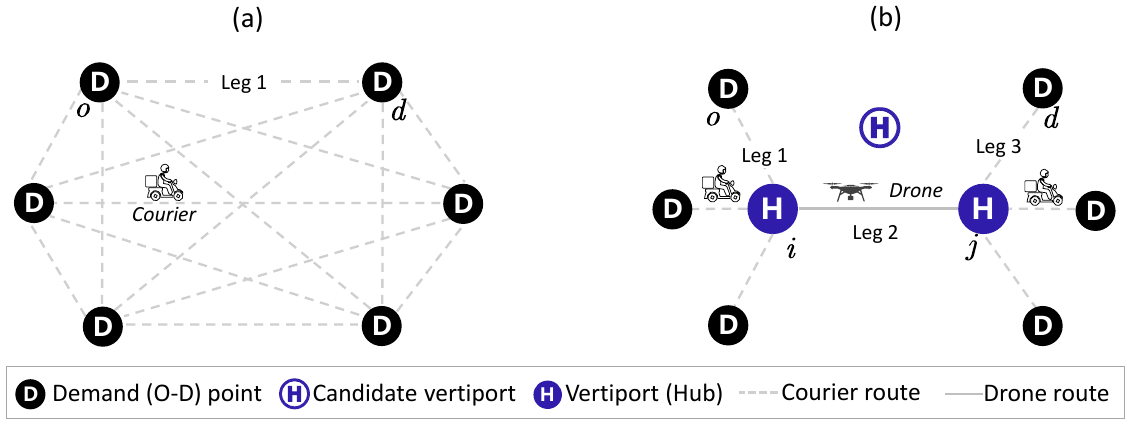}
    \label{fig.ICES} 
\end{figure}

As discussed, a door-to-door network (see Figure \ref{fig.ICES}(a)), which is adopted by the current dedicated courier service, is not viable for the application of drone delivery due to restrictions on vertiport deployment and low-altitude airspace usage.
We propose a coordinated logistics mode which features a hub-and-spoke network as illustrated in Figure \ref{fig.ICES}(b).
Operating such a coordinated drone-courier logistics system involves addressing four key issues: network planning, service level management, traffic flow management, and drone operations.  

Network planning entails determining vertiport locations and identifying service routes while considering restrictions imposed by the operator and the physical limitations of drones. 
Let $\mathcal{R}$ be the set of demand points, and let $\mathcal{N}$ denote the set of candidate vertiports. 
We consider that at most $P$ vertiports can be built from these $|\mathcal{N}|$ candidates.
In the service network, a service route $(o,i,j,d)$ comprises three itinerary legs: the collection leg from origin $o\in\mathcal{R}$ to vertiport $i\in\mathcal{N}$ by a courier, the flight leg from vertiport $i\in\mathcal{N}$ to $j\in\mathcal{N}$ by a drone, and the distribution leg from vertiport $j\in\mathcal{N}$ to destination $d\in\mathcal{R}$ by a courier.
To maintain a feasible service network, the service route $(o,i,j,d)$ can be activated if the following conditions are met: (1) the distances $l_{o,i}$ from origin $o$ to collecting vertiport $i$ and $l_{j,d}$ from distributing vertiport $j$ to destination $d$ are no longer than the service range $L^{\mathrm{S}}$ of the vertiports, and (2) the distance $l_{i,j}$ from vertiport $i$ to vertiport $j$ is within the reachable flight range $L^{\mathrm{D}}$ of the drone departing from vertiport $i$.
We define $\delta_{o,i,j,d}$ as follows, to indicate whether a service route $(o,i,j,d)$ is feasible,
\begin{eqnarray}
    \delta_{o,i,j,d} = \begin{cases}
        1, & \text{if }  l_{o,i}\leq L^{\mathrm{S}}, l_{j,d} \leq L^{\mathrm{S}}, l_{i,j}\leq L^{\mathrm{D}},\\
        0,  & \text{otherwise}.
    \end{cases}
\end{eqnarray}

We define the service level $\rho_i$ of vertiport $i$ as the proportion of demand served at vertiport $i$ to the demand attracted by vertiport $i$.
We customize a differentiated service level for each selected vertiport.
In reality, only a proportion $\xi$ of of customers require the fast delivery service, while others are unwilling to pay extra for this service \citep{savelsbergh201650th, dayarian2020same}.  
We therefore tailor the service level for each vertiport to ensure that drone-based delivery accounts for a minimum percentage $\xi$ (expected market share) of the total demand.

Traffic flow management involves regulating drone traffic across the service network.
We define drone flights with and without parcels as transit and repositioning flows, respectively. 
The transit flow is driven by the origin-destination (O-D) demand $D_{od}$. 
We assume that the arrival of parcels follows a Poisson process on each O-D pair.  
We also consider pooling parcels with the same vertiport pair. 
The operator arranges a flight service when the number of waiting parcels reaches the pooling size $Q$, with respect to the payload capacity $C$ of the drones.
The pooling process at each vertiport is depicted by a probabilistic splitting Poisson process.
This method propagates Poisson dynamics from parcel delivery demand to flight demand.
Additionally, we use standardized boxes compatible with the drones' containers.
At each vertiport, we ensure the balance between inbound and outbound flows, thereby driving the need to regulate the repositioning flow.
We assume that the repositioning processes also follow a Poisson process with an endogenous arrival rate.
This ensures the Poisson dynamics at the flow level.

Finally, we consider the operations-related problems, including fleet sizing, battery charging, and queueing issues. 
Considering the long-term nature of network planning, it is essential to consider the operations under stochasticity.
We leverage a queueing network to mirror the drone operations at the tactical level. 
To maintain the operations, a fleet size $\Gamma$ of drones traverses across the queueing network. 
Considering battery charging, we assume that drones start being charged as soon as they arrive at vertiports.
This assumption is aligned with the practice of SF Express (and many other drone delivery service providers), whose vertiports also function as charging stations. 
As a small vehicle, the battery capacity of a drone is small and the time required for a full charge is short, especially with advances in charging technology \citep{mourgelas2020autonomous}.
We therefore formulate the charging restriction in an aggregated manner.
Specifically, we ensure that the total amount of power charged at each vertiport is greater than the total amount of power required for all departing drones to reach their destinations, i.e., $b^{\mathrm{chg}}t_i^{\mathrm{chg}} \geq b^{\mathrm{fly}}t_i^{\mathrm{fly}}$, where $b^{\mathrm{chg}}$ and $t_i^{\mathrm{chg}}$ represent the charging rate and average charging time of drones at vertiport $i\in\mathcal{N}$, and $b^{\mathrm{fly}}$ and $t_i^{\mathrm{fly}}$ represent the electricity depletion rate and average flight time of drones launching at vertiport $i\in\mathcal{N}$. 
Considering the limitation of vertiports, we control the probability of the vertiport operating beyond its parking capacity by a small threshold, such as 5\%.
We assume that the parking capacity is limited and determined by the model, but the queue capacity for landing drones is assumed to be infinite.

\subsection{Model Formulation}

We develop a mixed-integer nonlinear programming model for the coordinated drone-courier logistics network planning problem. 
Corresponding to four key questions, the constraints can be classified into four categories: 
Constraints (\ref{eq.m1.P})--(\ref{eq.m1.y3}) govern the network design aspects; 
Constraints (\ref{eq.m1.rhoz})--(\ref{eq.m1.alphay}) manage the service level at each vertiport;
Constraints (\ref{eq.m1.psi})--(\ref{eq.m1.flow}) regulate traffic flows within the network;
Constraints (\ref{eq.m1.fleet})--(\ref{eq.m1.rhoCap1}) involve the operational aspects of the drones.

\subsubsection{Network Planning.} 
Let $\mathcal{H}$ denote a set of possible values for the number of parking aprons in a vertiport.
The following variables are defined to represent the service network:

\begin{tabular}{p{1cm}p{14.5cm}}
    $x_i$,       & binary variable, equals 1, if a vertiport is built at location $i\in\mathcal{N}$; 0, otherwise;   \\
    $z_{ih}$,       & binary variable, equals 1, if $h\in\mathcal{H}$ aprons are built at location $i\in\mathcal{N}$; 0, otherwise; \\
    $y_{o,i,j,d}$,  & binary variable, equals 1, if demand from $o\in\mathcal{R}$ to $d\in\mathcal{R}$ is served by drones from $i\in\mathcal{N}$ to $j\in\mathcal{N}$; 0, otherwise, \\
    $\Gamma$,  & integer variable, denoting the fleet size.  \\
\end{tabular}
These decision variables are subject to the following constraints:
\begin{subequations}
    \begin{eqnarray}
        && \sum_{i\in\mathcal{N}} x_{i} \leq P, \label{eq.m1.P}\\
        && \sum_{h\in\mathcal{H}} z_{ih} = x_i, \;\; \forall i\in\mathcal{N}, \label{eq.m1.zih1}\\ 
        && \sum_{i\in \mathcal{N}}\sum_{h\in\mathcal{H}} h z_{ih} \geq \Gamma, \label{eq.m1.Gamma}\\
        && y_{o,i,j,d}\leq x_i, \;\; \forall o,d\in\mathcal{R}, i,j\in\mathcal{N}, \label{eq.m1.yi}\\
        && y_{o,i,j,d}\leq x_j, \;\; \forall o,d\in\mathcal{R}, i,j\in\mathcal{N},\label{eq.m1.yj}\\
        && \sum_{i,j\in\mathcal{N}} y_{o,i,j,d} \leq 1, \;\; \forall o,d\in \mathcal{R}, \label{eq.m1.ysum}\\ 
        && y_{o,i,j,d} \leq \delta_{o,i,j,d}, \;\; \forall o,d\in\mathcal{R}, i,j\in\mathcal{N}. \label{eq.m1.y3}
    \end{eqnarray}
\end{subequations} 
Constraint (\ref{eq.m1.P}) ensures that no more than $P$ vertiports can be built.
Constraints (\ref{eq.m1.zih1}) mean that a capacity value $h$ should be selected if location $i\in\mathcal{N}$ is selected for vertiport construction.
Constraint (\ref{eq.m1.Gamma}) requires that there is sufficient capacity to park the drone fleet. 
Constraints (\ref{eq.m1.yi}) and (\ref{eq.m1.yj}) ensure that demand from $o$ to $d$ can be served by vertiport pair $(i,j)$ only if vertiports $i$ and $j$ are chosen for construction.
Constraints (\ref{eq.m1.ysum}) assign an O-D pair $(o,d)$ to at most one vertiport pair $(i,j)$.
Constraints (\ref{eq.m1.y3}) enforce that the selected service routes are feasible for operations.

\subsubsection{Service Level Management.}
Each constructed vertiport $i \in \mathcal{N}$ should uphold a satisfactory service level $\rho_i$. 
At this service level, the transported demand should not fall below a certain proportion (market share $\xi$) of the total demand. 
We define the service level $\rho_i$ as:

\begin{tabular}{p{1cm}p{14.5cm}}
    $\rho_i,$  & continuous variable, the service level, defined as percentage of parcels at $i\in\mathcal{N}$ that transport by drones.\\
\end{tabular}
Then the constraints for service level are given by:
\begin{subequations}
    \begin{eqnarray}
        && \rho_{i} \leq x_i, \;\; \forall i\in\mathcal{N}, \label{eq.m1.rhoz}\\
        && \sum_{o,d\in\mathcal{R}} \sum_{i,j\in\mathcal{N}} {\rho_i y_{o,i,j,d}} D_{od} \geq \xi \sum_{o,d\in\mathcal{R}} D_{od}. \label{eq.m1.demand}
    \end{eqnarray} 
\end{subequations}
Constraints (\ref{eq.m1.rhoz}) require that service levels can only be offered to the built vertiports. 
Constraint (\ref{eq.m1.demand}) ensures that the quantity of parcels shipped by drones should not be less than the expected amount. 
We linearize Constraint (\ref{eq.m1.demand}) as follows, where $\alpha_{o,i,j,d}$ is a continuous variable denoting the percentage of demand on O-D pair $(o,d)$ that is transported by drones from location $i$ to $j$:
\begin{subequations}
    \begin{eqnarray}
        && \sum_{o,d\in\mathcal{R}} \sum_{i,j\in\mathcal{N}}  \alpha_{o,i,j,d} D_{od} \geq \xi \sum_{o,d\in\mathcal{R}} D_{od}, \label{eq.m1.demandLinear}\\
        && \rho_i - (1-y_{o,i,j,d}) \leq \alpha_{o,i,j,d} \leq \rho_i, \;\; \forall o,d\in\mathcal{R}, i,j\in\mathcal{N}, \label{eq.m1.alpharho1}\\ 
        && \alpha_{o,i,j,d} \leq y_{o,i,j,d}, \;\; \forall o,d\in\mathcal{R}, i,j\in\mathcal{N}. \label{eq.m1.alphay}
    \end{eqnarray}
\end{subequations}

\subsubsection{Traffic Flow Management.}
There are transit and repositioning flows in the network.
We define the arrival rates of these traffic flows as:

\begin{tabular}{p{1cm}p{14.5cm}}
    $\psi_{i,j},$  & continuous variable,  the rate of transit trips from $i\in\mathcal{N}$ to $j\in\mathcal{N}$ per unit time;\\
    $\varphi_{i,j},$  & continuous variable,  the rate of repositioning trips from $i\in\mathcal{N}$ to $j\in\mathcal{N}$ per unit time. 
\end{tabular}
Then, the following constraints are applied to regulate traffic flows:
\begin{subequations}
    \begin{eqnarray}
        && \psi_{ij} = \sum_{o,d\in\mathcal{R}} \frac{D_{od}}{Q} {\alpha_{o,i,j,d}}, \;\; \forall i,j\in\mathcal{N}, \label{eq.m1.psi}\\
        && \sum_{j\in\mathcal{N}} \left(\psi_{j,i}+\varphi_{j,i} \right) = \sum_{j\in\mathcal{N}}{\left(\psi_{i,j}+\varphi_{i,j} \right)}, \;\; \forall i\in\mathcal{N}. \label{eq.m1.flow}
    \end{eqnarray}
\end{subequations}
Constraints (\ref{eq.m1.psi}) state that the transit flow is jointly determined by the potential demand $D_{od}$, the pooling size $Q$, and the realized percentage of demand $\alpha_{o,i,j,d}$. 
Constraints (\ref{eq.m1.flow}) control the flow balance between inbound and outbound flows at each vertiport.

\subsubsection{Drone Operations.} 
We employ a queueing network model to derive constraints for fleet sizing, battery charging, and queue overflow.

\begin{figure}[h]
    \centering
    \caption{Closed and Open queueing Network for Drone Operations}
    \includegraphics[width=0.6\linewidth]{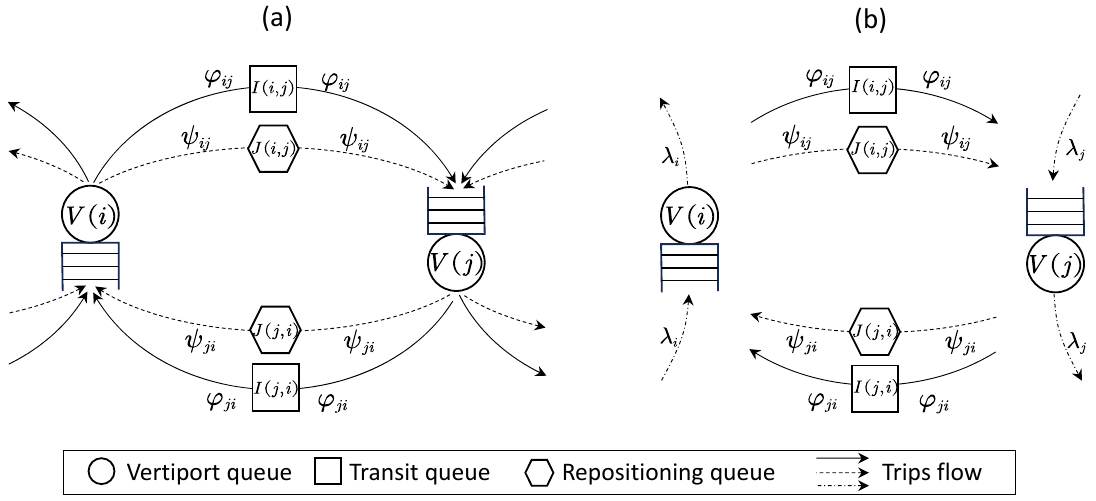}
    \label{fig.queueNetwork}
\end{figure}

\textbf{Queueing Network Model.} 
We use the closed queueing network model developed in \cite{pavone2015autonomous} and \cite{he2017service}.
Figure \ref{fig.queueNetwork}(a) details a closed queueing network.
A key feature of this model is that drones, rather than parcels, are considered as entities that flow within the queueing network.
This distinguishes it from the models proposed in the studies by \cite{elhedhli2010lagrangean} and \cite{bayram2023hub}, where parcels are treated as entities for each location queue. 
At any given time instant, a drone can be in one of three queues (states): (i) parked and charged at vertiport $i\in\mathcal{N}$, waiting for the next flight (\textit{vertiport queue} $V(i)$); (ii) traveling from vertiport $i$ to $j$ with parcels (\textit{transit queue} $I(i,j)$); and (iii) traveling from a vertiport $i$ to $j$ due to repositioning (\textit{repositioning queue} $J(i,j)$).
According to the assumptions, the arrival flows of these three queues follow Poisson processes. 
Therefore, the vertiport queues can be modeled as capacitated $M/M/1$ systems, where the service process (pooling process) also follows a Poisson process.
As for the transit and repositioning queues, they can be formulated as uncapacitated $M/G/\infty$ systems, where $G$ represents the general distribution of flight time (service time). 

We employ three metrics, i.e., the average number of drones in three queues, the average service time spent in three queues, and the probability of the vertiport queue operating over parking capacity, to formulate the constrains. 
Firstly, we use the fixed-population-mean method, which demonstrates how to approximate the number of entities in a closed queueing network (refer to Figure \ref{fig.queueNetwork}(a)) by using an open queueing network (refer to Figure \ref{fig.queueNetwork}(b)).
In the open queueing network, the number of drones in queue $V(i)$ in steady state is $f(\rho_i)=\frac{\rho_i}{1-\rho_i}$ \citep{he2017service}. 
From the Little's law, the number of drones in queues $I(i,j)$ and $J(i,j)$ can be formulated as $t_{ij}^f\psi_{ij}$ and $t_{ij}^f\varphi_{ij}$, respectively. 
Secondly, recalling the Little's law, the average service time $t_i^{\rm{chg}}$ spent at $i\in\mathcal{N}$ can be formulated as $t_i^{\rm{chg}}=f(\rho_i)/\lambda_i$, where $\lambda_{i}=\sum_{j\in\mathcal{N}}(\psi_{ij}+\varphi_{ij})$ denote the average arrival rate of drones at $i\in\mathcal{N}$. 
The average service times for the transit and repositioning queues are the flight time $t_{i,j}^f$. 
Finally, the probability that a vertiport $i\in\mathcal{N}$ is operating over its parking capacity $\mathrm{Prob}\{f(\rho_i)>\sum_{h\in \mathcal{H}}hz_{ih}\} = \rho_i^{\sum_{h\in\mathcal{H}}hz_{ih}} = \sum_{h\in\mathcal{H}}\rho_i^h z_{ih}$.
The following constraints are formulated:
\begin{subequations}
    \begin{eqnarray}
        && \sum_{i\in\mathcal{N}} {\frac{\rho_i}{1-\rho_i}} + \sum_{i,j\in\mathcal{N}} t_{i,j}^{f} \left(\psi_{i,j}+\varphi_{i,j} \right) \leq \Gamma, \label{eq.m1.fleet}\\
        && \frac{b^{\text{fly}}}{b^{\text{chg}}} \sum_{j\in\mathcal{N}}{t_{i,j}^{f}} \left(\psi_{i,j}+\varphi_{i,j}\right) - {\frac{\rho_i}{1-\rho_i}} \leq 0, \;\; \forall i\in\mathcal{N}, \label{eq.m1.charge}\\
        && \sum_{h\in\mathcal{H}_i}\gamma^{\frac{1}{h}}z_{ih} \rho_i \leq \sum_{h\in\mathcal{H}_i}{\gamma}^{\frac{1}{h+1}}z_{ih}, \;\; \forall i\in\mathcal{N}. \label{eq.m1.rhoCap1}
    \end{eqnarray}
\end{subequations}
Constraint (\ref{eq.m1.fleet}) ensures that the fleet size should not be less than the total number of drones in all vertiports $\sum_{i\in\mathcal{N}}f(\rho_i)$, transit queues $\sum_{i,j\in\mathcal{N}} t_{i,j}^{f} \psi_{i,j}$, and repositioning queues $\sum_{i,j\in\mathcal{N}} t_{i,j}^{f} \varphi_{i,j}$.
Constraints (\ref{eq.m1.charge}) is an equivalent of $b^{\mathrm{chg}}t_i^{\mathrm{chg}}\geq b^{\mathrm{fly}}t_i^{\mathrm{fly}}$.
From the queueing network model, we get $t_i^{\mathrm{chg}}=f(\rho_i)/\lambda_i$.
The average fly time $t_i^{\mathrm{fly}}$ can be obtained from the weighted average of flight times, i.e., $\frac{\sum_{j\in\mathcal{N}} {t_{i,j}^{f}} \left(\psi_{i,j}+\varphi_{i,j}\right)}{\sum_{j\in\mathcal{N}} \left(\psi_{i,j}+\varphi_{i,j}\right)} = \frac{\sum_{j\in\mathcal{N}} {t_{i,j}^{f}} \left(\psi_{i,j}+\varphi_{i,j}\right)}{\lambda_i}$. 
By multiplying both sides by $\lambda_i$, we can get Constraints (\ref{eq.m1.charge}). 
Constraints (\ref{eq.m1.rhoCap1}) are derived from $\mathrm{Prob}\{f(\rho_i)>\sum_{h\in \mathcal{H}}hz_{ih}\}\leq \gamma$, where $\gamma$ is a small threshold that caps the probability of the vertiport queues operating over their capacity.


\subsection{Full Formulation}
The objective function (\ref{eq.m1.obj}) minimizes the daily operating costs, which comprises the fleet cost, drone transport cost, and courier delivery cost (for the first- and last-mile deliveries).
$c^{\mathrm{E}}$, $p_{ij}^{\mathrm{F}}$, and $p_{oi}^{\mathrm{LT}}$ denote the unit costs for vehicles, drone flights, and courier deliveries, respectively.
\begin{eqnarray}
    \Pi = 
    \underset{\text{fleet\ cost}}{\underbrace{c^{\mathrm{E}}\Gamma }} 
    + \underset{\text{drone\ transport\ cost}}{\underbrace {\sum_{i,j\in \mathcal{N}}{p_{i,j}^{\mathrm{F}}} \left( \psi_{i,j}+\varphi_{i,j} \right) } }
    + \underset{\text{courier\ delivery\ cost}}{\underbrace{\sum_{o,d\in \mathcal{R}}{\sum_{i,j\in \mathcal{N}}{D_{od}}}\left( p_{o,i}^{\mathrm{LT}} + p_{j,d}^{\mathrm{LT}} \right)   \alpha_{o,i,j,d} }}. \label{eq.m1.obj}
\end{eqnarray}
By far, we can formulate the full model as follows:
\begin{eqnarray}
    [M^*] &\min & \Pi (\text{Equation} (\ref{eq.m1.obj})), \nonumber\\
    &\mathrm{s.t.} & \text{Network design constraints (\ref{eq.m1.P})--(\ref{eq.m1.y3})},\nonumber\\
    && \text{Service level constraints (\ref{eq.m1.rhoz}), (\ref{eq.m1.demandLinear})--(\ref{eq.m1.alphay})},\nonumber\\
    && \text{Traffic flow constraints (\ref{eq.m1.psi}) and (\ref{eq.m1.flow})},\nonumber\\
    && \text{Drone operations constraints (\ref{eq.m1.fleet})--(\ref{eq.m1.rhoCap1})},\nonumber\\
    && \bm{\alpha}, \bm{\rho}, \bm{\psi}, \bm{\varphi} \text{ non-negative};
    \Gamma \text{ integer}; \bm{x}, \bm{z}, \bm{y} \text{ binary}.\nonumber
\end{eqnarray}

\section{Solution Approach}\label{sec.algorithm}

The $[M^*]$ is a mixed-integer nonlinear program, compromising a tractable mixed-integer linear program (MILP) and nonlinear Constraints (\ref{eq.m1.fleet}) and (\ref{eq.m1.charge}).
It can be formulated in a general form:
\begin{eqnarray}\label{eq.general.model}
    &\underset{\bm{\rho},\bm{\nu}}{\min} & \bm{c}^{\top} \bm{\nu}\nonumber\\
    &\mathrm{s.t.}& \bm{F}\bm{\rho} + \bm{H}\bm{\nu} \leq \bm{d},\nonumber\\
    && \bm{f}(\bm{\rho})^{\top} \mathbf{1} + \bm{g}^{\top}\bm{\nu} \leq 0,\\
    && \bm{D}\bm{\nu}-\bm{f}(\bm{\rho}) \leq \mathbf{0}\nonumber,
\end{eqnarray}
where $\bm{f}(\bm{\rho})=\frac{\bm{\rho}}{\bm{1}-\bm{\rho}}$. $\bm{c}$ is the non-negative cost parameters corresponding to the variable vector $\bm{\nu}$ which contains all variables except $\bm{\rho}$. 
$\bm{F}$, $\bm{H}$, and $\bm{D}$ are constraint matrices.
$\bm{d}$ is a constant vector.

The queue length function $f(\rho)=\frac{\rho}{1-\rho}$ is strictly convex and non-decreasing in $\rho\in[0,1)$.
This means that the feasible region defined by the linear relaxation of Constraint (\ref{eq.m1.fleet}) is convex, whereas that defined by the linear relaxation of Constraints (\ref{eq.m1.charge}) is concave.
Therefore, $[M^*]$ is a MINLP with a feasible region that is partially convex and partially concave. 
We present general principles for tackling such a model, and then tailor them to our problem.

\subsection{General Principles} 

Let $[G^*]$ be a generic nonconvex optimization model where $g(\bm{\nu})$ is the linear objective function, $h(\chi)$ is convex over $\chi\in\mathcal{X}$, $q({\chi})$ is concave over ${\chi}\in\mathcal{X}$, and $\bm{A_1}$ and $\bm{A_2}$ are constraint matrices. 
\begin{eqnarray}
    [G^*] = \min\{ g(\bm{\nu}): \bm{h}(\bm{\chi}) + \bm{A}_1\bm{\nu} \leq \bm{0}, \bm{q}(\bm{\chi}) + \bm{A}_2\bm{\nu} \leq \bm{0}, \bm{\chi}\in{\mathcal{X}}, \bm{\nu}\in{\mathcal{V}} \}.
\end{eqnarray}
The nonlinear constraint functions make this model intractable.
To leverage the advancements of mixed-integer linear optimization solvers, we linearize this nonlinear model.
Based on the discretization of ${\mathcal{X}}$, we replace the nonlinear constraints with their ``pessimistic'' piecewise linear approximations which reduce the feasible region defined by $[G^*]$, thereby generating a \textit{conservative} model $[G_C]$.
Meanwhile, we introduce ``optimistic'' piecewise linear approximations that expand the feasible region defined by $[G^*]$, resulting in a \textit{relaxed} model $[G_R]$.
Let $(\mathcal{X}^k)_{k=1,\cdots,K}$ be a partition of $\mathcal{X}$ based on the discretization.
Let $\overline{h}^k({\chi})$ and $\overline{q}^k({\chi})$ be ``pessimistic'' linear approximations of  ${h}({\chi})$ and ${q}({\chi})$ over ${\chi}\in\mathcal{X}^k$ respectively. 
Correspondingly, $\underline{h}^k({\chi})$ and $\underline{q}^k({\chi})$ denote ``optimistic'' linear approximations of  ${h}({\chi})$ and ${q}(\chi)$ over ${\chi}\in\mathcal{X}^k$, respectively. 
We can define $[G_C]$ and $[G_R]$ as follows: 
\begin{eqnarray}
    &\left[ G_C \right] =\min \Bigg\{ g(\bm{\nu}) : & \sum_{k=1}^K{\overline{\bm{h}}^k}\left( \bm{\chi} \right) \mathbbm{1}\left(\bm{\chi}\in {\mathcal{X}}^k \right) + \bm{A}_1 \bm{\nu} \leq \bm{0}, \label{eq.general.mc}\\
    && \sum_{k=1}^K{\overline{\bm{q}}^k}\left( \bm{\chi} \right) \mathbbm{1}\left( \bm{\chi} \in {\mathcal{X}}^k \right) + \bm{A}_2 \bm{\nu} \leq \bm{0}, \bm{\chi}\in \bigcup_{k=1}^K{\mathcal{Q}_{C}^{k}},\bm{\nu }\in \bm{\mathcal{V}}\Bigg\} \nonumber,
\end{eqnarray}
where $\mathbbm{1}\left({\chi}\in {\mathcal{X}}^k \right)=1$ if ${\chi}\in {\mathcal{X}}^k$ is true, and 0 otherwise, and 
\begin{eqnarray}
    \mathcal{Q}_{C}^{k}=
    \begin{cases}
        \mathcal{X}^k,&		\text{if\,\,}\forall {\chi}\in \mathcal{X}^k: h\left( {\chi} \right) \leq \overline{h}^k\left( {\chi} \right) ,q\left( {\chi} \right) \leq \overline{q}^k\left( {\chi} \right) ,\\
        \varnothing ,&		\text{otherwise}.\\
    \end{cases}\nonumber
\end{eqnarray}
\begin{eqnarray}
    [G_R] = &\min \Bigg\{g(\bm{\nu}): & \sum_{k=1}^K\underline{\bm{h}}^k(\bm{\chi}) \mathbbm{1}(\bm{\chi}\in{\mathcal{X}}^k) + \bm{A}_1 \bm{\nu} \leq 0,\label{eq.general.mr}\\
    && \sum_{k=1}^K\underline{\bm{q}}^k(\bm{\chi}) \mathbbm{1}(\bm{\chi}\in{\mathcal{X}}^k) + \bm{A}_2 \bm{\nu} \leq 0, \bm{\chi}\in\bigcup_{k=1}^K \mathcal{Q}_R^k, \bm{\nu} \in {\mathcal{V}} \Bigg\} \nonumber,
\end{eqnarray}
where 
\begin{eqnarray}
   \mathcal{Q}_R^k = \begin{cases}
        \mathcal{X}^k, & \text{if } \forall {\chi}\in\mathcal{X}^k: h({\chi}) \geq \underline{h}^k({\chi}), q({\chi}) \geq \underline{q}^k({\chi}),\\
        \varnothing, & \mathrm{otherwise.}
    \end{cases} \nonumber
\end{eqnarray}

\begin{figure}[h]
    \centering
    \caption{Convex Function $h(\cdot)$ and Concave Function $g(\cdot)$, with Their Piecewise Linear Approximations}
    \includegraphics[width=0.6\linewidth]{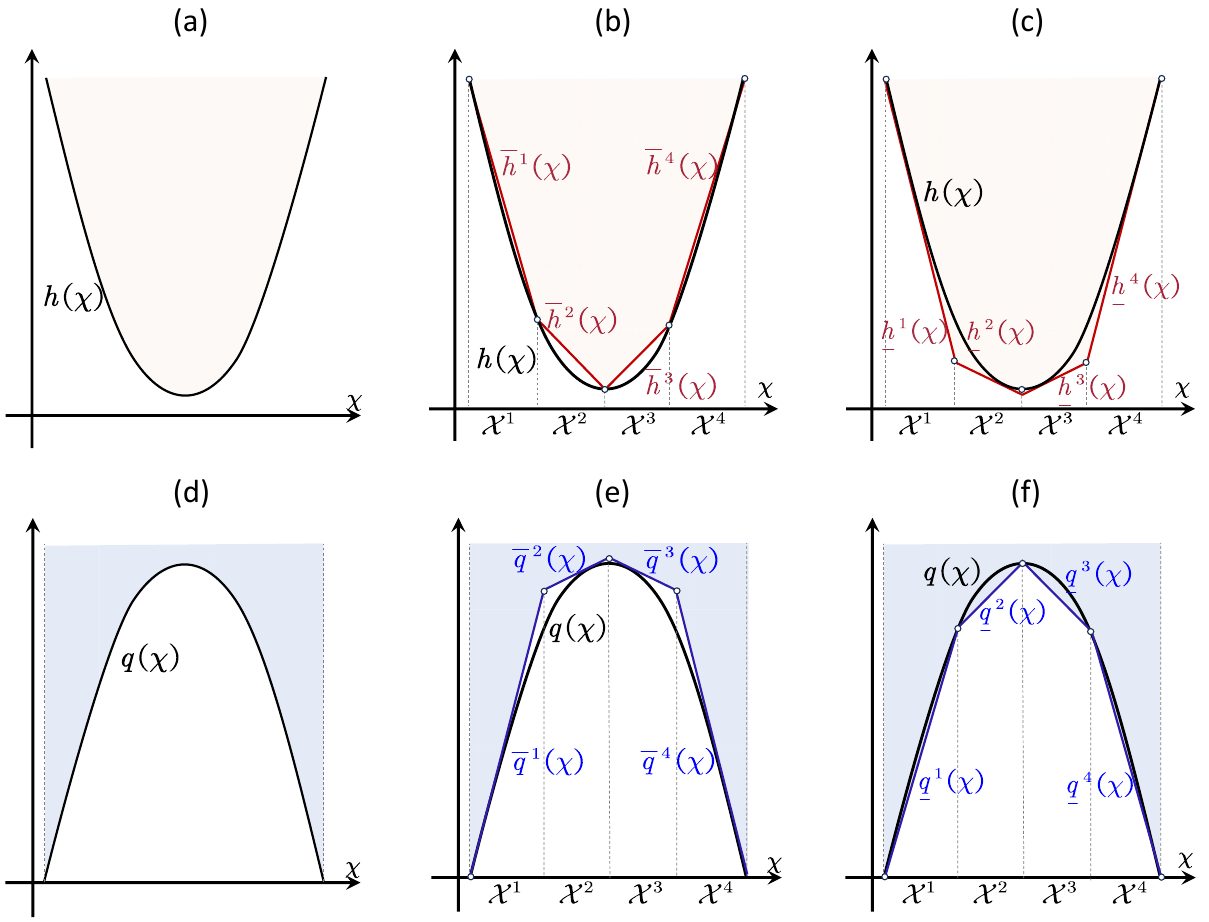}
    \label{fig.generalPrinciple}
    \begin{flushleft}
        {\footnotesize \textit{Note}. 
        (a) Feasible region defined by $h(\cdot)$. 
        (b) Pessimistic approximation of convex function $h(\cdot)$. 
        (c) Optimistic approximation of convex function $h(\cdot)$.
        (d) Feasible region defined by $g(\cdot)$.
        (e) Pessimistic approximation of concave function $g(\cdot)$.  
        (f) Optimistic approximation of concave function $g(\cdot)$.} 
    \end{flushleft}
\end{figure}

As shown in Figure \ref{fig.generalPrinciple}, $\overline{h}^k(\chi)$ and $\overline{q}^k(\chi)$ (resp., $\underline{h}^k(\chi)$ and $\underline{q}^k(\chi)$) are locally pessimistic-approximation (resp., optimistic-approximation) cutting planes. 
The resulting MILP models $[G_C]$ and $[G_R]$ yield upper and lower bounds for $[G^*]$ as stated in Proposition \ref{Pro.1}. 

\begin{myPro} \label{Pro.1}
    Let $Z[\cdot]$ denote an optimum a model.
    $[G_C]$ provides an upper bound on $Z[G^*]$, while $[G_R]$ provides a lower bound on $Z[G^*]$, i.e., $Z[G_R] \leq Z[G^*] \leq Z[G_C]$.
\end{myPro}

Proposition \ref{Pro.1} implies that an optimal solution of $[G_C]$ is a feasible solution for $[G^*]$.
The gap between $Z[G_C]$ and $Z[G_R]$ provides a measure of how well an optimal solution of $[G_C]$ can yield a feasible solution for $[G^*]$.
The upper bound $Z[G_C]$ and the lower bound $Z[G_R]$ will asymptotically approach the optimal value $Z[G^*]$ as the granularity of discretization is refined.
Given a very granular discretization, if $Z[G_C]=Z[G^*]=Z[G_R]$, the optimal solution of $[G_C]$ can be regarded as a global optimal solution.  
Nevertheless, a globally granular discretization can result in a large number of cutting planes (constraints) and auxiliary variables, which can hinder the efficiency of solving $[G_C]$ and $[G_R]$.
Our motivation is to find a way to add fewer but stronger cutting planes that not only guarantee a high-quality gap but also do not hinder efficiency.
Theorem \ref{Pro.2} states that if we take a granular discretization only locally around an optimal solution of $[G_R]$ and take a coarse discretization elsewhere, the corresponding relaxed model can guarantee a $\varepsilon$-optimal.
It also indicates the role of $[G_R]$ in providing certificates of optimality.

\begin{myTheo}\label{Pro.2}
    Assume that $h(\cdot)$ and $q(\cdot)$ are $L_h$- and $L_q$-Lipschitz continuous over $\mathcal{X}$.
    Let $L=\max\{|L_h|,|L_q|\}$.
    Define $\varepsilon > 0$.
    Let a set $\tilde{\bm{k}}$ be such that the set $\mathcal{X}^{\tilde{\bm{k}}}$ contains an optimal solution of $[G_R]$.
    If there exists a positive constant $S$ such that the diameter of $\mathcal{X}^{\tilde{\bm{k}}}$ is less than $S\varepsilon$, then there exists a feasible solution in $\mathcal{X}^{\tilde{\bm{k}}}$ within $\varepsilon$ of the optimum of $[G^*]$.
\end{myTheo}

\subsection{Conservative and Relaxed Models}

We introduce a conservative model $[M_C]$ and a relaxed model $[M_R]$ in which new piecewise linear constraints replace the nonlinear Constraints (\ref{eq.m1.fleet}) and (\ref{eq.m1.charge}).
We define a set $\bm{\varrho_i}=(\varrho_i^k)_{k\in\mathcal{K}_i}$ to discretize the domain of variable $\rho_i$, indexed in increasing values with $\varrho_i^1=0$ and $\varrho_i^{|\mathcal{K}_i|+1}=1$. 
Let $\overline{\varrho}_i^k$ represent the abscissa of the intersection point of the two tangents of $f(\rho)$ passing through points $(\varrho_i^{k-1},f(\varrho_i^{k-1}))$ and $(\varrho_i^{k},f(\varrho_i^{k}))$, i.e.,
\begin{eqnarray}
	&& \overline{\varrho}_i^k = \begin{cases}
		\frac{\varrho_i^{k-1}+\varrho_i^{k}-2\varrho_i^{k-1}\varrho_i^{k}}{2-\varrho_i^{k-1} - \varrho_i^{k}}, & k\in\mathcal{K}_i\setminus\{1,|\mathcal{K}_i|+1\}, \\
		\varrho_i^k, & \mathrm{otherwise},
	\end{cases}  \;\; \forall i\in\mathcal{N}.
\end{eqnarray}

\begin{figure}[H]
    \centering
    \caption{Illustration of Discretization}
    \includegraphics[width=0.4\linewidth]{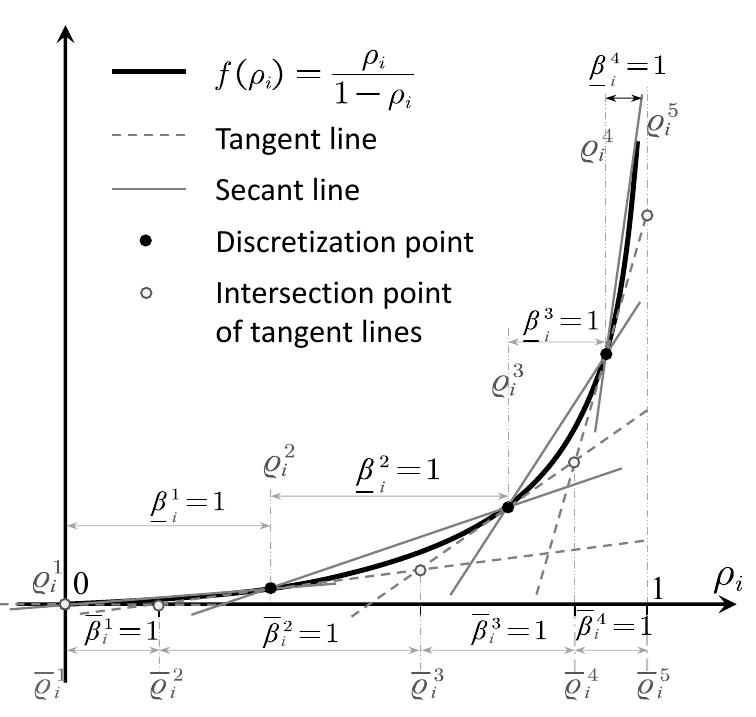}
    \label{fig.exampleOfDiscretization}
\end{figure}

Figure \ref{fig.exampleOfDiscretization} depicts the discretization scheme adopted in formulating $[M_C]$ and $[M_R]$.
Secant and tangent approximations defined by the discretization set $\bm{\varrho}$ are utilized as the linear approximations for the nonlinear terms $f(\rho)$.
The following binary variables are introduced to linearize model $[M^*]$:
\begin{tabular}{p{1cm}p{14cm}}
    $\underline{\beta}_i^k,$  & equals 1 if $\varrho_i^{k} \leq \rho \leq \varrho_i^{k+1}$; 0, otherwise;\\
    $\overline{\beta}_i^k,$  & equals 1 if $\overline{\varrho}_i^{k} \leq \rho \leq \overline{\varrho}_i^{k+1}$; 0, otherwise.\\
\end{tabular}

$[M_C]$ can be formulated as follows:
\begin{subequations}
    \begin{eqnarray}
    [M_C]
    && \overline{\Theta}_i \geq \frac{\rho_i -\varrho_i^k\varrho_i^{k+1}}{(1-\varrho_i^k)(1-\varrho_i^{k+1})}, \;\; \forall i\in\mathcal{N}, k\in\mathcal{K}_i, \label{eq.mc.Theta}\\
    && \sum_{i\in \mathcal{N}} \overline{\Theta}_i + \sum_{i,j\in \mathcal{N}}{t_{i,j}^{f}}\left( \psi_{i,j}+\varphi_{i,j} \right) \leq \Gamma, \label{eq.mc.fleet}\\
    && \frac{b^{\text{dpt}}}{b^{\text{chg}}}\sum_{j\in \mathcal{N}}{t_{i,j}^{f}}\left( \psi _{i,j}+\varphi _{i,j} \right) - \sum_{k\in\mathcal{K}_i}\frac{\overline{\pi}_i^k-(\varrho_i^k)^2 \overline{\beta}_{i}^k}{(1-\varrho_i^k)^2} \leq 0, \;\; \forall i\in\mathcal{N} \label{eq.mc.charge},\\
    && \sum_{k\in\mathcal{K}}\overline{\varrho}_i^k \overline{\beta}_{i}^k \leq \rho_i \leq \sum_{k\in\mathcal{K}}\overline{\varrho}_i^{k+1} \overline{\beta}_{i}^k, \;\; \forall i\in\mathcal{N}, \label{eq.mc.rhobeta1}\\
    && \sum_{k\in\mathcal{K}_i} \overline{\beta}_{i}^k = 1,\;\; \forall i\in\mathcal{N}, \label{eq.mc.sumbeta}\\
    && \rho_i - (1-\overline{\beta}_{i}^{k}) \leq \overline{\pi}_{i}^{k} \leq \rho_i, \;\; \overline{\pi}_{i}^{k} \leq \overline{\beta}_{i}^{k}, \;\; \forall i\in\mathcal{N}, k\in\mathcal{K}_i, \label{eq.mc.pirho1}\\
    && \overline{\bm{\Theta}}, \overline{\bm{\pi}} \;\; \text{non-negative}, \overline{\bm{\beta}} \;\; \text{binary}. \label{eq.mc.value}
\end{eqnarray}
\end{subequations}
Constraint (\ref{eq.m1.fleet}) in $[M^*]$ is replaced by Constraints (\ref{eq.mc.Theta}) and (\ref{eq.mc.fleet}) in which $\overline{\Theta}_i$ is a non-negative continuous decision variable representing an upward approximation of the number of drones at location $i\in\mathcal{N}$.
The right-hand side of Constraints (\ref{eq.mc.Theta}) is the secant defined by the points $(\varrho_i^k,f(\varrho_i^{k}))$ and $(\varrho_i^{k+1},f(\varrho_i^{k+1}))$.
Constraints (\ref{eq.m1.charge}) in $[M^*]$ are replaced by (\ref{eq.mc.charge}).
The second term on the left-hand side of Constraints (\ref{eq.mc.charge}) is the tangent determined by point ($\varrho_i^k,\frac{\varrho_i^k}{1-\varrho_i^k}$).
$\overline{\pi}_i^k$ is a non-negative continuous decision variable, which is equal to $\rho_i\overline{\beta}_i^k$.
Constraints (\ref{eq.mc.rhobeta1}) ensure that $\overline{\beta}_{i}^k$ is active only when $\overline{\varrho}_i^{k} \leq \rho_i \leq \overline{\varrho}_i^{k+1}$. 
Constraints (\ref{eq.mc.pirho1}) linearize the bilinear term $\overline{\pi}_{i}^k = \rho_i \overline{\beta}_{i}^k$.

\begin{myPro} \label{Pro.4}
    Constraints (\ref{eq.mc.Theta})--(\ref{eq.mc.value}) satisfy the general principles in (\ref{eq.general.mc}) and are pessimistic approximations to Constraints (\ref{eq.m1.fleet}) and  (\ref{eq.m1.charge}).
    Any feasible solutions of $[M_C]$ are feasible for $[M^*]$.
\end{myPro}

$[M_R]$ can be formulated as follows:
\begin{subequations}
    \begin{eqnarray}
        [M_R]  
        && \underline{\Theta}_i \geq \frac{\rho_i-(\varrho_i^k)^2}{(1-\varrho_i^k)^2},\;\; \forall i\in\mathcal{N}, k\in\mathcal{K}_i, \label{eq.mr.Theta}\\
        && \sum_{i\in \mathcal{N}} \underline{\Theta}_i  + \sum_{i,j\in \mathcal{N}}{t_{i,j}^{f}} \left(\psi_{i,j}+\varphi_{i,j}\right) \leq \Gamma, \label{eq.mr.fleet}\\
        && \frac{b^{\text{dpt}}}{b^{\text{chg}}}\sum_{j\in \mathcal{N}}{t_{i,j}^{f}} \left(\psi_{i,j}+\varphi_{i,j}\right) - \sum_{k\in\mathcal{K}_i} \frac{\underline{\pi}_i^k -\varrho_i^k\varrho_i^{k+1}\underline{\beta}_{i}^k}{(1-\varrho_i^k)(1-\varrho_i^{k+1})} \leq 0, \;\; \forall i\in\mathcal{N} \label{eq.mr.charge},\\
        && \sum_{k\in\mathcal{K}_i}\varrho_{i}^k \underline{\beta}_{i}^k \leq \rho_i \leq \sum_{k\in\mathcal{K}_i}\varrho_{i}^{k+1} \underline{\beta}_{i}^k, \;\; \forall i\in\mathcal{N}, \label{eq.mr.rhobeta1}\\
        && \sum_{k\in\mathcal{K}_i} \underline{\beta}_{i}^k = 1,\;\; \forall i\in\mathcal{N}, \label{eq.mr.sumbeta}\\
        && \rho_i - (1-\underline{\beta}_{i}^{k}) \leq \underline{\pi}_{i}^{k} \leq \rho_i, \;\; \underline{\pi}_{i}^{k} \leq \underline{\beta}_{i}^{k}, \;\; \forall i\in\mathcal{N}, k\in\mathcal{K}_i, \label{eq.mr.pirho0}\\
        && \underline{\bm{\Theta}}, \underline{\bm{\pi}} \;\; \text{non-negative}, \underline{\bm{\beta}} \;\; \text{binary}. \label{eq.mr.value}
    \end{eqnarray}
\end{subequations}
The right-hand side of Constraints (\ref{eq.mr.Theta}) is the tangent through point $(\varrho_i^k,f(\varrho_i^k))$.
$\underline{\Theta_i}$ is a non-negative continuous decision variable denoting a downward approximation of the number of drones at location $i\in\mathcal{N}$. 
The second term on the left-hand side of Constraints (\ref{eq.mr.charge}) is the secant determined by the points $(\varrho_i^k,f(\varrho_i^{k}))$ and $(\varrho_i^{k+1},f(\varrho_i^{k+1}))$.
$\underline{\pi}_i^k$ is a non-negative continuous decision variable, which is equal to $\rho_i\underline{\beta}_i^k$.
Constraints (\ref{eq.mr.rhobeta1}) define the activated discretization domain of $\rho_i$. 
Constraints (\ref{eq.mr.pirho0}) are linear equivalents of $\underline{\pi}_{i}^k = \rho_i \underline{\beta}_{i}^k$.

\begin{myPro} \label{Pro.5}
    Constraints (\ref{eq.mr.Theta})--(\ref{eq.mr.value}) are optimistic approximations to Constraints (\ref{eq.m1.fleet}) and  (\ref{eq.m1.charge}), satisfying the general principles in (\ref{eq.general.mr}).
    Any feasible solutions of $[M^*]$ are feasible for $[M_R]$.
\end{myPro}

Propositions \ref{Pro.4} and \ref{Pro.5} suggest that the optima of $[M_C]$ and $[M_R]$ offer lower and upper bounds on $Z[M^*]$.
Proposition \ref{Pro.6} additionally indicates that these bounds become tighter as the discretization becomes finer.
This motivates the adaptive discretization algorithm, which iteratively refines the discretization set and tightens the bounds.
\begin{myPro} \label{Pro.6}
    Let $[M_C]$ and $[M_C']$ (resp., $[M_R]$ and $[M_R']$) refer to the conservative (reps., relaxed) model with two discretizations $\mathcal{K}_i'$ and  $\mathcal{K}_i$ such that $\{\varrho_i^{k'}\}_{k'\in\mathcal{K}_i'} \subseteq \{\varrho_i^{k}\}_{k\in\mathcal{K}_i}$ for all $i\in\mathcal{N}$. 
    It holds that $Z[M_C'] \geq Z[M_C] \geq Z[M^*] \geq Z[M_R] \geq Z[M_R']$.
\end{myPro}

\begin{myremark}
    In $[M_C]$ and $[M_R]$, we approximate both $f(\rho)$ and $-f(\rho)$ using linear approximations. 
    Considering the convexity of $f(\rho)$, \cite{he2021charging} propose a SOCP reformulation for it. 
    While this reformulation provides an equivalent representation for $f(\rho)$ without approximation, its computational burden is considerably higher than our linear approximation approach.
    Moreover, the SOCP reformulation is not suitable for addressing $-f(\rho)$.
    Preliminary tests conducted on this SOCP reformulation solely for $f(\rho)$ suggest that computational times are unreasonable even for small-scale cases (refer to Online Appendix \ref{sec.socp} for more details).
\end{myremark}

\subsection{Adaptive Discretization Approach}

\subsubsection{Overview.}
Our adaptive discretization approach is an iterative algorithm that dynamically and ``wisely'' expands the set $\bm{\varrho}$ of discretization points around the current incumbent. 
The algorithm alternates between solving $[M_C]$ and $[M_R]$ to gradually tighten the upper bound and lift the lower bound. 
Specifically, after solving $[M_C]$ and obtaining its optimal solution $\bm{\rho}^C$, we apply an adaptive discretization procedure to expand set $\bm{\varrho}$ around $\bm{\rho}^C$. 
The added discretization points yield new cuts (Constraints (\ref{eq.mr.Theta})--(\ref{eq.mr.value})), which reduce the feasible region of $[M_R]$ in the next iteration, thereby improving the lower bound.
Subsequently, a new optimal solution $\bm{\rho}^R$ obtained from solving $[M_R]$ is also inputted into the adaptive discretization procedure to update set $\bm{\varrho}$. 
The expanded $\bm{\varrho}$ leads to an expanded feasible region of $[M_C]$ and results in a tighter upper bound.
The algorithm terminates until the optimal gap is less than a threshold $\epsilon$.


Theorem \ref{theorem.2} highlights the finite convergence of the algorithm.
Our algorithm alternates between $[M_C]$ and $[M_R]$, leveraging the potential of $[M_C]$ in finding ``promising'' feasible solution and of $[M_R]$ in providing optimality guarantees. 
The algorithm converges as the discretization is refined.

\begin{myTheo} \label{theorem.2}
    Define $\varepsilon > 0$. Let vector $\bm{\nu}$ contains all variables except $\bm{\rho}$. Let $(\bm{\rho}^R, \bm{\nu}^R)$ be an optimal solution of $[M_R]$.
    For each $i\in\mathcal{N}$, let $k_i$ be such that $\varrho_{i}^{k_i} \leq \rho_i \leq \varrho_i^{k_i+1}$.
    There exists $S>0$ such that, if $\varrho_{i}^{k_i+1} - \varrho_i^{k_i} \leq S\varepsilon$ $\forall i\in\mathcal{N}$, we can construct a feasible solution $(\bm{\rho}^R, \tilde{\bm{\nu}})$ within $\varepsilon$ of the global optimum. 
\end{myTheo}

Theorem \ref{theorem.2} also implies the potential for ``wise'' expansion of the discretization set $\bm{\varrho}$.
It shows that the solution of $[M_R]$ can serve as a certificate of optimality with only a ``locally'' granular discretization.
We therefore propose an adaptive discretization procedure that not only guarantees the near-optimality---by local granular discretization around $[M_R]$ solutions---but also enhances the tractability---by coarse discretization elsewhere. 

\subsubsection{Adaptive Discretization Procedure.} \label{sec.procedure}
Let $\bm{\rho}$ denote the optimal solution of $[M_C]$ or $[M_R]$ in a given iteration.
We add new discretization points ``around''  $\bm{\rho}$. 
For a location $i\in\mathcal{N}$, let $k$ be such that $\varrho_i^k \leq \rho_i \leq \varrho_i^{k+1}$.
Generally, inserting $\rho_i$ directly into set $\bm{\varrho}_i$ is a straightforward approach to refine the discretization.
This method resembles the one proposed by \cite{elhedhli2006service}.
After adding this point, the feasible region of $[M_C]$ is expanded, while that of $[M_R]$ is reduced. 

In addition to adding $\rho_i$, we insert midpoints $\frac{\varrho_i^k+\rho_i}{2}$ and $\frac{\rho_i+\varrho_i^{k+1}}{2}$ into $\bm{\varrho}_i$.
As per Proposition \ref{Pro.6}, these additional midpoints enable the algorithm to tighten the upper bound obtained from $[M_C]$ and raise the lower bound obtained from $[M_R]$, thereby speeding up convergence.
However, if the points to be added are very close to $\varrho_i^k$ or $ \varrho_i^{k+1}$, the deduced approximation may not expand (reps., reduce) the current feasible region of $[M_C]$ (reps., $[M_R]$) efficiently. 
Therefore, we only add points that are more than a threshold $\Delta_{\min}$ away from the existing points.


\subsubsection{Acceleration Strategy.} Two acceleration strategies are implemented.

\textbf{Warm Start.} According to Proposition \ref{Pro.6}, any feasible solutions of $[M_C']$ are feasible for $[M_C]$ and $[M_R]$.
This indicates that an optimal solution of $[M_C]$ in the $l$-th iteration is feasible for $[M_R]$ in the $(l+1)$-th iteration and $[M_C]$ in the $(l+2)$-th iteration.
Therefore, we warm start the subsequent $[M_R]$ and $[M_C]$ by using an optimal solution obtained from $[M_C]$ in the last iteration.

\textbf{Neighborhood Search.} After solving $[M_C]$ in a given iteration, we get an optimal solution of location variable $\bm{x}^C$ and service level variable $\bm{\rho}^C$.
We can fix the location decision in $[M_C]$ as $\bm{x}=\bm{x}^C$ and use a globally granular discretization $\bm{\varrho}^N$, for instance, a discretization with a unit $\Delta\leq 0.05$.
This leads to a new model $[M_C^N]$. 
Solving $[M_C^N]$ can further explore the more ``promising'' operational decisions in given location decisions. 
Additionally, the incumbent solution $\bm{\rho}^C$ is inserted into the discretization set $\bm{\varrho}^N$.
This ensures that the optimal solution obtained from $[M_C]$ is also feasible for $[M_C^N]$, and that the $[M_C^N]$ optimum is no worse than the $[M_C]$ optimum. 
Note that the incumbent solution of $[M_C]$ is also applicable to the warm start of $[M_C^N]$.

Algorithm \ref{alg.exact} in Online Appendix \ref{sec.alg.exact} summarizes and details the proposed adaptive discretization approach. 
In addition, we devise an algorithm to obtain high-quality solutions of large-scale cases (see Online Appendix \ref{sec.heu.appendix}). 

\subsection{Contribution to the Literature of Location-Queueing Models} 

The adaptive discretization algorithm presents a versatile solution approach for a diverse range of location-queueing models, where queueing-related formulations, such as queue length and waiting time, commonly introduce nonlinearity and complexity \citep{elhedhli2006service,bayram2023hub}. 
In our model $[M^*]$, we employ the adaptive discretization algorithm to address the nonlinearity and nonconvexity arising from the queue length function $f(\rho)$ and its negative $-f(\rho)$.
More generally, our algorithm is capable of handling the nonlinearity inherent in formulating queue length or waiting time across various queueing models such as the $M/D/1$ model \citep{marianov2003location}, $M/M/1$ model \citep{bayram2023hub}, and $M/G/K$ model \citep{lejeune2024drone}. 
Methodologically, our algorithm offers a globally optimal solution approach for location-queueing models, with bound guarantees and finite convergence.
Computationally, compared to prior SOCP-based reformulations \citep{he2021charging,bayram2023hub} and a similar piecewise linear approximation-based approach \citep{elhedhli2006service}, the computational results in Section \ref{sec.computionalResults} and Online Appendix \ref{sec.socp} will demonstrate the superior performance of our algorithm.

\section{Computational Results and Managerial Implications} \label{sec.results}


We collaborate with SF Express, a prominent logistics company in China, to explore coordinated drone-courier logistics.
SF Express provides a detailed dataset of intra-city express orders from Shenzhen, a megacity and logistics hub, recording origin and destination details, loading and unloading times, and parcel volumetric weight.
We aggregate this data to compute the arrival rate for each O-D pair and validate its adherence to a Poisson arrival process.
Experimental instances are generated using this dataset and the parameters of the drones provided by SF Express.
Further details on instance generation can be found in Online Appendix \ref{ec.expSetup}.

\subsection{Algorithm Efficacy} \label{sec.computionalResults}


\subsubsection{Benchmarks and Computational Setup.}

We denote a static run of Model $A$ as $Sta(A)$ and an adaptive discretization of Model $A$ as $Ada(A)$.
When running Models $A$ and $B$ alternatively, we refer to it as $(A,B)$. 
When running Models $A$ and $B$ independently, we refer to it as $\{A,B\}$.
When running Model $A$ along with its neighborhood search procedure, we refer to it as $A+A^N$.
Applying this nomenclature, we run following algorithms:
\begin{itemize}[leftmargin=*]
    \item $Sta\{M_C,M_R\}$: Run $[M_C]$ and $[M_R]$ statically and independently with a given discretization;
    \item $Ada\{M_C,M_R\}$: Run $[M_C]$ and $[M_R]$ independently within the adaptive discretization framework;
    \item $Ada(M_C,M_R)$: Run $[M_C]$ and $[M_R]$ iteratively within the adaptive discretization framework;
    \item $Ada(M_C+M_C^N,M_R)$: Run $Ada(M_C,M_R)$ along with the neighborhood search (Algorithm \ref{alg.exact}).
\end{itemize}

Let a string $\#\rm{OD}$-$|\mathcal{N}|$-$P$ denote the settings (the numbers of O-D pairs, candidate vertiports, and selected vertiports) of an instance.
We generate 15 groups of instances, with each group containing five randomly generated instances with identical network configurations.
The time limit for each approach was set to 7,200 seconds, with a single run for each iteration set to 3,600 seconds.
The optimality gap threshold was set to 1\%.  
All algorithms were coded in C++ calling CPLEX 22.1 on a PC equipped with an Intel Core i7 2.10 GHz processor and 32 GB of RAM. 

\subsubsection{Values of Algorithmic Design.} 
The proposed algorithm offers benefits in terms of the finite convergence with optimality guarantees and computational efficiency.

\begin{figure}[h]
    \centering
    \caption{Convergence of the the Adaptive Discretization Algorithms:  $Ada(M_C,M_R)$ and $Ada(M_C+M_C^N,M_R)$}
    \includegraphics[width=0.8\linewidth]{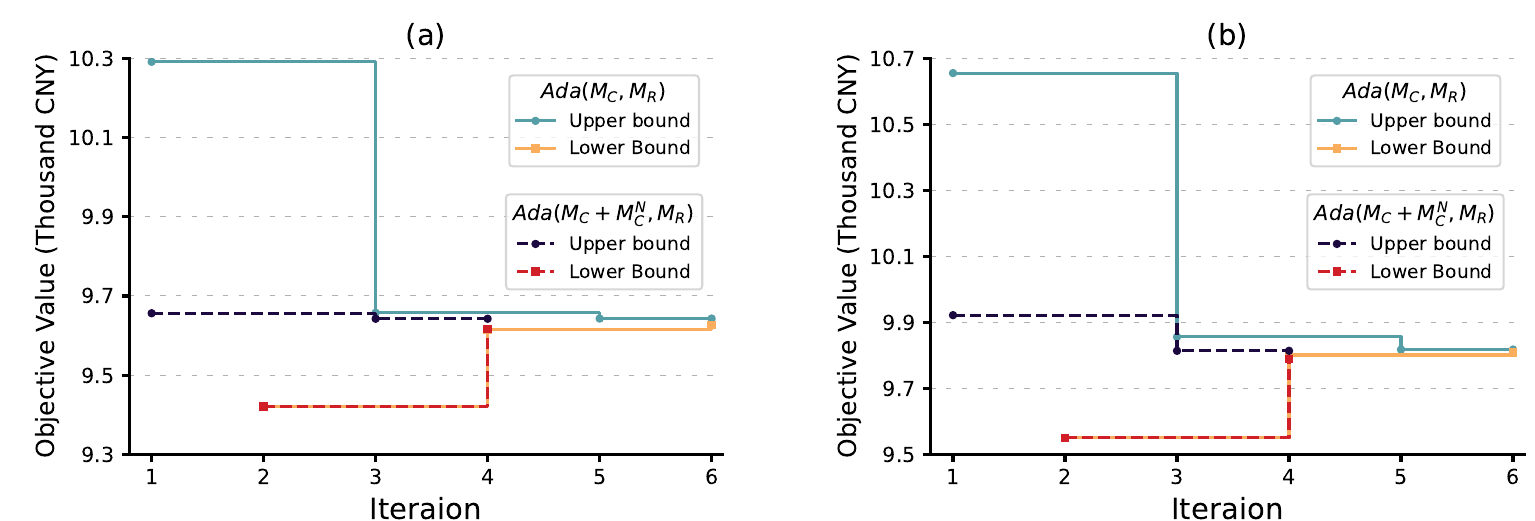}
    \label{fig.iter}
    \begin{flushleft}
        {\footnotesize \textit{Note}. (a) An instance in group 200-20-10; (b) An instance in group 200-20-12.}
    \end{flushleft}
\end{figure}

\textbf{Convergence with Optimality Guarantee}.
In all instances, our generic adaptive discretization algorithm $Ada(M_C,M_R)$ and full algorithm $Ada(M_C+M_C^N,M_R)$ converge to optimality gaps of 1\% with a finite number of iterations. 
Figure \ref{fig.iter} demonstrates the convergence in two instances.
Both algorithms converge to the threshold gap within no more than 6 iterations, showcasing their capability to achieve the global optimality for the nonconvex model $[M^*]$. 

By contrast, as shown in Table \ref{tab.dynamicDiscretization}, the static benchmark algorithm $Sta\{M_C,M_R\}$ using coarse discretization units, such as $\Delta=0.2$ or $0.1$, may not converge to the required gaps, with an average gap of 2.49\% or 1.09\%.
Although the static algorithm with a granular discretization, such as $\Delta=0.05$, can guarantee threshold gaps, the increased complexity due to the expansion of constraints and binary variables may limit its scalability. 
In addition, $Ada\{M_C,M_R\}$ does not converge to the required optimality gap, as Model $[M_C]$ cannot avoid local optima by itself without the ``global exploration'' enabled by $[M_R]$, highlighting the value of $[M_R]$ in guaranteeing convergence.

\begin{table}[h]
  \centering
  \tiny
  \caption{Comparison of Average Computational Results of Our Algorithm and Benckmarks}
  \begin{threeparttable}
    \begin{tabular}{lrrrrrrrrrrrrrrrrrr} 
    \toprule
    & \multicolumn{9}{c}{$Sta\{M_C,M_R\}$} & \multicolumn{3}{c}{$Ada\{M_C,M_R\}$}      & \multicolumn{3}{c}{$Ada(M_C,M_R)$} & \multicolumn{3}{c}{$Ada(M_C+M_C^N,M_R)$}   \\
    \cmidrule(r){2-10}    
    & \multicolumn{3}{c}{$\Delta=0.2$}   & \multicolumn{3}{c}{$\Delta=0.1$}        & \multicolumn{3}{c}{$\Delta=0.05$}      &    \\
    \cmidrule(r){2-4}    \cmidrule(r){5-7}   \cmidrule(r){8-10}  \cmidrule(r){11-13}   \cmidrule(r){14-16} \cmidrule(r){17-19}      \#OD-$|\mathcal{N}|$-$P$     & \multicolumn{1}{r}{UB} & \multicolumn{1}{r}{Gap} & \multicolumn{1}{r}{CPU}  & \multicolumn{1}{r}{UB} & \multicolumn{1}{r}{Gap} & \multicolumn{1}{r}{CPU}  & \multicolumn{1}{r}{UB} & \multicolumn{1}{r}{Gap} & \multicolumn{1}{r}{CPU}  & \multicolumn{1}{r}{UB} & \multicolumn{1}{r}{Gap} & \multicolumn{1}{r}{CPU}  & \multicolumn{1}{r}{UB} & \multicolumn{1}{r}{Gap} & \multicolumn{1}{r}{CPU} &\multicolumn{1}{r}{UB}  & \multicolumn{1}{r}{Gap} & \multicolumn{1}{r}{CPU} \\
    \midrule
    200-20-10 & 9.67   & 2.93  & 4     & 9.57    & 1.15  & 3     & 9.53    & 0.22  & 3     & 9.70  & 2.11  & 8     & 9.53  & 0.33  & 3     & 9.53  & 0.20  & {3} \\
    200-20-12 & 9.49   & 2.52  & 12    & 9.41    & 1.03  & 9     & 9.37    & 0.20  & 12    & 9.70  & 3.94  & 15    & 9.40  & 0.74  & 7     & 9.39  & 0.63  & {7} \\
    200-20-14 & 9.41   & 2.34  & 14    & 9.34    & 0.99  & 14    & 9.30    & 0.17  & 17    & 9.71  & 4.71  & 20    & 9.31  & 0.64  & 12    & 9.32  & 0.72  & {11} \\
    \cmidrule(r){1-1} \cmidrule(r){2-4}    \cmidrule(r){5-7}   \cmidrule(r){8-10}  \cmidrule(r){11-13}   \cmidrule(r){14-16} \cmidrule(r){17-19}
    300-30-10 & 12.54 & 2.68  & 11    & 12.42  & 1.03  & 9     & 12.37  & 0.16  & 8     & 12.38 & 0.75  & 15    & 12.39 & 0.54  & 8     & 12.37 & 0.38  & {8} \\
    300-30-12 & 12.15 & 2.32  & 46    & 12.07 & 0.92  & 59    & 12.03  & 0.17  & 65    & 12.33 & 3.28  & 58    & 12.04 & 0.87  & 38    & 12.03 & 0.73  & {29} \\
    300-30-14 & 11.98  & 2.36  & 81   & 11.90  & 0.92  & 129  & 11.86  & 0.16  & 1,072   & 12.38 & 5.27  & 132   & 11.86 & 0.65  & 108   & 11.86 & 0.83  & {69} \\
    \cmidrule(r){1-1} \cmidrule(r){2-4}    \cmidrule(r){5-7}   \cmidrule(r){8-10}  \cmidrule(r){11-13}   \cmidrule(r){14-16} \cmidrule(r){17-19}
    400-40-10 & 15.02 & 2.66  & 37    & 14.88 & 0.92  & 26    & 14.83 & 0.20  & 40     & 14.83 & 0.67  & 58    & 14.84 & 0.41  & {18}  & 14.82 & 0.25  & {19} \\
    400-40-12 & 14.56 & 2.55  & 342   & 14.43 & 0.9   & 387   & 14.38 & 0.15  & 1,074  & 14.73 & 3.50  & 323   & 14.38 & 0.71  & 238   & 14.38 & 0.76  & {177} \\
    400-40-14 & 14.26 & 2.17  & 736   & 14.16 & 0.8   & 2,090  & 14.12 & 0.52  & 6,396  & 14.72 & 5.31  & 663   & 14.12 & 0.93  & 1,126  & 14.11 & 0.88  & {369} \\
    \cmidrule(r){1-1} \cmidrule(r){2-4}    \cmidrule(r){5-7}   \cmidrule(r){8-10}  \cmidrule(r){11-13}   \cmidrule(r){14-16} \cmidrule(r){17-19}
    500-50-10 & 16.87 & 2.23  & 98    & 16.75 & 0.94  & 108   & 16.69 & 0.16  & 197   & 16.70 & 0.79  & 109   & 16.70 & 0.88  & 42    & 16.69 & 0.90  & {26} \\
    500-50-12 & 16.32 & 2.35  & 1,119  & 16.20 & 0.90  & 1,128  & 16.15 & 0.26  & 2,784  & 16.45 & 3.02  & 623   & 16.16 & 0.87  & 510   & 16.14 & 0.76  & {367} \\
    500-50-14 & 15.98 & 2.81  & 6,133  & 15.89 & 2.30  & 7,200  & 15.84 & 2.13  & 7,200 & 16.50 & 6.13  & 5,760  & 15.84 & 0.96  & 5,060  & 15.83 & 0.95  & {1,980} \\
    \cmidrule(r){1-1}  \cmidrule(r){2-4}    \cmidrule(r){5-7}   \cmidrule(r){8-10}  \cmidrule(r){11-13}   \cmidrule(r){14-16} \cmidrule(r){17-19}
    600-60-10 & 19.02 & 2.22  & 46    & 18.88 & 0.95  & 53    & 18.82 & 0.17  & 62    & 18.83 & 0.53  & 126   & 18.81 & 0.71  & 41    & 18.81 & 0.67  & {26} \\
    600-60-12 & 18.24 & 2.37  & 1719  & 18.10 & 0.91  & 1,826  & 18.04 & 0.13  & 1,913  & 18.06 & 1.30  & 1,044  & 18.07 & 0.82  & 2,010  & 18.04 & 0.68  & {503} \\
    600-60-14 & 17.63 & 2.53  & 6,541  & 17.52 & 1.75  & 7,200  & 17.47 & 1.38  & 7,200  & 17.59 & 2.70  & 6,602  & 17.47 & 0.94  & 5,624  & 17.46 & 0.88  & {2,898} \\
    \cmidrule(r){1-1}  \cmidrule(r){2-4}    \cmidrule(r){5-7}   \cmidrule(r){8-10}  \cmidrule(r){11-13}   \cmidrule(r){14-16} \cmidrule(r){17-19}
    Total & 14.21 & 2.49  & 1,129  & 14,10 & 1.09  & 1,349  & 14.05 & 0.41  & 1,870  & 14.31 & 2.93  & 1,037  & 14.06 & 0.73  & 990 & 14.05 & 0.68  & {433} \\
    \bottomrule
    \end{tabular}
    \begin{tablenotes}
        \item \textit{Note}. This table includes average values of the upper bound (UB, in thousand CNY), optimality gap (Gap, in percentage terms), and computational time (CPU, in CUP seconds) for each instance group.
    \end{tablenotes}
  \end{threeparttable}
  \label{tab.dynamicDiscretization}%
\end{table}%

\textbf{Computational Efficiency}.
In the comparison with $Sta\{M_C,M_R\}$, Table \ref{tab.dynamicDiscretization} shows that our proposed algorithm $Ada(M_C+M_C^N,M_R)$ saves, on average, over 60\% CPU time compared to the static benchmark with $\Delta=0.2$, with an almost 80\% reduction when compared to the static benchmark with $\Delta=0.05$.
These results highlight the value of adaptive discretization in improving efficiency and verify the possibility of ``locally'' granular discretization as stated in Theorem \ref{theorem.2}.

In the comparison with $Ada\{M_C,M_R\}$, our algorithm also shows an average computational time reduction of nearly 60\%. 
$Ada\{M_C,M_R\}$ can be regarded as the method proposed in \cite{elhedhli2006service}.
Table \ref{tab.dynamicDiscretization} shows that $Ada\{M_C,M_R\}$ yields upper bounds that are, on average, 1.7\% inferior to those obtained by our algorithm, indicating a deviation towards local optima.
In contrast, the lower bounds obtained from $Ada\{M_C,M_R\}$ are close to those obtained from $Ada(M_C+M_C^N,M_R)$, with only a marginal deviation of 0.2\% lower.
These results further confirm Theorem \ref{theorem.2}, highlighting the role of the relaxed model $[M_R]$ in providing certificates of optimality and guiding the conservative model $[M_C]$ to avoid local optima.

In the comparison with $Ada(M_C,M_R)$, the results indicate that the neighborhood search can significantly reduce the computational time by more than half on average, and this reduction is particularly noticeable in instances with large values of $P$.
By incorporating the neighborhood search phase, the upper bounds can be substantially reduced by leveraging the incumbent.
In addition, the number of iterations required to achieve the desired optimality gap can be reduced, demonstrating the role of neighborhood search in reducing the computational burden.


Our algorithm also showcases the scalability. 
In large networks (with large numbers of O-D pairs $\#$OD and candidate vertiports $|\mathcal{N}|$) or relaxed resource requirements (with large $P$), other algorithms may exceed a 2-hour CPU time limit and/or fail to provide optimality-guaranteed solutions.
The static benchmark algorithm adds numerous valid yet unnecessary cuts, needlessly increasing computational load, particularly for Model $[M_R]$.
In contrast, our algorithm can achieve a fine-grained discretization around the $[M_R]$ solutions while maintaining a coarse discretization elsewhere, thus facilitating scalability.

In summary, our algorithm $Ada(M_C+M_C^N,M_R)$ can provide optimality-guaranteed solutions for Model $[M^*]$ with shorter CPU times and fewer iterations, demonstrating the algorithm efficacy.

\subsection{Managerial Insights} \label{sec.insights}

We discuss three practical questions regarding the benefits of the coordination, the features of coordinated networks, and the drivers of potential benefits.

\subsubsection{Benefits of Coordination.} \label{sec.networkBenefit}

We measure the value of our coordinated drone-courier logistics system in terms of cost-effectiveness and operational efficiency.
The mean operating cost per delivery serves as the metric for assessing cost-effectiveness, while the operational efficiency is evaluated using the mean lead time (see Online Appendix \ref{ec.insights} for calculation details).

In our coordinated system, drones handle the transit tasks between vertiports, while couriers ride delivery motorcycles for parcel collection and distribution, featuring a hub-and-spoke network with coordination of drones and couriers (H\&S-D\&C).
The benefits of this system come from two aspects: (1) the implementation of the hub-and-spoke network; and (2) the employment of drones.
We develop two benchmarks as control groups to illustrate the benefits from these aspects:
\begin{itemize}[leftmargin=*]
    \item D2D-C: Couriers with delivery motorcycles provide dedicated door-to-door services on a fully connected network (Figure \ref{fig.ICES}(a));
    \item H\&S-C: Couriers with delivery motorcycles provide both hub-to-hub and hub-to-spoke/spoke-to-hub services (Figure \ref{fig.ICES}(b)). 
\end{itemize}
The corresponding mathematical models (see Online Appendix \ref{ec.insights}) for these benchmarks run on the same experimental settings.
Table \ref{tab.costTime} lists the mean operating costs and mean lead times under different pooling sizes $Q$.

\begin{table}[h]
  \centering
  \caption{Mean Operating Cost and Mean lead time for Our Proposed Network and Two Benchmarks}
    \scriptsize
    \begin{tabular}{lcccccc}
    \toprule
          & \multicolumn{3}{c}{Mean operating cost (CNY)} & \multicolumn{3}{c}{Mean lead time (min)} \\
          \cmidrule(r){2-4}    \cmidrule(r){5-7}
    \multicolumn{1}{l}{Pooling size} & D2D-C & H\&S-C & H\&S-D\&C & D2D-C & H\&S-C & H\&S-D\&C \\
    \midrule
    1     & 20.1  & 15.0  & \textbf{9.4}   & 53.1  & 47.9  & \textbf{22.7} \\
    2     & 10.7  & 9.5   & \textbf{6.8}   & 68.6  & 62.1  & \textbf{32.5} \\
    4     & 5.8   & 6.7   & \textbf{5.0}   & 103.2 & 81.4  & \textbf{40.1} \\
    6     & 4.1   & 5.8   & \textbf{4.3}   & 134.5 & 87.9  & \textbf{49.4} \\
    8     & 3.3   & 5.8   & \textbf{4.0}   & 166.4 & 87.9  & \textbf{58.3} \\
    10    & 2.7   & 5.1   & \textbf{3.8}   & 196.4 & 98.0  & \textbf{67.2} \\
    12    & 2.4   & 4.4   & \textbf{3.6}   & 227.1 & 116.4 & \textbf{76.6} \\
    \bottomrule
    \end{tabular}%
  \label{tab.costTime}%
\end{table}%

\textbf{Value of Network.}
This value can be understood by comparing D2D-C with H\&S-C. 
We first compare their operating costs.
Implementing a hub-and-spoke network even increase the operating cost by 15\%--83\% when the pooling size exceeds 4.
For the H\&S-C mode, the first- and last-mile delivery costs constitute a larger portion, up to 58\%, of the total cost as the pooling size increases.
The dedicated door-to-door delivery without the need for first- and last-mile services becomes more cost effective when the cost for dedicated couriers decreases with the larger pooling size. 
In addition, the proposed network does not necessarily yield cost savings due to the detour.

The comparison of mean lead times between D2D-C and our H\&S-D\&C suggests that, despite the detour, the proposed network allows operators to reduce lead times by 10\%--49\%.
This is because hubs consolidate parcels from different O-D pairs and aggregate flows on hub pairs, increasing the possibility of meeting the pooling threshold and reducing pooling times.
Conversely, dedicated courier services can experience three times longer pooling times due to the absence of consolidation.

\begin{myFind}
    The coordinated network enables efficiency gains by consolidating dispersed demand, while costs may not be necessarily reduced due to additional first- and last-mile trips.
\end{myFind}

\textbf{Value of Drones.}
This value can be identified by comparing H\&S-C with our H\&S-D\&C.
Table \ref{tab.costTime} indicates that after introducing drones on the hub pairs, the H\&S-D\&C offers 18\%--37\% lower mean operating costs compared to H\&S-C. 
The previous argument against the cost-effectiveness of drones in the logistics industry stems from the relatively high acquisition costs of drones. 
In contrast, our results demonstrate that drones can offset their ownership costs through their efficiency advantages. 
Drones have higher turnover efficiency than couriers, enabling operators to fulfill the demand with a small fleet, which is more economical than employing dedicated couriers.

The comparison of mean lead times reveals that drones enable an efficiency improvement of 34\%--53\%.
This advantage also stems from the drones' ability to move swiftly and directly.
On the other hand, couriers on motorcycles may face traffic congestion and detours, resulting in an average travel time that is nearly three times the flight time of drones.

\begin{myFind}
    Drones can simultaneously promote the cost-efficiency and timeliness of deliveries through their superior turnover efficiency within the closed operations network.
\end{myFind}

\begin{figure}[H]
    \centering
    \caption{Cost-Efficiency Analysis of Different Modes}
    \includegraphics[width=0.5\linewidth]{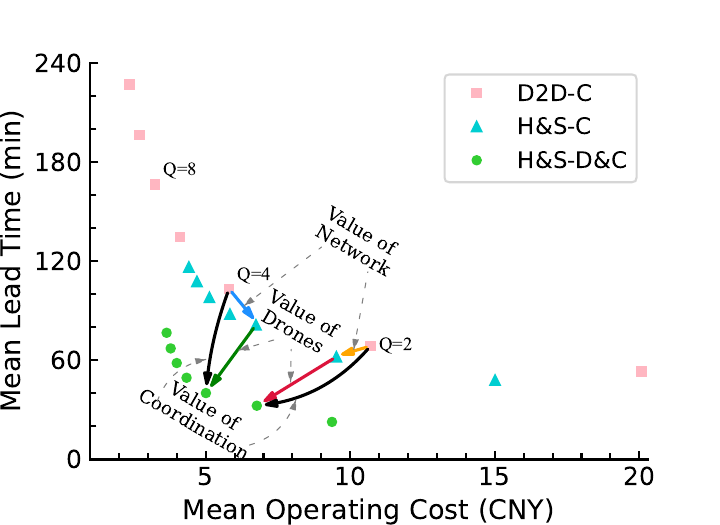}
    \label{fig.cost-efficiency}
    \begin{flushleft}
        {\footnotesize \textit{Note}. We generate scatter plots for each mode by varying the pooling size $Q$. From left to right, the points for each group correspond to the pooling sizes $Q=12, 10, 8, 6, 4, 2, 1.$}
    \end{flushleft}
\end{figure}

\textbf{Value of Coordination.}
By combining the benefits of the coordinated network and drones, the value of coordination can be recognized.
Figure \ref{fig.cost-efficiency} illustrates a cost-efficiency analysis of our proposed H\&S-D\&C and two benchmarks. 
In this figure, an arc pointing towards the bottom left corner signifies a reduction in both lead time and operating cost, whereas an arc directed towards the bottom right corner indicates a decrease in lead time but an increase in operating cost.
The value of coordination, represented by the dark arcs, is the combined sum of the value of network (the blue/orange arcs) and the value of drones (the green/red arcs).
Leveraging drones consistently enhances the efficiency benefits of the coordinated network, as demonstrated by the ``downward'' dark arcs.
Moreover, the operating cost reduction (the orange arc) associated with a small pooling size ($Q \leq 2$), can be further augmented by drones (the red arc), while drones can help mitigate potential operating cost increase (the blue arc) related to a large pooling size ($Q \geq 4$).

Overall, Figure \ref{fig.cost-efficiency} demonstrates that the H\&S-D\&C is not dominated by the D2D-C and H\&S-C, offering a Pareto improvement in both lead time and operating cost.
Although the door-to-door courier-based delivery may provide a lower operating cost than our coordinated delivery when the pooling size $Q$ exceeds 8, average lead times exceeding 2.5 hours are impractical, considering the two-hour time window usually set by customers.


\begin{myFind}
    Within the coordinated network, the coordination between drones and couriers can stride a win-win situation, enhancing delivery efficiency and reducing operating costs.
\end{myFind}

\subsubsection{Features of Coordinated Network.}  \label{sec.networkStructure}

\begin{figure}
    \centering
    \caption{Network Structure, Selected Vertiports, and Selected Service Routes}
    \includegraphics[width=0.95\linewidth]{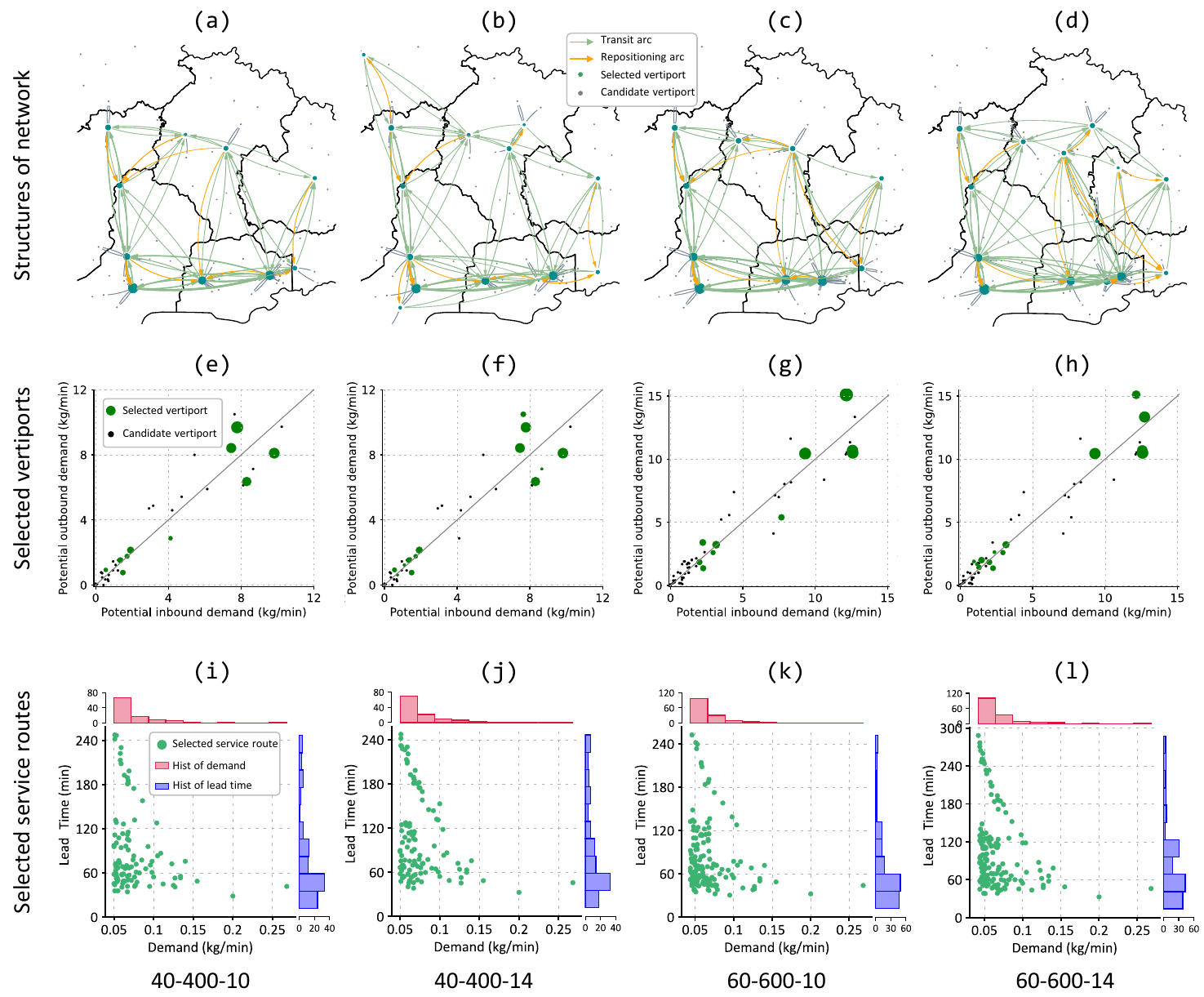}
    \label{fig.networkStructure}
    \begin{flushleft}
        {\footnotesize \textit{Note}. (a)--(d) Network structures, with line thickness representing the volume of flow and the size of green points indicating the total amount of demand across the vertiport. (e)--(h) Selected vertiports, with the size of point representing the total amount of demand across the vertiport. (i)--(l) selected service routes with their lead time and demand. The first column ((a), (e), and (i)) is for an instance in 40-400-100. The second, third, and fourth columns are for instances in 40-400-14, 60-600-10, and 60-600-14, respectively. }
    \end{flushleft}
\end{figure}

We analyze the coordinated networks focusing on network structure, vertiports, and service routes, as illustrated in Figure \ref{fig.networkStructure}, along with an examination of the interplay between network design and operations.

\textbf{Network Structure.}
The coordinated network represents a ``concentrated structure'' (Figures \ref{fig.networkStructure}(a)--(d)), where all selected vertiports are synergistically connected, rather than scattered in separate sub-networks.
One benefit of a concentrated network is increased consolidation. 
In such a network, vertiports have the potential to collect more parcels, thereby enhancing the efficiency of parcel pooling and reducing the lead time.
Moreover, the concentrated network can lead to a synergistic effect.
Interconnected vertiports allow the operator to coordinate transit and repositioning flows collaboratively, reducing deadhead repositioning costs.
Additionally, this network provides opportunities for spatial coupling of different O-D pairs.  

\textbf{Selected Vertiports.}
Let $D_i^I=\sum_{o,d\in\mathcal{R},j\in\mathcal{N}}D_{o,d}\delta_{o,i,j,d}$ and $D_i^O= \sum_{o,d\in\mathcal{R},i\in\mathcal{N}}D_{o,d}\delta_{o,i,j,d}$ denote the potential inbound and outbound demands of locations $i\in\mathcal{N}$, and let $D_{\max}^P=\max\{D_i^I,D_i^O, \forall i\in\mathcal{N}\}$ denote the maximum potential value.
We define the low-, medium-, and high-potential location sets as $\mathcal{N}^L:=\{i: D_i^I\leq\frac{1}{3}D_{\max}^P, D_i^O\leq\frac{1}{3}D_{\max}^P,\forall i\in\mathcal{N}\}$, $\mathcal{N}^M:=\{i: D_i^I\leq\frac{2}{3}D_{\max}^P, D_i^O\leq\frac{2}{3}D_{\max}^P,\forall i\in\mathcal{N}\setminus\mathcal{N}^L\}$, and $\mathcal{N}^H:=\mathcal{N}\setminus\{\mathcal{N}^L\cup\mathcal{N}^M\}$, respectively.
Figures \ref{fig.networkStructure}(e)--(h) indicate that candidates in the high-potential set $\mathcal{N}^H$ and the low-potential set $\mathcal{N}^L$ are more likely to be selected, while those in the medium-potential set $\mathcal{N}^M$ exhibit the lowest chance of selection.
This result is surprising, in that the likelihood of selection does not correspond to the attractiveness of the vertiports.

To validate this counterintuitive result, we conduct experiments on additional instances.
Table \ref{tab.demandClass} further indicates that medium-potential vertiports surprisingly cater to no more than 5\% of the demand (out of nearly 80\% of instances), a figure significantly lower than the average of 25\% for low-potential vertiports.
Geographically, medium-potential vertiports are often spread out around the high-potential vertiports, leading to a loss of their competitive edges.
High-potential vertiports exhibit a ``siphon effect'' on the market in their vicinity.  
These results also suggest a heuristic for reducing the pool of candidate locations.

\begin{table}[h]
  \centering
  \scriptsize
  \caption{Proportions of Demand Served by Vertiports Falling into Low, Medium, and High Potential Classes}
    \begin{tabular}{lccc}
    \toprule
    Percentage & \multicolumn{1}{l}{Low-potential} & \multicolumn{1}{l}{Medium-potential} & \multicolumn{1}{l}{High-potential} \\
    \midrule
    (20,200,10) & 21.4\% & 10.8\% & 67.9\% \\
    (20,200,14) & 22.2\% & 10.7\% & 67.1\% \\
    (30,300,10) & 21.8\% & 0.0\% & 78.2\% \\
    (30,300,14) & 27.9\% & 0.0\% & 72.1\% \\
    (40,400,10) & 25.1\% & 3.2\% & 71.7\% \\
    (40,400,14) & 25.1\% & 3.2\% & 71.7\% \\
    (50,500,10) & 26.2\% & 0.0\% & 73.8\% \\
    (50,500,14) & 28.3\% & 0.0\% & 71.7\% \\
    (60,600,10) & 24.9\% & 5.0\% & 70.0\% \\
    (60,600,14) & 30.8\% & 0.0\% & 69.2\% \\
    \bottomrule
    \end{tabular}%
  \label{tab.demandClass}%
\end{table}%

\begin{myFind}
    The vertiports selected tend to have higher or lower potential to attract demand, while medium-potential candidates are rarely selected due to the competition from high-potential ones.
\end{myFind}

\textbf{Selected Service Routes.} The selected service routes can be coupled efficiently to guarantee a rapid delivery process.
Figures \ref{fig.networkStructure}(i)--(l) illustrate that more than 70\% of the service routes can provide a 2-hour delivery, even when the demand is very low.
Through the consolidation facilitated by vertiports, parcels on some low-demand routes are pooled together, thereby aggregating the flow on vertiport pairs and enhancing the pooling process and delivery efficiency.

On the other hand, there are still a few service routes, about 10\%, where the lead times exceed 3 or even 4 hours.
We find that transit flows on these routes serve the dual purpose of repositioning.
To demonstrate this, we exclude service routes with lead times $t^L$ exceeding 240, 180, and 150 minutes, respectively. 
Table \ref{tab.repositioning} shows the rise in the proportion of deadhead repositioning trips after imposing time restrictions.
The tightest restriction even leads to a 36\% surge in repositioning, highlighting the role of transit flows in repositioning within these routes. 


\begin{table}[H]
  \centering
  \scriptsize
  \caption{Changes in Repositioning Trips after Imposing Restriction on Long-duration Service Routes}
    \begin{tabular}{p{5em}p{7em}<{\centering}p{5em}<{\centering}p{5em}<{\centering}p{5em}<{\centering}}
    \toprule
          & \multicolumn{1}{c}{Baseline} & \multicolumn{3}{c}{Changes (\%) after restricting routes with }   \\
          \cmidrule(r){2-2}    \cmidrule(r){3-5}        & Repositioning trips per day & \multicolumn{1}{c}{ $t^L\geq$240} & \multicolumn{1}{c}{$t^L\geq$180} & \multicolumn{1}{c}{$t^L\geq$150} \\
    \midrule
    40-400-10 & 61.7  & 5.2   & 24.2  & 36.0 \\
    40-400-14 & 86.2  & 13.9  & 20.8  & 27.0 \\
    60-600-10 & 67.0  & 3.8   & 9.3   & 21.0 \\
    60-600-14 & 69.1  & 3.0   & 10.9  & 10.1 \\
    \bottomrule
    \end{tabular}%
  \label{tab.repositioning}%
\end{table}%

\begin{myFind}
    Some long lead-time transit trips serve the dual purpose of repositioning drones, effectively reducing the need for deadhead repositioning trips.
\end{myFind}

\textbf{Interplay between Network Design and Operations.}
Understanding how network design intertwines with operations performance can guide future operators in crafting long-term network planning and striking a balance between cost and efficiency.
Table \ref{tab.iterplay} illustrates the interplay between network design and operations.
The operating costs decrease with the network size as larger networks enhance service accessibility, thereby lowering first- and last-mile delivery costs.
A larger network, however, does not always lead to positive outcomes. 
It may require a 2.7\%--16.2\% larger fleet size to sustain service, resulting in decreased drone utilization rates by 6.4\%--17.5\%.
The lead time could increase by up to 32.4\% as a decentralized network could undermine the consolidation function of vertiports.
Lastly, the expanded network also results in a lower service level.
These results highlight the role of network design in balancing the operating costs and performance.

\begin{table}[h]
  \centering
  \scriptsize
  \caption{Interplay between Network Design and Operations}
    \begin{tabular}{p{12em}p{5em}<{\centering}p{5em}<{\centering}p{5em}<{\centering}p{5em}<{\centering}p{5em}<{\centering}p{5em}<{\centering}p{5em}<{\centering}}
    \toprule
          & \multicolumn{1}{c}{Baseline} & \multicolumn{5}{c}{Changes (\%)} \\
    \cmidrule(r){2-2}    \cmidrule(r){3-8}             & Nine vertiports & Ten vertiports & Eleven vertiports & Twelve vertiports & Thirteen vertiports & Fourteen vertiports & Fifteen vertiports\\
    \midrule
    Operating cost (CNY) & 20,586 & -3.6  & -5.6  & -7.4  & -9.0  & -10.4 & -11.6 \\
    First/Last mile cost (CNY) & 14,574 & -4.6  & -8.8  & -11.3 & -13.2 & -16.3 & -18.6 \\
    Fleet cost (CNY) & 2,664  & 2.7   & 8.1   & 8.1   & 8.1   & 13.5  & 16.2 \\
    Fleet utilization rate & 65.2\% & -6.4  & -9.8  & -10.0 & -10.8 & -15.1 & -17.5 \\
    Lead time (min) & 72.5  & 5.7   & 12.9  & 19.6  & 27.2  & 28.3  & 32.4 \\
    Avg. Service level & 0.76  & -1.8  & -2.9  & -5.6  & -8.0  & -8.6  & -10.2 \\
    \bottomrule
    \end{tabular}%
  \label{tab.iterplay}%
\end{table}%

\subsubsection{Drivers of Potential Benefits.}
\label{sec.profitability}

The coordinated drone-courier city logistics are currently in the early stages of pilot programs. 
The potential benefits of this coordinated system, in the midst of the rapid evolution of the delivery industry \citep{he2022smart} and drone technology \citep{baloch2020strategic}, remain largely unknown.
To illuminate this novel logistics mode, we assess its benefits in light of market (demand), technical, and operational changes.

We concern benefits in terms of economical efficiency (operating cost) and operational efficiency (lead time).
Table \ref{tab.sensitiveAnalysis} concludes the operating costs and lead times with demand parameter (market share $\xi$), technical parameters (flight range $L^{\rm D}$ and flight speed $l/t^f$), and operational parameters (drone acquisition cost $c^{\mathrm{E}}$, flight cost $p^{\mathrm{F}}$, and courier cost $p^{\mathrm{LT}}$).

\begin{table}[htbp]
  \centering
  \scriptsize
  \caption{Impact of Demand, Technical, and Operations Parameters on Operating Cost and Lead Time}
    \begin{tabular}{rllcccc}
    \toprule
          & Baseline  &       & \multicolumn{4}{c}{Changes (\%) after} \\
\cmidrule(r){2-3}  \cmidrule(r){4-7}
& Parameters & Metrics (values) & \multicolumn{2}{c}{Decreasing parameter by} & \multicolumn{2}{c}{Increasing parameter by} \\
\cmidrule(r){4-5}  \cmidrule(r){6-7}          &       &       & 50\%    & 25\%    & 25\%    & 50\% \\
    \midrule
    \multicolumn{1}{l}{Demand side} & Market share $\xi$ & Avg. Operating cost (3.78 CNY) & -20.6 & -4.1  & 12.1  & 13.5 \\
          &       & Avg. Lead time (72.8 min) & -97.2 & -19.5 & 2.9   & 25.9 \\
          &       & Min. \#vertiports (9) & -66.7 & -33.3 & 22.2  & 88.9 \\
    \midrule
    \multicolumn{1}{l}{Technical side} & Flight range $L_{\rm{D}}$ & Avg. Operating cost (3.78 CNY) & -     & 20.2  & 3.5   & 1.9 \\
          &       & Avg. Lead time (72.8 min) & -     & 35.5  & -2.8  & -0.6 \\
          &       & Min \#vertiports (9) & -     & 111.1 & -22.2 & -22.2 \\
    \cmidrule{2-7}          & Flight speed $l/t^f$ & Avg. Operating cost (3.78 CNY) & 3.0   & 0.9   & -0.7  & -1.0 \\
          &       & Avg. Lead time (72.8 min) & 17.7  & 5.8   & -3.6  & -5.9 \\
          &       & Min \#vertiports (9) & 0.0   & 0.0   & 0.0   & 0.0 \\
    \midrule
    \multicolumn{1}{l}{Operational side} & Vehicle cost $c^{\mathrm{E}}$ & Avg. Operating cost (3.78 CNY) & -6.5  & -3.3  & 3.2   & 6.3 \\
          & Flight cost $p^{\mathrm{F}}$ & Avg. Operating cost (3.78 CNY) & -8.2  & -4.1  & 4.0   & 8.0 \\
          & Courier cost $p^{\mathrm{LT}}$ & Avg. Operating cost (3.78 CNY) & -35.8 & -17.8 & 17.7  & 35.3 \\
    \bottomrule
    \end{tabular}%
  \label{tab.sensitiveAnalysis}%
\end{table}%

On the demand side, a large market size $\xi$ has potentially adverse effects on both economic and operational efficiency.
To maintain a significant market share, the operator has to operate a more extensive network and cater to some unprofitable O-D pairs (involving lengthy first- and last-mile transportation), potentially increasing operating costs and lengthening lead times.
In contrast, the impacts of parameters on the technical side are positive but marginal.
Even a 50\% increase in a drone's flight range or speed results in a negligible change (less than 2\%) in operating costs.
Fortunately, technical advancements are beneficial for enhancing operating efficiency and maintaining a small network. 
Overall, costs (especially expenses for couriers) incurred on the operational side have a more significant impact on cost-effectiveness.
A 25\% reduction in courier delivery costs results in a 17.8\% saving in operating costs, and this figure increases to 35.8\% when the courier transport costs are reduced by 50\%.
In contrast, the impacts of vehicle and flight costs on operating expenses are relatively modest (within 8.2\%).

From a managerial perspective, these findings suggest that a critical driver in improving cost-effectiveness is managing expenses in the courier delivery. 
Optimizing delivery routes and courier assignment in the first- and last-mile sectors could be promising strategies for cost reduction. 
On the other hand, technical innovations can enhance operational efficiency. 
Creating fast and long-range drones, for example, can boost network turnover efficiency and reduce lead times.

\section{Conclusion}\label{sec.conclusion}

This paper introduced a coordinated drone-courier logistics system for the intra-city express service.
To shed light on this system, we developed a unified optimization model that combines a hub location model, which optimizes the network design (specially planning vertiport locations), with a queueing network model that captures the interplay between network design and tactical operations.
Due to the queueing-related formulations, the model exhibits a feasible region that is partially convex and partially concave.
To address the intractability, we employed the piecewise linear approximation to develop an exact algorithm, the adaptive discretization algorithm, that provides optimality guarantees. 
The computational experiments show the scalability of algorithms.

We demonstrated the application of our methodology to reveal the coordinated drone-courier logistics using real data. 
Our results uncovered several managerial implications towards the this system.
First, this system can simultaneously save operating costs and promote the delivery efficiency of city logistics, thanks to the implementation of the coordinated network and the use of drones. 
Second, the network presents a concentrated structure with interconnected vertiports and collaborative traffic dispatching, enhancing the consolidation of scattered demand and reducing deadhead trips. 
Ultimately, to enhance the sustainable development, efforts should be directed not only towards technical advancements but also towards optimizing costs in courier delivery.

Looking ahead, we identify two promising research directions to extend this study. 
One avenue is to apply the proposed methodology to unlock generalized coordinated unmanned vehicle-human logistics systems, such as the coordinated autonomous vehicle-courier delivery system and the drone-based emergency medical care system.
Another potential avenue is to combine existing services, such as truck delivery or dedicated courier service, to design differentiated services for customers with different time and cost sensitivities.
Customizing service options with different prices and delivery times is essential for operators to cater to diverse customer demographics.

\baselineskip=14pt
\def\bibfont{\small}
\bibliographystyle{informs2014}
\bibliography{mybibfile01}

\ECSwitch 
\nociteannex{*}
\ECHead{E-Companion to ``On Coordinated Drone-Courier Logistics for Intra-city Express Services''}

\section{Details on Modelling (Section \ref{sec.model})}
\label{ec.notations}
Notations used in this paper. 

{\small
\renewcommand\arraystretch{1.3}
\begin{longtable}{p{1.5cm}p{14cm}} 
\toprule
Notation       & Description\\
\midrule
    \multicolumn{2}{l}{\textit{\textbf{Sets:}}} \\
    $\mathcal{N}$       & Set of nodes (candidate vertiport locations);\\
    $\mathcal{H}$       & Set of possible capacity values; \\
    $\mathcal{R}$       & Set of origins and destinations. \\
    \multicolumn{2}{l}{\textit{\textbf{Input Parameters:}}} \\
    $c^{\mathrm{E}}$            & Daily cost for each drone;\\
    $p_{ij}^{\mathrm{F}}$       & Drone transport cost (per kg) on vertiport pair $(i,j)$; \\
    $D_{od}$       & Demand on O-D pair $(o,d)$;\\
    $p_{o,i}^{\mathrm{LT}}, p_{j,d}^{\mathrm{LT}}$   & First/Last-mile courier delivery cost (per kg) from origin $o$ to vertiport $i$ or from vertiport $j$ to destination $d$;\\
    $P$       & The maximum number of vertiports to be selected;\\
    $\xi$     & Expected market share of drone delivery;\\
    $l_{o,i}, l_{j,d}$ & Courier travel distance from origin $o$ to vertiport $i$ or from vertiport $j$ to destination $d$;\\
    $l_{i,j}$ & Drone flight distance between vertiports $i$ and $j$;\\
    $L^{\mathrm{S}}$ & Maximum service range of vertiports;\\
    $L^{\mathrm{D}}$ & Maximum flight distance of drones;\\
    $Q$              & Pooling size of parcels in a single flight; \\
    $t_{i,j}^f$      & Flight time between vertiports $i$ and $j$;\\
    $b^{\mathrm{chg}},b^{\mathrm{fly}}$ & Battery charging and depletion rate.\\
    \multicolumn{2}{l}{\textit{\textbf{Decision Variables:}}} \\
    $x_{i}$   &   Binary variable, equals 1, if a vertiport is built in location $i \in \mathcal{N}$, and 0, otherwise; \\
    $z_{i,h}$   &   Binary variable, equals 1, if a vertiport is built in location $i \in \mathcal{N}$ with capacity $h$, and 0, otherwise;\\
    $y_{o,i,j,d}$     & Binary variable, equals 1, if demand between O-D pair $o, d \in \mathcal{R}$ is served by drones between vertiports $i,j \in \mathcal{N}$ , and 0, otherwise;\\
    $\Gamma$     & Integer variable, fleet size;\\
    $\rho_i$     & Continuous variable, service level demonstrating the percentage of parcels at vertiport $i$ that can be transported by drones;\\
    $\psi_{i,j}$       & Continuous variable, rate of transit flow from vertiports $i$ to $j$;\\
    $\varphi_{i,j}$    & Continuous variable, rate of repositioning trips from vertiports $i$ to $j$;\\
    $\alpha_{o,i,j,d}$ & Auxiliary continuous variable; \\
    $\beta_{i}^k$      & Auxiliary binary variable; \\
    $\Theta_{i}$       & Auxiliary continuous variable; \\
    $\pi_{i}^{k}$      & Auxiliary continuous variable. \\
    \bottomrule
\end{longtable}
}

\section{Details on Algorithm (Section \ref{sec.algorithm})}

\subsection{Mathematical Proofs}






\subsubsection{Proof of Proposition \ref{Pro.1}}~


To prove $Z[G_R]\leq Z[G^*] \leq Z[G_C]$, we can initially prove that: 
\begin{enumerate}[fullwidth,itemindent=2em,label=(\arabic*)]
    \item  Model $[G^*]$ is a relaxation of model $[G_C]$;
    \item  Model $[G_R]$ is a relaxation of model $[G^*]$.
\end{enumerate}

A model $[A]$ can be considered a relaxation of a model $[B]$ if the feasible region of $[B]$ is entirely contained within that of $[A]$, and the objective function value of $[A]$ is consistently lower than that of $[B]$ for the same solution (in minimization problems, and vice versa in maximization problems) \citep{wolsey2020integer}.
Since the objective functions of $[G^*]$, $[G_C]$, and $[G_R]$ are identical, we only need to prove the containment relationship of the feasible regions.

We first prove that $Z[G^*] \leq Z[G_C]$.
Let $(\bm{\chi}^C,\bm{\nu}^C)$ be an arbitrary feasible solution of $[G_C]$.
$(\bm{\chi}^C,\bm{\nu}^C)$ satisfies that $\sum_{k=1}^K\overline{\bm{h}}^k(\bm{\chi}^C)\mathbbm{1}(\bm{\chi}^C\in\mathcal{X}^k) + \bm{A}_1\bm{\nu}^C \leq \bm{0}$ and $\sum_{k=1}^K\overline{\bm{q}}^k(\bm{\chi}^C)\mathbbm{1}(\bm{\chi}^C\in\mathcal{X}^k) + \bm{A}_2\bm{\nu}^C \leq \bm{0}$.
$\bm{h}(\bm{\chi}^C) \leq \sum_{k=1}^K\overline{\bm{h}}^k(\bm{\chi}^C)\mathbbm{1}(\bm{\chi}^C\in\mathcal{X}^k)$ always holds true, as within each partition $\mathcal{X}^k$, $\overline{h}^k(\chi) \geq h(\chi)$ is always satisfied.
Therefore, $\bm{h}(\bm{\chi}^C) + \bm{A}_1\bm{\nu}^C \leq \bm{0}$ must be satisfied if $\sum_{k=1}^K\overline{\bm{h}}^k(\bm{\chi}^C)\mathbbm{1}(\bm{\chi}^C\in\mathcal{X}^k) + \bm{A}_1\bm{\nu}^C \leq \bm{0}$ is satisfied.
Similarly, $\bm{q}(\bm{\chi}^C) + \bm{A}_2\bm{\nu}^C \leq \bm{0}$ must be satisfied as $\bm{q}(\bm{\chi}^C) + \bm{A}_2\bm{\nu}^C \leq \sum_{k=1}^K\overline{\bm{q}}^k(\bm{\chi}^C)\mathbbm{1}(\bm{\chi}^C\in\mathcal{X}^k) + \bm{A}_2\bm{\nu}^C \leq \bm{0}$ is satisfied.
Therefore, any feasible solution of $[G_C]$ is feasible for $[G^*]$.
This means that the feasible region of $[G_C]$ is entirely contained within that of $[G^*]$.
$[G^*]$ is a relaxation of $[G_C]$, and its optimum $Z[G^*]$ provides a lower bound on $Z[G_C]$, i.e., $Z[G^*] \leq Z[G_C]$.
In others words, $Z[G_C]$ provides an upper bound on $Z[G^*]$.

Next, we prove that ${Z}[G_R]\leq Z[G^*]$.
Let $(\bm{\chi}^*,\bm{\nu}^*)$ be an arbitrary feasible solution of $[G^*]$.
$(\bm{\chi}^*,\bm{\nu}^*)$ satisfies that $\bm{h}(\bm{\chi}^*)+\bm{A}_1\bm{\nu}^*\leq \bm{0}$ and $\bm{q}(\bm{\chi}^*)+\bm{A}_2\bm{\nu}^*\leq \bm{0}$.
As defined in (\ref{eq.general.mr}),  we can deduce that $\bm{h}(\bm{\chi}^*)\geq \sum_{k=1}^K\underline{\bm{h}}^k(\bm{\chi}^*)\mathbbm{1}(\bm{\chi}^*\in\mathcal{X}^k)$ and $\bm{q}(\bm{\chi}^*)\geq \sum_{k=1}^K\underline{\bm{q}}^k(\bm{\chi}^*)\mathbbm{1}(\bm{\chi}^*\in\mathcal{X}^k)$.
Therefore, we can get that $\sum_{k=1}^K\underline{\bm{h}}^k(\bm{\chi}^*)\mathbbm{1}(\bm{\chi}^*\in\mathcal{X}^k) + \bm{A}_1\bm{\nu}^*\leq \bm{0}$ and $\sum_{k=1}^K\underline{\bm{q}}^k(\bm{\chi}^*)\mathbbm{1}(\bm{\chi}^*\in\mathcal{X}^k) + \bm{A}_2\bm{\nu}^*\leq \bm{0}$, showing that any feasible solution $(\bm{\chi}^*,\bm{\nu}^*)$ of $[G^*]$ is also feasible for $[G_R]$.
We can thus infer that the feasible region of $[G^*]$ is completely contained within that of $[G_R]$.
Therefore, the optimum of $[G_R]$ provide a lower bound on $Z[G^*]$, i.e., $Z[G_R]\leq Z[G^*]$. 
This completes the proof. 
\hfill\Halmos

\subsubsection{Proof of Proposition \ref{Pro.2}}~

Let $(\bm{\chi}^R,\bm{\nu}^R)$ be a maximizer of $g(\cdot)$ over $\mathcal{X}^{\tilde{\bm{k}}}$.
Let us fix a particular element $\bm{\chi}^R \in \mathcal{X}^{\tilde{\bm{k}}}$ in $[G^*]$. 
Model $[G^*]$ thus can be converted to a linear program: 
\begin{eqnarray}
    [\tilde{G}] = \min\{\bm{c}^{\top}\bm{\nu}:\bm{A}_1\bm{\nu} \leq -\bm{h}(\bm{\chi}^R), \bm{A}_2\bm{\nu} \leq -\bm{q}(\bm{\chi}^R), \bm{\nu} \in \mathcal{V}\}. \nonumber
\end{eqnarray}

Suppose that there exists a non-degenerate primal optimal basic solution.
Let matrix $\bm{B}$ and vector $\bm{\nu}_{\bm{B}}$ denote the corresponding basis matrix and basic variables of model $[\tilde{G}]$. 
The optimal solution of $[\tilde{G}]$ can be denoted by $(\bm{\chi}^R,\tilde{\bm{\nu}})$, in which $\tilde{\bm{\nu}}=[\tilde{\bm{\nu}}_{\bm{B}};\tilde{\bm{\nu}}_{\bm{N}}]$, $\tilde{\bm{\nu}}_{\bm{B}} = \bm{B}^{-1}{\bm{b}}$ is the optimal solution of basic vector, $\tilde{\bm{\nu}}_{\bm{N}}=\bm{0}$ is non-basic vector, and $\bm{b}=[-\bm{h}(\bm{\chi}^R), -\bm{q}(\bm{\chi}^R)]^{\top}$ is the RHS constant vector.
The optimal objective function value of $[\tilde{G}]$ is given by:
\begin{eqnarray}
    g(\tilde{\bm{\nu}})=\bm{c}_{\bm{B}}^{\top}\bm{B}^{-1}\left[
    \begin{array}{c}   
        - {\bm{h}}(\bm{\chi}^R)\\
        - {\bm{q}}(\bm{\chi}^R)
    \end{array}
    \right]. \nonumber
\end{eqnarray}

Let $(\bm{\chi}^R,\bm{\nu}^R)$ denote an optimal solution of $[G_R]$. 
If we fix $\bm{\chi}=\bm{\chi}^R$ in $[G_R]$, we get a mixed-integer linear program: 
\begin{eqnarray}
    [\hat{G}] = \min\{\bm{c}^{\top}\bm{\nu}:\bm{A}_1\bm{\nu} \leq -\underline{\bm{h}}(\bm{\chi}^R), \bm{A}_2\bm{\nu} \leq -\underline{\bm{q}}(\bm{\chi}^R), \bm{\nu} \in \mathcal{V}\}. \nonumber
\end{eqnarray}
$\bm{\nu}^R$ is still an optimal solution of $[\hat{G}]$.
Model $[\hat{G}]$ is quite similar to $[\tilde{G}]$ except for the difference in the RHS constant vector.
Let $\bm{b}^R=[-\underline{\bm{h}}(\bm{\chi}^R), -\underline{\bm{q}}(\bm{\chi}^R)]^{\top}$ denote the RHS constant vector of $[\hat{G}]$.
If the difference between $\bm{b}$ and $\bm{b}^R$ is small enough such that $\bm{B}^{-1}\bm{b}$ remains positive and we still have a basic feasible solution.
The reduced cost $\bm{c}_{\bm{B}}^{\top}\bm{B}^{-1}$ is not affected by the change from  $\bm{b}$ to $\bm{b}^R$ and remain non-negative. 
It means that $\bm{B}$ remains an optimal basis for the new problem as well. 
The corresponding optimal solution can be obtained by $\bm{\nu}^R = [\bm{\nu}_{\bm{B}}^R; \bm{\nu}_{\bm{N}}^R]$, in which $\bm{\nu}_{\bm{B}}^R = \bm{B}^{-1}\bm{b}^R$, and $\bm{\nu}_{\bm{N}}^R=\bm{0}$.
The optimal objective value $g(\bm{\nu}^{R})$ for the new problem is given by: 
\begin{eqnarray}
    && g(\bm{\nu}^R) = \bm{c}_{\bm{B}}^{\top}\bm{B}^{-1} \left[
    \begin{array}{c}   
        - \underline{\bm{h}}(\bm{\chi}^R)\\
        - \underline{\bm{q}}(\bm{\chi}^R)
    \end{array}
    \right]. \nonumber
\end{eqnarray}

The difference between $g(\tilde{\bm{\nu}})$ and $g(\bm{\nu}^R)$ is given by:
\begin{eqnarray}
     g(\tilde{\bm{\nu}}) - g(\bm{\nu}^R) &=& \bm{c}_{\bm{B}}^{\top}\bm{B}^{-1} \left[
    \begin{array}{c}   
        \underline{\bm{h}}(\bm{\chi}^R) - \bm{h}(\bm{\chi}^R) \\
        \underline{\bm{q}}(\bm{\chi}^R) - \bm{q}(\bm{\chi}^R)
    \end{array}
    \right]  \nonumber\\
    &\leq &
    \Vert \bm{c}_{\bm{B}}^{\top}\bm{B}^{-1} \Vert_2 \cdot \left\Vert \left[
    \begin{array}{c}   
        \bm{h}(\bm{\chi}^R) - \underline{\bm{h}}(\bm{\chi}^R) \\
        \bm{q}(\bm{\chi}^R) - \underline{\bm{q}}(\bm{\chi}^R)
    \end{array}
    \right] \right\Vert_2. \nonumber
\end{eqnarray}
The inequality  comes from the Cauchy–Schwarz inequality.
Since $h(\chi)$ and $q(\chi)$ are $L_h$- and $L_q$-Lipschitz continuous and $L=\max\{|L_q|, |L_h|\}$, we can easily get ${h}(\bm{\chi}^R) - \underline{{h}}(\bm{\chi}^R)) \leq L \Vert \mathcal{X}^{\tilde{\bm{k}}}\Vert_{\infty}$ and ${q}(\bm{\chi}^R) - \underline{{q}}(\bm{\chi}^R)) \leq L \Vert \mathcal{X}^{\tilde{\bm{k}}} \Vert_{\infty}$.    
Let $\Delta=\Vert \mathcal{X}^{\tilde{\bm{k}}} \Vert_{\infty}$ denote the maximum diameter of discretization, and let $n_1$ and $n_2$ represent the dimensions of vector $\bm{h}(\cdot)$ and $\bm{q}(\cdot)$. 
We can get:
\begin{eqnarray}
    \left\Vert \left[
    \begin{array}{c}   
        \bm{h}(\bm{\chi}^R) - \underline{\bm{h}}(\bm{\chi}^R) \\
        \bm{q}(\bm{\chi}^R) - \underline{\bm{q}}(\bm{\chi}^R)
    \end{array}
    \right] \right\Vert_2 &=& \sqrt{\left(\bm{h}(\bm{\chi}^R) - \underline{\bm{h}}(\bm{\chi}^R)\right)^\top \left(\bm{h}(\bm{\chi}^R) - \underline{\bm{h}}(\bm{\chi}^R)\right) + \left(\bm{q}(\bm{\chi}^R) - \underline{\bm{q}}(\bm{\chi}^R)\right)^\top \left(\bm{q}(\bm{\chi}^R) - \underline{\bm{q}}(\bm{\chi}^R)\right)}\nonumber\\
    &\leq &\sqrt{n_1(L\Delta)^2 + n_2(L\Delta)^2} =  L\Delta\sqrt{n_1+n_2}.\nonumber
\end{eqnarray}
It holds that: 
\begin{eqnarray}
    && g(\tilde{\bm{\nu}}) - g(\bm{\nu}^R)  \leq L\Delta \sqrt{n_1+n_2}\Vert \bm{c}_{\bm{B}}^{\top}\bm{B}^{-1} \Vert_2. \nonumber
\end{eqnarray}
If $S\leq \frac{1}{L\sqrt{n_1+n_2}\Vert \bm{c}_{\bm{B}}^{\top}\bm{B}^{-1} \Vert_2}$ and $\Delta=S\varepsilon$, we can get $g(\tilde{\bm{\nu}}) - g(\bm{\nu}^R) \leq \varepsilon$.
Since the optimum of $[G_R]$ is achieved over $\mathcal{X}^{\tilde{\bm{k}}}$, $g(\bm{\nu}^R)$ is a lower bound of $[G^*]$, i.e., $g(\bm{\nu}^R)\leq Z[G^*]$.
Thus, $g(\tilde{\bm{\nu}}) \leq Z[G^*] + \varepsilon$.
Therefore, there exists a feasible solution $(\bm{\chi}^R,\tilde{\bm{\nu}})$ in $\mathcal{X}^{\tilde{\bm{k}}}$ within $\varepsilon$ of the optimum of $[G^*]$. \hfill \Halmos

\subsubsection{Proof of Proposition \ref{Pro.4}}~

We first prove that Constraints (\ref{eq.mc.Theta})--(\ref{eq.mc.value}) satisfy the general principles described in (\ref{eq.general.mc}) and are pessimistic approximations to Constraints (\ref{eq.m1.fleet}) and  (\ref{eq.m1.charge}). 

Let $\overline{f}_S^k(\rho_i)$ represents the secant determined by two adjacent points $\left(\varrho_i^k, f(\varrho_i^k)\right)$ and $\left(\varrho_i^{k+1}, f(\varrho_i^{k+1})\right)$. 
\begin{eqnarray}
    \overline{f}_S^k(\rho_i) = \frac{\rho_i - \varrho^k\varrho^{k+1}}{(1-\varrho^k)(1-\varrho^{k+1})}
    \label{eq.mc.secant}
\end{eqnarray}
Due to the convexity of $f(\rho_i)$, $f(\rho_i) \leq \overline{f}_S^k(\rho_i)$ is true when $\varrho_i^k \leq \rho_i \leq \varrho_i^{k+1}$.
Therefore, $f(\rho_i) \leq \sum_{k\in\mathcal{K}_i} \overline{f}_S^k(\rho_i) \beta_i^{k}$ always holds if $\varrho_i^k \leq \rho_i \leq \varrho_i^{k+1}$, i.e., $\overline{\beta}_i^k=1$ is guaranteed.
This proves that replacing $f(\rho_i)$ by $\sum_{k\in\mathcal{K}_i} \overline{f}_S^k(\rho_i) \overline{\beta}_i^{k}$ satisfies the general principles presented in (\ref{eq.general.mc}). 

Note that we can avoid the usage of the binary variable $\overline{\beta}_{i}^k$ due to the convexity of feasible region defined by piecewise linear terms $\overline{f}_S^k(\rho_i), \forall k\in\mathcal{K}_i$. 
Let a continuous variable $\Theta_i$ approximate $f(\rho_i)$ such that $\overline{\Theta}_i \geq \overline{f}_S^k(\rho_i), \forall k\in\mathcal{K}_i$, i.e., Constraints (\ref{eq.mc.Theta}).
We can get $\overline{\Theta}_i \geq \underset{k\in\mathcal{K}_i}{\max}\{\overline{f}_S^k(\rho_i)\}$.
In the other word, $\overline{\Theta}_i$ is bounded by the secant which has the maximum function value at $\rho_i$.
For a specific $\rho_i$, the secant $\overline{f}_S^k(\rho_i)$ corresponding to $\overline{\beta}_i^k=1$ always provides the tightest bound for $\overline{\Theta}_i$ because for $k'$ be such that $\overline{\beta}_i^{k'}=1$.
It holds that $\overline{f}_S^{k'}(\rho_i) \geq f(\rho_i) \geq \left(\overline{f}_S^{k}(\rho_i)\right)_{k\in\mathcal{K}_i\setminus {k'}}$.
Thus we can get $\underset{k\in\mathcal{K}_i}{\max}\{\overline{f}_S^k(\rho_i)\}=\sum_{k\in\mathcal{K}_i}\overline{f}_S^k(\rho_i)\overline{\beta}_i^k$.
We can also get $\overline{\Theta}_i \geq  \underset{k\in\mathcal{K}_i}{\max}\{\overline{f}_S^k(\rho_i)\} = \sum_{k\in\mathcal{K}_i}\overline{f}_S^k(\rho_i) \beta_i^{k}$.
Since $f(\rho_i) \leq \sum_{k\in\mathcal{K}_i} \overline{f}_S^k(\rho_i) \beta_i^{k}$, it holds that $f(\rho_i) \leq \overline{\Theta}_i$.
Therefore, from the constraints $\overline{\Theta}_i \geq \overline{f}_S^k(\rho_i), \forall k\in\mathcal{K}_i$, we can deduce $f(\rho_i) \leq \overline{\Theta}_i$.
This means that $\overline{\Theta}_i$ is a pessimistic approximation of $f(\rho)$ and is in line with the general principles in replacing convex function presented in (\ref{eq.general.mc}).
This completes the proof that Constraints (\ref{eq.mc.Theta})--(\ref{eq.mc.fleet}) are pessimistic approximation of Constraint (\ref{eq.m1.fleet}).

In linearizing Constraint (\ref{eq.m1.charge}), we choose tangents to approximate $-f(\rho_i)$ because tangents $\overline{f}_T(\cdot)$ yield downward approximations for $f(\rho_i)$ and $-\overline{f}_T(\cdot)$ yield upward approximations to $-f(\rho_i)$. 
We choose the tangent $\overline{f}_T^k(\rho_i)$ determined by point $(\rho_i,\frac{\rho_i}{1-\rho_i})$ (see Equation (\ref{eq.mc.chargeAppro})).
As shown in Figure \ref{fig.exampleOfDiscretization}, we discretize the domain of $\rho_i$ based on the intersection points of tangents.
$\overline{\varrho}_i^k$ and $\overline{\varrho}_i^{k+1}$ are intersection points of $\overline{f}_T^k$ with $\overline{f}_T^{k-1}$ and $\overline{f}_T^{k+1}$. 
Due to convexity of $f(\rho_i)$, we can get $f(\rho_i) \geq \overline{f}_T^k(\rho_i)$ and $-f(\rho_i) \leq -\overline{f}_T^k(\rho_i)$ when $\overline{\varrho}_i^{k} \leq \rho_i \leq \overline{\varrho}_i^{k+1}$.
Correspondingly, the concave term $ -f(\rho_i) \leq -\sum_{i\in\mathcal{N}} \sum_{k\in\mathcal{K}_i} \overline{f}_T^k(\rho_i) \overline{\beta}_i^{k}$ always holds if Constraints (\ref{eq.mc.rhobeta1})--(\ref{eq.mc.sumbeta}) are guaranteed.
However, replacing $ -f(\rho_i) $ by $-\sum_{i\in\mathcal{N}} \sum_{k\in\mathcal{K}_i} \overline{f}_T^k(\rho_i) \overline{\beta}_i^{k}$ may lead to a nonlinear term $\overline{\beta}_{i}^k \rho_i$.
We can linearize it as shown in Constraints (\ref{eq.mc.pirho1})--(\ref{eq.mc.value}). 
Finally, Constraints (\ref{eq.mc.charge})--(\ref{eq.mc.value}) are in line with general principle in (\ref{eq.general.mc}).
Constraints (\ref{eq.mc.charge})--(\ref{eq.mc.value}) are pessimistic approximation of Constraints (\ref{eq.m1.charge}).
\begin{eqnarray}
    \overline{f}_T^k(\rho_i) = \frac{\rho_i-(\varrho_i^k)^2}{(1-\varrho_i^k)^2}, \;\; \forall i\in\mathcal{N}. \label{eq.mc.chargeAppro}
\end{eqnarray}

However, the binary variable $\overline{\beta}_i^k$ cannot be avoided in linearizing Constraints (\ref{eq.m1.charge}).
Suppose there exists $\varXi_i \geq 0$. 
If  $ -f(\rho_i) \leq -\sum_{i\in\mathcal{N}} \sum_{k\in\mathcal{K}_i} \overline{f}_T^k(\rho_i) \overline{\beta}_i^{k} \leq -\varXi_i$, $-\varXi_i$ is a pessimistic approximation of $-f(\rho_i)$.
We can get $\varXi_i \leq \sum_{i\in\mathcal{N}} \sum_{k\in\mathcal{K}_i} \overline{f}_T^k(\rho_i) \overline{\beta}_i^{k}, \forall k\in\mathcal{K}_i$, i.e., $\varXi_i \leq \underset{k\in\mathcal{K}_i}{\min}\{\overline{f}_T^k(\rho_i)\}$.
$\varXi_i \leq \underset{k\in\mathcal{K}_i}{\min}\{\overline{f}_T^k(\rho_i)\}$ is a piecewise linear concave constraints which cannot be replaced directly by $\varXi_i \leq \overline{f}_T^k(\rho_i), \forall k\in\mathcal{K}_i$.

Next, we will prove that any feasible of Model $[M_C]$ are feasible for $[M^*]$. 
The only difference between $[M_C]$ and $[M^*]$ is that linear Constraints (\ref{eq.mc.Theta})--(\ref{eq.mc.value}) in $[M_C]$ replace nonlinear Constraints (\ref{eq.m1.fleet}) and (\ref{eq.m1.charge}) in $[M^*]$. 
For the linear approximations in line with general principles described in (\ref{eq.general.mc}), the yielded constrains reduce the feasible region of original nonlinear model. 
The original model $[M^*]$ is a relaxation of the conservative model $[M_C]$. 
This means that any feasible solution of $[M_C]$ is a feasible solution of $[M^*]$.
This completes the proof.  \hfill\Halmos


\subsubsection{Proof of Proposition \ref{Pro.5}}~

We first prove that Constraints (\ref{eq.mr.Theta}) and (\ref{eq.mr.fleet}) are optimistic approximations to Constraint (\ref{eq.m1.fleet}).
The RHS on Constraints (\ref{eq.mr.Theta}) is the tangent function $\underline{f}_T^k(\rho_i)$ (same as (\ref{eq.mc.chargeAppro})) encompassing point $(\varrho_i^k,f(\varrho_i^k))$.
Due to the convexity of $f(\rho_i)$, it holds that $\underline{f}_T^k(\rho_i) \leq f(\rho_i), \forall \rho_i\in[0,1), k\in\mathcal{K}_i$, i.e., $f(\rho_i) \geq \underset{k\in\mathcal{K}_i}{\max}\{\underline{f}_T^k(\rho_i)\}$.
For a specific $\rho_i$, $f(\rho_i)$ is actually bounded by only one tangent which have the largest function value at here.
Let $\underline{\Theta}_i$ be a continuous variable bounded by all tangents $\underline{f}_T^k(\rho_i), \forall k\in\mathcal{K}_i$, i.e., $\underline{\Theta}_i \geq \underset{k\in\mathcal{K}_i}{\max}\{\underline{f}_T^k(\rho_i)\}$.
However, we cannot get that $\underline{\Theta}_i$ is downward approximation of $f(\rho_i)$, i.e., $\underline{\Theta}_i \leq f(\rho_i)$, from $\underline{\Theta}_i \geq \underset{k\in\mathcal{K}_i}{\max}\{\underline{f}_T^k(\rho_i)\}$ and $f(\rho_i) \geq \underset{k\in\mathcal{K}_i}{\max}\{\underline{f}_T^k(\rho_i)\}$.
Constraint (\ref{eq.mr.fleet}) imposes an upper bound on $\underline{\Theta}_i$. 
If the imposed upper bound exactly equal to its lower bound, $\underline{\Theta}_i = \underset{k\in\mathcal{K}_i}{\max}\{\underline{f}_T^k(\rho_i)\} \leq f(\rho_i)$.
Otherwise, $\underline{\Theta}_i$ is not tightly bounded by Constraint (\ref{eq.mr.fleet}).
This means that $\underline{\Theta}_i$ is free in its feasible range and $\underline{\Theta}_i > \underset{k\in\mathcal{K}_i}{\max}\underline{f}_T^k(\rho_i)$ is possible. 
In this case, there is a slack $\underline{\Theta}_i-\underset{k\in\mathcal{K}_i}{\max}\{\underline{f}_T^k(\rho_i)\}$ that does not necessarily reduce the objective function. 
In summary, regardless of whether Constraint (\ref{eq.mr.fleet}) is relaxed or not, $\Tilde{\underline{\Theta}}_i = \underset{k\in\mathcal{K}_i}{\max}\{\underline{f}_T^k(\rho_i)\}$ is the real part of $\underline{\Theta}_i$ that contributes to objective function. 
Recall $f(\rho_i) \geq \underset{k\in\mathcal{K}_i}{\max}\{\underline{f}_T^k(\rho_i)\}$, the real activated part $\Tilde{\underline{\Theta}}_i \leq f(\rho_i)$. 
Constraints (\ref{eq.mr.Theta})--(\ref{eq.mr.fleet}) is in line with general principles presented in (\ref{eq.general.mr}) even though these may lead to slacks.
This proves that Constraints (\ref{eq.mr.Theta})--(\ref{eq.mr.fleet}) are optimistic approximations to Constraint (\ref{eq.m1.fleet}).

Here, we focus on the proof that Constraints (\ref{eq.mr.charge})--(\ref{eq.mr.value}) are optimistic approximations to Constraints (\ref{eq.m1.charge}).
We choose the negative function $-\underline{f}_S^k(\rho_i)$ of tangent  ${f}_S^k(\rho_i)$ (see Equation (\ref{eq.mc.secant})) to locally approximate concave item $-f(\rho_i)$, as $-\underline{f}_S^k(\rho_i)$ yields a downward approximation of $-f(\rho_i)$ within $[\varrho_i^k, \varrho_i^{k+1}]$. 
It always holds that $-\sum_{k\in\mathcal{K}_i} \underline{f}_S^k(\rho_i) \underline{\beta}_i^k \leq -f(\rho_i)$ due to the convexity of $f(\rho_i)$. 
Similar to $[M_C]$, we can introduce continuous $\underline{\pi}_i^k$ and Constraints (\ref{eq.mr.charge})--(\ref{eq.mr.value}) to linearize $\rho_i \underline{\beta}_i^k$ caused by $\underline{f}_S^k(\rho_i) \underline{\beta}_i^k$. 
This proves that Constraints (\ref{eq.mr.charge})--(\ref{eq.mr.value}) are optimistic approximations to Constraints (\ref{eq.m1.charge}).

Next, we prove that any feasible of $[M^*]$ are feasible for $[M_R]$.
From Proposition \ref{Pro.2}, optimistic approximations of nonlinear constraints yield a relaxation. 
$[M_R]$ is a relaxation of $[M^*]$. 
Therefore, any feasible solution of $[M_R]$ is feasible for $[M^*]$.
This completes the proof. \hfill \Halmos

\subsubsection{Proof of Proposition \ref{Pro.6}} \label{proof.Pro.6}~

We first prove that $Z[M_C] \geq Z[M^*] \geq Z[M_R]$.
From Propositions \ref{Pro.4} and \ref{Pro.5}, we know that $[M^*]$ is a relaxation of $[M_C]$, and $[M_R]$ is a relaxation of $[M^*]$. 
Therefore, $Z[M^*]$ can offer a lower bound on $Z[M_C]$, and $Z[M_R]$ can offer a lower bound on $Z[M^*]$, i.e., $Z[M_C] \geq Z[M^*] \geq Z[M_R]$.


To prove $Z[M_C'] \geq Z[M_C]$, we can first prove that $[M_C]$ is a relaxation of $[M_C']$. 
This can be proved if the following properties are satisfied:
\begin{enumerate}[fullwidth,itemindent=2em,label=(\alph*)]
    \item The objective values of $[M_C']$ and $[M_C]$ are identical;
    \item Any feasible solutions of $[M_C']$ are also feasible for $[M_C]$.
\end{enumerate}

Property (a) is obviously satisfied.
Therefore, we focus on proving Property (b).
Let us consider a feasible solution $(\overline{\bm{\Theta}}', {\bm{\rho}}', \bm{\nu}')$ of $[M_C']$. 
As compared to $[M_C]$, Constraints (\ref{eq.mc.Theta}), and (\ref{eq.mc.charge})--(\ref{eq.mc.value}) are different due to the different discretization schemes. 
Here we focus on that a feasible solution $(\overline{\bm{\Theta}}', {\bm{\rho}}', \bm{\nu}')$ of $[M_C']$ satisfies Constraints (\ref{eq.mc.Theta}) in $[M_C]$.
For a given vertiport $i\in\mathcal{N}$, let $\tilde{k}'$ be such that $\varrho_i^{\tilde{k}'} \leq \rho_i' \leq \varrho_i^{\tilde{k}'+1}$ and $\tilde{k}'\in\mathcal{K}_i'$.
Let $\tilde{k}$ be such that $\varrho_i^{\tilde{k}} \leq \rho_i' \leq \varrho_i^{\tilde{k}+1}$ and $\tilde{k}\in\mathcal{K}_i$.
Since $\{\varrho_i^{k'}\}_{k'\in\mathcal{K}_i'} \subseteq \{\varrho_i^{k}\}_{k\in\mathcal{K}_i}$, we have: $\varrho_i^{\tilde{k}'} \leq \varrho_i^{\tilde{k}} \leq \rho_i' \leq \varrho_i^{\tilde{k}+1} \leq \varrho_i^{\tilde{k}'+1}$.
For $[M_C']$, the bounded secant for $\overline{\Theta}_i'$ is $\overline{f}_{S'}^{\tilde{k}'}(\rho_i')$ determined by points $(\varrho_i^{\tilde{k}'},f(\varrho_i^{\tilde{k}'}))$ and $(\varrho_i^{\tilde{k}'+1},f(\varrho_i^{\tilde{k}'+1}))$.
For $[M_C]$, at $\rho_i'$, the bounded secant for $\overline{\Theta}_i'$ is $\overline{f}_S^{\tilde{k}}(\rho_i')$ determined by $(\varrho_i^{\tilde{k}},f(\varrho_i^{\tilde{k}}))$ and $(\varrho_i^{\tilde{k}+1},f(\varrho_i^{\tilde{k}+1}))$. 
From the convexity of $f(\rho_i)$, we have: $\overline{f}_{S'}^{\tilde{k}'}(\rho_i') \geq \overline{f}_S^{\tilde{k}}(\rho_i')$ for all $\rho_i' \in [\varrho_i^{\tilde{k}'}, \varrho_i^{\tilde{k}'+1}]$ (this can be proved using the property of convex functions).
This means that for all $i\in\mathcal{N}$, if $\overline{\Theta}_i' \geq \overline{f}_{S'}^{\tilde{k}'}(\rho_i')$, it always hold that $\overline{\Theta}_i' \geq \overline{f}_{S}^{\tilde{k}}(\rho_i')$. 
This completes proof that $(\overline{\bm{\Theta}}', {\bm{\rho}}', \bm{\nu}')$ satisfies Constraints (\ref{eq.mc.Theta}) in $[M_C]$. 

Here we focus on proving that a feasible solution of $[M_C']$ satisfies Constraints (\ref{eq.mc.charge}) in $[M_C]$.
We ignore $\overline{\bm{\pi}}$ in linearizing $\overline{\bm{\beta}}\bm{\rho}$ and remain the nonlinear form of Constraints (\ref{eq.mc.fleet}) as $\mathcal{D}\bm{\nu}' \leq \overline{\bm{f}}_T(\bm{\rho})\overline{\bm{\beta}}$.
For a given vertiport $i\in\mathcal{N}$, let $\tilde{k}'$ be such that $\overline{\varrho}_i^{\tilde{k}'} \leq \rho_i' \leq \overline{\varrho}_i^{\tilde{k}'+1}$ and $\tilde{k}'\in\mathcal{K}_i'$, i.e., $\overline{\beta}_i^{\tilde{k}'}=1$.
Let $\tilde{k}$ be such that $\overline{\varrho}_i^{\tilde{k}} \leq \rho_i' \leq \overline{\varrho}_i^{\tilde{k}+1}$ and $\tilde{k}\in\mathcal{K}_i$, i.e., $\overline{\beta}_i^{\tilde{k}}=1$.
Similarly, we can get $\overline{\varrho}_i^{\tilde{k}'} \leq \overline{\varrho}_i^{\tilde{k}} \leq \rho_i' \leq \overline{\varrho}_i^{\tilde{k}+1} \leq \overline{\varrho}_i^{\tilde{k}'+1}$.
Let row vector $\mathcal{D}_i$ represent the $i$-th row of $\mathcal{D}$.
For $[M_C']$, at $\rho_i'$, the bounded tangent for $\mathcal{D}_i\bm{\nu}'$ is $-\overline{f}_{T'}^{\tilde{k}'}(\rho_i')$ within $[\overline{\varrho}_i^{\tilde{k}'}, \overline{\varrho}_i^{\tilde{k}'+1}]$. 
For $[M_C]$, at $\rho_i'$, the bounded tangent for $\mathcal{D}_i\bm{\nu}'$ is $-\overline{f}_T^{\tilde{k}}(\rho_i')$ within $[\overline{\varrho}_i^{\tilde{k}}, \overline{\varrho}_i^{\tilde{k}+1}]$. 
From the convexity of $f(\rho)$, its piecewise tangents satisfy that $\overline{f}_{T'}^{\tilde{k}'}(\rho_i') \leq \overline{f}_{T}^{\tilde{k}}(\rho_i')$ for  $\rho_i\in[\overline{\varrho}_i^{\tilde{k}},\overline{\varrho}_i^{\tilde{k}+1}]$ (this can be proved using the property of concave functions). 
Therefore, if $\mathcal{D}_i\bm{\nu}' \leq \overline{f}_{T'}^{\tilde{k}'}(\rho_i') \overline{\beta}_i^{\tilde{k}'}$, it always holds that $\mathcal{D}_i\bm{\nu}' \leq \overline{f}_{T}^{\tilde{k}}(\rho_i')\overline{\beta}_i^{\tilde{k}}$.
Linearizing $\overline{\bm{\beta}}\bm{\rho}$ does not change this conclusion.
This completes the proof of that a feasible solution of $[M_C']$ satisfies Constraints (\ref{eq.mc.charge}) in $[M_C]$. 
So far, Property (b) is proved to be valid in our conservative model.
Finally, this proves that $Z[M_C'] \geq Z[M_C]$.

Similar to the proof of $Z[M_C'] \geq Z[M_C]$, we can prove that $Z[M_R] \geq Z[M_R']$ if the following properties are satisfied:
\begin{enumerate}[resume,fullwidth,itemindent=2em]
    \item[(c)] The objective values of $[M_R]$ and  $[M_R']$ are identical;
    \item[(d)] Any feasible solutions of $[M_R]$ are also feasible for $[M_R']$.
\end{enumerate}

Property (c) is straightforward and let us focus on proving (d).
Let us consider a feasible solution $(\underline{\bm{\Theta}}, {\bm{\rho}}, \bm{\nu})$ of $[M_R]$. 
As compared to $[M_R']$, Constraints (\ref{eq.mr.Theta}), and (\ref{eq.mr.charge})--(\ref{eq.mr.value}) in $[M_R]$ are different due to the different discretization schemes. 
We first focus on the proof that $(\underline{\bm{\Theta}}, {\bm{\rho}}, \bm{\nu})$ of $[M_R]$ is feasible for Constraints (\ref{eq.mr.Theta}).
For a given vertiport $i\in\mathcal{N}$, let $k$ be such that $\overline{\varrho}_i^k \leq \rho_i \leq \overline{\varrho}_i^{k+1}$. 
Let $k'$ be such that $\overline{\varrho}_i^{k'} \leq \rho_i \leq \overline{\varrho}_i^{k'+1}$. 
For $[M_R]$, the bounded tangent for $\underline{{\Theta}}_i$ is $\underline{f}_T^k(\rho_i)$ within $[\overline{\varrho}_i^k, \overline{\varrho}_i^{k+1}]$. 
For $[M_R']$, the bounded tangent for $\underline{{\Theta}}_i$ is $\underline{f}_T^{k'}(\rho_i)$ within $[\overline{\varrho}_i^{k'}, \overline{\varrho}_i^{k'+1}]$. 
From the convexity of $f(\rho)$, we know that $\underline{f}_T^{k}(\rho_i) \geq \underline{f}_T^{k'}(\rho_i)$ within $[\overline{\varrho}_i^{k'}, \overline{\varrho}_i^{k'+1}]$.
If  $\underline{{\Theta}}_i \geq \underline{f}_T^{k}(\rho_i)$, it always hold that $\underline{{\Theta}}_i \geq \underline{f}_T^{k'}(\rho_i)$ within $[\overline{\varrho}_i^{k'}, \overline{\varrho}_i^{k'+1}]$.
This completes the proof that $(\underline{\bm{\Theta}}, {\bm{\rho}}, \bm{\nu})$ is feasible for Constraints (\ref{eq.mr.Theta}) in $[M_R']$. 

To prove the feasibility of $(\underline{\bm{\Theta}}, {\bm{\rho}}, \bm{\nu})$ to Constraints (\ref{eq.mr.charge})--(\ref{eq.mr.value}), we also ignore $\underline{\pi}$.
The nonlinear form of Constraints (\ref{eq.mr.charge}) can be formulated as: $\mathcal{D}\bm{\nu} \leq -\underline{\bm{f}}_S(\bm{\rho}) \underline{\bm{\beta}}$.
For a given vertiport $i\in\mathcal{N}$, let $k$ be such that ${\varrho}_i^k \leq \rho_i \leq {\varrho}_i^{k+1}$. 
Let $k'$ be such that ${\varrho}_i^{k'} \leq \rho_i \leq {\varrho}_i^{k'+1}$. 
For $[M_R]$, the bounded secant for $\mathcal{D}\bm{\nu}$ is $-\underline{f}_S^k(\rho_i)$ within $[{\varrho}_i^k, {\varrho}_i^{k+1}]$.
For $[M_R']$, the bounded secant for $\mathcal{D}\bm{\nu}$ is $-\underline{f}_S^{k'}(\rho_i)$ within $[{\varrho}_i^{k'}, {\varrho}_i^{k'+1}]$.
From the convexity of $f(\rho)$, we can get $\underline{f}_S^{k}(\rho_i) \leq \underline{f}_S^{k'}(\rho_i)$ for $\rho_i\in [{\varrho}_i^k, {\varrho}_i^{k+1}]$.
Since $(\underline{\bm{\Theta}}, {\bm{\rho}}, \bm{\nu})$ is feasible for $[M_R]$, we have $\mathcal{D}_i\bm{\nu} \leq \underline{f}_S^{k}({\rho_i}) \underline{{\beta}_i^k}$.
Therefore, it always holds that $\mathcal{D}_i\bm{\nu} \leq \underline{f}_S^{k'}({\rho_i}) \underline{{\beta}_i^{k'}}$.
Linearizing Constraints (\ref{eq.mr.rhobeta1}) does not effect this conclusion.
Therefore, the feasible solution $(\underline{\bm{\Theta}}, {\bm{\rho}}, \bm{\nu})$ of $[M_R]$ is feasible for $[M_R']$. 
This means that Property (d) holds, and $[M_R']$ is a relaxation of $[M_R]$.
We can finally get $Z[M_R]\geq Z[M_R']$.
This completes the proof. 
\hfill\Halmos

\subsubsection{Proof of Theorem \ref{theorem.2}}~

Let us define $\varepsilon>0$ and define $(\bm{\rho}^R,\bm{\nu}^R)$ as a $[M_R]$ optimum.
From this solution, we can retrieve a feasible solution of $[M^*]$.
We fix $\bm{\rho}^R$ and substitute it into the general model (\ref{eq.general.model}).
The model can be reformulated as a MILP model as follows:
\begin{eqnarray}
    &\underset{{\bm{\nu}}}{\min} & \bm{c}^{\top} {\bm{\nu}} \nonumber\\
    &\mathrm{s.t.}& \bm{A}{\bm{\nu}} \leq {\bm{b}},\label{eq.general.mstar}
\end{eqnarray}
where
\begin{eqnarray}
    \bm{A} = \left[ \begin{array}{c}
	\bm{H}\\
	\bm{g}^{\top}\\
	\bm{D}\\
    \end{array} \right], \quad \quad 
    {\bm{b}} = \left[ \begin{array}{c}
	\bm{d} - \bm{F}\bm{\rho}^R\\
	-\bm{f}(\bm{\rho}^R)^{\top} \mathbf{1}\\
	\bm{f}(\bm{\rho}^R)\\
    \end{array} \right]. \nonumber
\end{eqnarray}

Suppose that there exits a non-degenerate primal optimal basic solution for Model (\ref{eq.general.mstar}).
Let matrix ${\bm{B}}$ and vector ${\bm{\nu}}_{\bm{B}}$ represent the basis matrix and basic variables for above model.
Then, its optimal solution $\tilde{\bm{\nu}}$ is given by $\tilde{\bm{\nu}}=[{\bm{\nu}}_{{\bm{B}}};{\bm{\nu}}_{{\bm{N}}}]$, where ${\bm{\nu}}_{{\bm{B}}}={\bm{B}}^{-1}{\bm{b}}$ and ${\bm{\nu}}_{{\bm{N}}} = \bm{0}$.
By construction, we can get a feasible solution $(\bm{\rho}^R,\tilde{\bm{\nu}})$ for $[M^*]$.

Next, we will prove that the solution $(\bm{\rho}^R,\tilde{\bm{\nu}})$ lies within $\varepsilon$ of $[M^*]$ optimum. 
The objective function value of $[M^*]$ at $(\bm{\rho}^R,\tilde{\bm{\nu}})$ satisfies:
\begin{eqnarray}
    \bm{c}^{\top} \tilde{\bm{\nu}} &=& \bm{c}_{\bm{B}}^{\top} {\bm{B}}^{-1}{\bm{b}} \nonumber\\ 
    &=& \bm{c}_{\bm{B}}^{\top} {\bm{B}}^{-1} 
    \left[\begin{array}{c}
         \bm{b} - \bm{F}\bm{\rho}^R \\
         -\underline{\bm{f}}_T(\bm{\rho}^R)^{\top} \bm{1}\\
         \underline{\bm{f}}_S(\bm{\rho}^R) 
    \end{array}\right] + \bm{c}_{\bm{B}}^{\top} {\bm{B}}^{-1} 
    \left[\begin{array}{c}
         0 \\
         \left( \underline{\bm{f}}_T(\bm{\rho}^R) -\bm{f}(\bm{\rho}^R) \right)^{\top}\bm{1}\\
         \bm{f}(\bm{\rho}^R) - \underline{\bm{f}}_S(\bm{\rho}^R)
    \end{array}\right]. \nonumber
\end{eqnarray}

$(\bm{\rho}^R,\bm{\nu}^R)$ is an optimal solution of $[M_R]$.
Let us fix $\bm{\rho}^R$ in $[M_R]$ and separate the tractable part as follows: 
\begin{eqnarray}
    &\underset{{\bm{\nu}^R}}{\min} & \bm{c}^{\top} {\bm{\nu}^R}  \label{eq.general.mrr}\\
    &\mathrm{s.t.}& \bm{A}{\bm{\nu}^R} \leq {\bm{b}^R},\nonumber
\end{eqnarray}
where
\begin{eqnarray}
    {\bm{b}^R} = \left[ \begin{array}{c}
	\bm{b} - \bm{F}\bm{\rho}^R\\
	-\underline{\bm{f}}_T(\bm{\rho}^R)^{\top} \mathbf{1}\\
	\underline{\bm{f}}_S(\bm{\rho}^R)\\
    \end{array} \right]. \nonumber
\end{eqnarray}

It can be seen that Model (\ref{eq.general.mrr}) is quite similar to Model (\ref{eq.general.mstar}) except the RHS vectors ${\bm{b}}$ and ${\bm{b}^R}$.
Suppose that at $\bm{\rho}^R$, $\underline{\bm{f}}_{T}(\bm{\rho}^R)$ (respectively, $-\underline{\bm{f}}_{S}(\bm{\rho}^R)$) approximates $\bm{f}(\bm{\rho})$ well and the difference between them is sufficient small such that the optimal basis, i.e., ${\bm{B}}$ for model (\ref{eq.general.mstar}) and (\ref{eq.general.mrr}) are the same, and the reduced cost $\bm{c}^{\top} {\bm{B}}^{-1}$ is not affected by the slight difference between RHS vectors.
An optimal solution of $[M_R]$ is given by ${\bm{\nu}}^R=[{\bm{\nu}}_{{\bm{B}}}^R;{\bm{\nu}}_{{\bm{N}}}^R]$, where ${\bm{\nu}}_{{\bm{B}}}^R={\bm{B}}^{-1}{\bm{b}}^R$ and ${\bm{\nu}}_{{\bm{N}}}^R = \bm{0}$.
The optimal objective function value of $[M_R]$ is given by $\bm{c}^{\top}\bm{\nu}^R=\bm{c}_{\bm{B}}^{\top} {\bm{B}}^{-1}\bm{b}^R$.
Thus, we can obtain:
\begin{eqnarray}
    \bm{c}^{\top}\tilde{\bm{\nu}}  
    &=& \bm{c}^{\top} {\bm{\nu}}^R
    +  \bm{c}_{\bm{B}}^{\top} {\bm{B}}^{-1} 
    ({\bm{b}} - \bm{b}^R) \nonumber \\
    &\leq & \bm{c}^{\top}\bm{\nu}^R + \Vert \bm{c}_{\bm{B}}^{\top} {\bm{B}}^{-1} \Vert_2 \cdot 
    \left \Vert \bm{b} - {\bm{b}}^R\right\Vert_2, \nonumber
\end{eqnarray}
where
\begin{eqnarray}
    &\left\Vert \bm{b} - {\bm{b}^R}\right\Vert_2 &= \sqrt{(\bm{b} - {\bm{b}}^R)^{\top}\cdot(\bm{b} -  {\bm{b}}^R)} \label{eq.b-br}\\
    &&= \sqrt{\left(\sum_{i\in\mathcal{N}}\left(f({\rho}_i^R) - \underline{f}_T({\rho}_i^R)\right)\right)^2 + \sum_{i\in\mathcal{N}}\left(\underline{f}_S({\rho}_i^R)-f({\rho}_i^R))\right)^2}.\nonumber
\end{eqnarray}

Let us focus on $\underline{f}_S({\rho}_i^R)-f({\rho}_i^R)$. 
Let the maximal diameter of discretization $\underset{i\in\mathcal{N}, k\in\mathcal{K}_i}{\max}\{\varrho_i^{k+1}-\varrho_i^{k},\overline{\varrho}_i^{k+1}-\overline{\varrho}_i^{k}\}=\Delta$.
Let $k$ be such that $\varrho_i^{k} \leq \rho_i^R \leq \varrho_i^{k+1}$.
At $\rho_i^R$, we have:
\begin{eqnarray}
    \underset{\rho_i\in[\varrho_i^{k},\varrho_i^{k+1}]}{\inf} f(\rho_i) \leq f(\rho_i^R) \leq \underline{f}_{S}(\rho_i^R) \leq \underset{\rho_i\in[\varrho_i^{k},\varrho_i^{k+1}]}{\sup} f(\rho_i) \nonumber.
\end{eqnarray}
From the Lipschitz continuity and the monotonicity of $f(\rho)$, we can get:
\begin{eqnarray}
    \underline{f}_{S}(\rho_i^R)-f(\rho_i^R) &\leq& \underset{\rho_i\in[\varrho_i^{k},\varrho_i^{k+1}]}{\sup} f(\rho_i)-\underset{\rho_i\in[\varrho_i^{k},\varrho_i^{k+1}]}{\inf} f(\rho_i) \nonumber \\
    &\leq& f(\varrho_i^{k+1})-f(\varrho_i^{k})  \leq L\Delta \nonumber.
\end{eqnarray}

Regarding $f({\rho}_i^R) - \underline{f}_T({\rho}_i^R)$, let $k$ be such that $\overline{\varrho}_i^{k} \leq \rho_i^R \leq \overline{\varrho}_i^{k+1}$.
As shown in Figure \ref{fig.exampleOfDiscretization}, the tangent point $\varrho^k$ for activated tangent line $\underline{f}_T(\cdot)$ lies within $[\overline{\varrho}_i^{k},\overline{\varrho}_i^{k+1}]$. 
$f({\rho}_i^R) - \underline{f}_T({\rho}_i^R)$ is indeed the error between $f(\rho)$ and its tangent approximation.
It holds that: 
\begin{eqnarray}
    f({\rho}_i^R) - \underline{f}_T({\rho}_i^R) \leq \max\{f({\overline{\varrho}_i^{k}}) - \underline{f}_T(\overline{\varrho}_i^{k}), f({\overline{\varrho}_i^{k+1}}) - \underline{f}_T(\overline{\varrho}_i^{k+1})\}.\nonumber
\end{eqnarray}
We can get:
\begin{eqnarray}
    \underset{\rho_i\in[\overline{\varrho}_i^{k},\overline{\varrho}_i^{k+1}]}{\sup}\{f({\overline{\varrho}_i^{k}}) - \underline{f}_T(\overline{\varrho}_i^{k})\} 
    & \leq &   f({\overline{\varrho}_i^{k}}) - \frac{\overline{\varrho}_i^{k} - (\overline{\varrho}_i^{k+1})^2}{(1-\overline{\varrho}_i^{k+1})^2} \nonumber \\
    & \leq & \frac{1}{(1-\overline{\varrho}_i^{k+1})^2} \Delta \leq L\Delta, \nonumber\\
    \underset{\rho_i\in[\overline{\varrho}_i^{k},\overline{\varrho}_i^{k+1}]}{\sup}\{f({\overline{\varrho}_i^{k+1}}) - \underline{f}_T(\overline{\varrho}_i^{k+1})\} 
    & \leq &  f({\overline{\varrho}_i^{k+1}}) - \frac{\overline{\varrho}_i^{k+1} - (\overline{\varrho}_i^{k})^2}{(1-\overline{\varrho}_i^{k})^2}  \nonumber\\
    & \leq & \frac{1}{(1-\overline{\varrho}_i^{k})^2} \Delta \leq L\Delta. \nonumber
\end{eqnarray}
Therefore, we have:
\begin{eqnarray}
    f({\rho}_i^R) - \underline{f}_T({\rho}_i^R) \leq L\Delta, \forall i\in\mathcal{N}. \nonumber
\end{eqnarray}

Recalling (\ref{eq.b-br}), we can obtain: 
\begin{eqnarray}
    \left\Vert \bm{b}^R - {\bm{b}}\right\Vert_2 & \leq & \sqrt{\left(\sum_{i=1}^{|\mathcal{N}|} L\Delta\right)^2 + \sum_{i=1}^{|\mathcal{N}|}(L\Delta)^2} \nonumber \\
    &=& \sqrt{(|\mathcal{N}|L\Delta)^2 + |\mathcal{N}|(L\Delta)^2} \nonumber\\
    &=& L\Delta \sqrt{|\mathcal{N}|^2+|\mathcal{N}|}. \nonumber
\end{eqnarray}
If $L\Delta \sqrt{|\mathcal{N}|^2+|\mathcal{N}|}$ is less than $\frac{\varepsilon}{\Vert \bm{c}_{\bm{B}}^{\top} {\bm{B}}^{-1} \Vert_2}$, i.e., 
\begin{eqnarray}
    \Delta \leq \frac{\varepsilon}{L \sqrt{|\mathcal{N}|^2 + |\mathcal{N}|} \Vert \bm{c}_{\bm{B}}^{\top} {\bm{B}}^{-1} \Vert_2}, \nonumber
\end{eqnarray}
we can get: $\bm{c}^{\top} \tilde{\bm{\nu}} \leq \bm{c}^{\top}\bm{\nu}^R + \varepsilon$.
Since the $[M_R]$ optimum provides a lower bound for $[M^*]$, it holds that $\bm{c}^{\top}\bm{\nu}^R \leq Z[M^*]$.
Finally, we can obtain $\bm{c}^{\top} \tilde{\bm{\nu}} \leq Z[M^*]+\varepsilon$. 
This means that if the constant $S=\frac{1}{L \sqrt{|\mathcal{N}|^2 + |\mathcal{N}|} \Vert \bm{c}_{\bm{B}}^{\top} {\bm{B}}^{-1} \Vert_2}$ and the diameter of discretization is less than $S\varepsilon$,  there exists a feasible solution $(\bm{\rho}^R,\tilde{\bm{\nu}})$ within $\varepsilon$ of the global optimum.
This completes the proof. 
\hfill\Halmos

\subsection{Algorithm Design}
\subsubsection{Adaptive Discretization Algorithm.}~ \label{sec.alg.exact}
The adaptive discretization algorithm is summarized in Algorithm \ref{alg.exact}.

\begin{algorithm}[H]
\caption{Adaptive discretization algorithm}
\label{alg.exact}
    \renewcommand{\KwIn}{\textbf{Parameter initialization: }}
    \renewcommand{\KwOut}{\textbf{Discretization initialization: }}
    \renewcommand{\Return}{\textbf{Return }}
    \SetKwBlock{Begin}{Step}{end}
    \SetKwFor{ForPar}{Step}{}{end foreach}%
    
    \KwIn{Upper bound $UB=\infty$, lower bound $LB=0$, optimality gap threshold $\epsilon$, iteration counter $r=0$, the maximum CPU time $T_{\max}$}
    
    \KwOut{Initialize the discretization point set $\bm{\varrho}$

    \BlankLine
    
    \textbf{Iterate between Step 1 and Step 2, until} $\frac{UB-LB}{LB} \leq \epsilon$ \textbf{or} CPU time $\geq T_{\max}$
    
    \ForPar{\textnormal{1. Improve the upper bound from $[M_C]$ solution:}}{
        \BlankLine
        \textbf{Step} 1.1. ($r \gets r+1$) Solve $[M_C]$ (with initial solution from neighborhood search phase in the last iteration if $r>1$); Get incumbent solution $(\bm{x}^C,\bm{\rho}^C)$
        
        \ForPar{\textnormal{1.2. Neighborhood search phase:}}{
            \BlankLine
            Solve $[M_C^N]$ with
                \begin{itemize}
                    \item ${x}_i = {x}_i^C, \;\; \forall i\in\mathcal{N}$
                    \item A granular discretization $\bm{\varrho}_i^L$, such that the interval $\Delta\leq0.05$ and $\rho_i^C\in\bm{\varrho}_i^N, \forall i\in\mathcal{N}$
                \end{itemize}
                
            Get the incumbent service level variable $\bm{\rho}^N$; Record $UB = Z[M_C^N]$
        }

        \textbf{Step} 1.3. Call the adaptive discretization procedure to generate new discretization points around $\bm{\rho}^N$; Update $\bm{\varrho}$
        }
    } 

    \ForPar{\textnormal{2. Improve the lower bound from $[M_R]$ solution:}}{

        \textbf{Step} 2.1. ($r \gets r+1$) Solve $[M_R]$ with initial solution from neighborhood search phase; Get incumbent solution $(\bm{x}^R,\bm{\rho}^R)$; Record $LB=Z[M_R]$
        
        \textbf{Step} 2.2. Call the adaptive discretization procedure to generate new discretization points around $\bm{\rho}^R$; Update $\bm{\varrho}$
    }

    \Return{\textnormal{The optimal incumbent found in $[M_C^N]$}}
\end{algorithm}

\subsubsection{Heuristic Algorithm.} \label{sec.heu.appendix}

We develop a matheuristic algorithm based on the philosophy of the exact algorithm.
The matheuristic algorithm inherits the iterative framework of the adaptive discretization algorithm. 
We develop an easy-to-solve surrogate model, denoted as $[M_M]$, which can substitute $[M_R]$ in exploring promising vertiports. 
The relaxed model $[M_C]$ is exclusively executed over the promising vertiport set identified by $[M_M]$, thereby reducing the computational time. 
Let $(\bm{x}^M,\bm{y}^M)$ denote an optimal solution of $[M_M]$. 
Vertiports selected by $[M_M]$ are regarded as promising and inserted into a \textit{promising vertiport set}, i.e., $\mathcal{N}^{\dag}: = \{i: x_i^M=1, i\in\mathcal{N}\}$.
\begin{eqnarray}\label{model.MaxCoverage}
    [{M_M}]&\max & \sum_{o,d\in\mathcal{R}}\sum_{i,j\in\mathcal{N}}D_{od}y_{o,i,j,d}^M, \nonumber\\
    &\rm{s.t.} & \sum_{i\in\mathcal{N}} x_i^M \leq P, \nonumber \\
    && y_{o,i,j,d}^M \leq x_{i}^M,\;\; y_{o,i,j,d}^M \leq x_{j}^M,\;\; y_{o,i,j,d}^M \leq \delta_{o,i,j,d}, \;\; \forall o,d\in\mathcal{R}, i,j\in\mathcal{N}, \nonumber\\ 
    && \sum_{i,j\in\mathcal{N}}y_{o,i,j,d}^M \leq 1, \;\; \forall o,d\in\mathcal{R}, \nonumber \\ 
    && \bm{x}^M, \bm{y}^M \text{ integer}. \nonumber
\end{eqnarray}

After solving $[M_M]$, we can execute $[M_C]$ over the promising vertiport set $\mathcal{N}^{\dag}$ instead of the entire candidate set $\mathcal{N}$.
The reduced candidate set can significantly save the computational time of $[M_C]$.
The matheuristic algorithm alternates between $[M_M]$ and $[M_C]$ and adds discretization point for $[M_C]$ according to the adaptive discretization procedure described in Section \ref{sec.procedure}. 
Relying solely on $[M_M]$ over $\mathcal{N^{\dag}}$ is unable to avoid local optima.
After the first iteration, we introduce the following constraint into $[M_M]$ to recognize new promising vertiports: 
\begin{eqnarray}
    \sum_{i\in\mathcal{N}}x_i^M(1-\omega_i) +  \sum_{i\in\mathcal{N}}(1-x_i^M)\omega_i \geq N_S\label{eq.neigh}, \nonumber
\end{eqnarray}
where $\omega_i$ is a parameter that equals 1 if vertiport $i$ is in the set $\mathcal{N}^{\dag}$, and 0 otherwise.
$N_S$ is the number of new promising vertiports. 
Let $\mathcal{N}^{\#}:=\{i:x_i^M=1,i\in\mathcal{N}\setminus\mathcal{N}^{\dag}\}$ denote the set of new identified vertiports.
Let $\mathcal{N}^{\dag}\gets \mathcal{N}^{\dag}\cup\mathcal{N}^{\#}$.
The enlarged set $\mathcal{N}^{\dag}$ and the implementation of the adaptive discretization procedure for $[M_C]$ ensure a decrease in the upper bound.
Algorithm \ref{alg.heu} summarizes the procedure of the proposed matheuristic algorithm. 
\begin{algorithm}[H]
\caption{Matheuristic algorithm}
\label{alg.heu}
    \renewcommand{\KwIn}{\textbf{Parameter initialization: }}
    \renewcommand{\KwOut}{\textbf{Discretization initialization: }}
    \renewcommand{\Return}{\textbf{Return }}
    \SetKwBlock{Begin}{Step}{end}
    \SetKwFor{ForPar}{Step}{}{end foreach}%
    
    \KwIn{Upper bound $UB^0=\infty$, iteration counter $r=0$, the threshold $\epsilon$ of improvement in upper bound, the maximum CPU time $T_{\max}$}
    
    \KwOut{Initialize the discretization point set $\bm{\varrho}=\varnothing$

    \BlankLine
    \textbf{Iterate between Step 1 and Step 2, until} $\frac{UB^{r}-UB^{r-1}}{UB^r} \leq \epsilon$ \textbf{or} CPU time $\geq T_{\max}$
    
    \ForPar{\textnormal{1. Generate the promising vertiport set using $[M_M]$:}}{
        \BlankLine
        
        \If{$r<1$}{Solve $[M_M]$ and get incumbent location variables $\bm{x}^M$
        
        Record the promising vertiport set $\mathcal{N}^{\dag}: = \{i: x_i^M=1, i\in\mathcal{N}\}$

        For all $i\in\mathcal{N}^{\dag}$, initialize discretization point set $\bm{\varrho}_i$
        }
        \Else{Solve $[M_M]$ with Constraint (\ref{eq.neigh}) and get incumbent location variables $\bm{x}^M$

        Record new identified vertiport set $\mathcal{N}^{\#}$; Let $\mathcal{N}^{\dag}\gets \mathcal{N}^{\dag} \cup \mathcal{N}^{\#}$

        For all $i\in\mathcal{N}^{\#}$, initialize discretization point set $\bm{\varrho}_i$
        }
        
        }
    } 

    \ForPar{\textnormal{2. Improve the upper bound from $[M_C]$:}}{

        \textbf{Step} 2.1. ($r \gets r+1$) Solve $[M_C]$ only on the current promising vertiport set $\mathcal{N}^{\dag}$; Get incumbent solution $(\bm{x}^C,\bm{\rho}^C)$; Record $UB^r=Z[M_C]$
        
        \textbf{Step} 2.2. Call the adaptive discretization procedure to generate new discretization points around $\bm{\rho}^C$; Update $\bm{\varrho}$
    }

    \Return{\textnormal{The optimal incumbent found in $[M_C]$}}
\end{algorithm}

\section{Details on Experimental Setup and Computational Results (Section \ref{sec.computionalResults})}

\subsection{Experimental Setup} \label{ec.expSetup}

\subsubsection{Dataset.}

According to the data provided by SF Express, there are a total of 535 service stations operating in Shenzhen, which are considered as the origin or destination of each parcel.
There are 148 service stations which have the potential to construct vertiports.

After data cleansing, there are 358,801 orders with a volumetric weight of at least one kg (i.e., the minimum volumetric weight), distributed across 33,960 O-D pairs.
We obtain the O-D demand by calculating the average volumetric weight per minute between 9:00 and 21:00.
After filtering out O-D pairs with a volumetric weight of less than five kg (which is unlikely to be profitable for drone-based logistics), we have identified a total of 6,801 O-D pairs. 
These pairs are distributed across 397 collection service stations and 414 distribution service stations.
Considering the long-term nature of facility network location decisions, we scale the demand by a factor of 1.5 to represent the forecasted demand for the next 10-year period.


\subsubsection{Parameter Setup.} 
In the benchmark case, we employ a logistics drone, H4, developed by SF Express, as the test drone. 
The H4 drone has a maximum payload capacity $Q$ of 12 kg, a cruise speed $v$ of 15 m/s, a coverage range $L^{\rm D}$ of 15 km, and it is priced at 130,800 CNY.
Accordingly, we calibrate the flying time $t_{ij}^f$ by calculating $l_{ij}/v + 1$, where $l_{ij}$ represents the flight distance between $i$ and $j$, and one minute accounts for the take-off and landing process.
The daily amortized cost of drone ownership is set to 71.67 CNY over a lifespan of five years.
We use fares of 0.51 CNY/(km$\cdot$kg) and 1.25 CNY/(km$\cdot$kg) to calibrate the costs of drone transport ($p_{ij}^{\mathrm{F}}$) and courier delivery ($p_{o,i}^{\mathrm{LT}}$ and $ p_{j,d}^{\mathrm{LT}}$) based on data provided by SF Express. 
The maximum service range of a vertiport ($L^{\mathrm{S}}$) is set to five km, beyond which there would be a significant cost for ground transport.
The benchmark market share $\xi$ is set to 20\%.

\subsubsection{Instance Generation.} 
Based on the dataset mentioned in above, we generate parameters in our testing instances as follows.
Let $S$ and $|\mathcal{N}|$ denote the number of O-D pairs and candidate vertiports, respectively.
We select the $S$ busiest O-D pairs with the highest round-trip demand.
The candidate locations of vertiports are chosen from 148 feasible locations based on the $|\mathcal{N}|$ busiest locations with the highest collection and distribution volumetric weight of parcels.
The configuration of the candidate facility network is then determined.
To test the robustness of proposed algorithm, we generate five random instances under each network configuration.
We scale the projected demand of each O-D pair by a random factor following a uniform distribution of $U(0.5, 1.5)$.
For the large-scale instances, all 6,801 O-D pairs are considered. 
We generate the candidate network based on the top-$S$ O-D pairs and the top-$|\mathcal{N}|$ candidate locations.

\subsection{SOCP Reformulation and Computational Results} \label{sec.socp}
\cite{he2021charging} propose a SOCP-based reformulation of convex function $f(\rho_i)$.
In our problem, we can use this reformulation to construct a conservative model $[M_C^S]$: 
\begin{eqnarray}
    [M_C^S]&& \left\Vert 
    \begin{array}{c}
        2 \\
        \overline{\theta}_i-(1-\rho_i)
    \end{array}
    \right\Vert_2 \leq \overline{\theta}_i+(1-\rho_i),\;\;\forall i \in \mathcal{N}, \nonumber \\
    && \sum_{i\in\mathcal{N}}(\overline{\theta}_i-1)  + \sum_{i,j\in \mathcal{N}}{t_{i,j}^{f}} \left(\psi_{i,j}+\varphi_{i,j}\right) \leq \Gamma, \nonumber\\
    && \text{Constraints } (\ref{eq.mc.charge})\text{--} (\ref{eq.mc.value}), \bar{\bm{\theta}} \text{ non-negative}, \nonumber
\end{eqnarray}
where $\overline{\theta}_i-1$ is an equivalent to $f(\rho_i)$. 
Table \ref{tab.socp} shows the computational comparison between a static run of $[M_C^S]$, $Sta(M_C^S)$, with our algorithm $Ada(M_C+M_C^N,MR)$.
Although the upper bound obtained from $Sta(M_C^S)$ is close to that of our algorithm, the CPU time is more than four times longer than our algorithm.
A small-scale instance, 30-300-12, which our algorithm can solve in 25 seconds, takes 3,600 seconds to solve using $Sta(M_C^S)$. 
Additionally, this reformulation only provides feasible solutions without optimality guarantees, while our algorithm provides both lower and upper bounds in only, on average, one-fourth of the CPU time.

\begin{table}[htbp]
  \centering
  \scriptsize
  \caption{Comparison between SOCP-based Reformulation and Our Algorithm}
  \begin{threeparttable} 
    \setlength{\tabcolsep}{5mm}{\begin{tabular}{lrrrr}
    \toprule
          & \multicolumn{2}{c}{$Sta(M_C^S)$} & \multicolumn{2}{c}{$Ada(M_C+M_C^N,MR)$} \\
        \cmidrule(r){2-3}  \cmidrule(r){4-5}          & \multicolumn{1}{r}{UB} & \multicolumn{1}{r}{CPU} & \multicolumn{1}{r}{UB} & \multicolumn{1}{r}{CPU} \\
    \midrule
    20-200-10 & 10.12 & 6     & 10.12 & 3 \\
    20-200-12 & 9.81  & 6     & 9.81  & 5 \\
    20-200-14 & 9.73  & 72    & 9.78  & 9 \\
    \cmidrule(r){1-1} \cmidrule(r){2-3}  \cmidrule(r){4-5}
    30-300-10 & 12.88 & 11    & 12.89 & 7 \\
    30-300-12 & 12.44 & 3,600  & 12.44 & 25 \\
    30-300-14 & 12.22 & 3,600  & 12.23 & 56 \\
    \cmidrule(r){1-1} \cmidrule(r){2-3}  \cmidrule(r){4-5}
    40-400-10 & 14.63 & 64    & 14.63 & 22 \\
    40-400-12 & 14.24 & 3,600  & 14.16 & 186 \\
    40-400-14 & 13.97 & 3,600  & 13.91 & 362 \\
    \cmidrule(r){1-1} \cmidrule(r){2-3}  \cmidrule(r){4-5}
    50-500-10 & 16.67 & 318    & 16.68 & 33 \\
    50-500-12 & 16.16 & 3,600  & 16.02 & 280 \\
    50-500-14 & 15.96 & 3,600  & 15.70 & 2,393 \\
    \cmidrule(r){1-1} \cmidrule(r){2-3}  \cmidrule(r){4-5}
    60-600-10 & 18.89 & 244    & 18.90 & 19 \\
    60-600-12 & 18.22 & 3,600  & 18.19 & 439 \\
    60-600-14 & 17.84 & 3,600  & 17.61 & 3,288 \\
    \cmidrule(r){1-1} \cmidrule(r){2-3}  \cmidrule(r){4-5}
    Total &   14.25    &  1,968   & 14.20 & 475 \\
    \bottomrule
    \end{tabular}} %
  \label{tab.socp}%
    \begin{tablenotes}
        \item \textit{Note}. UB, in thousand CNY; CPU, in CUP seconds.
    \end{tablenotes}
    \end{threeparttable}
\end{table}%

\subsection{Results of Matheuristic Algorithm for Large-scale Cases} 
\label{sec.matheuriscResult}

\subsubsection{Benchmark Approaches for Real-world Cases.} \label{ec.benchmark} 
To showcase the effectiveness of our matheuristic algorithm (Algorithm 2) in solving real-world cases, we compare it with $Sta(M_C)$ and two benchmark heuristic methods:
\begin{itemize}
    \item $Sta(M_C)$: Running $[M_C]$ statically with a discretization unit of 0.05; 
    \item $Sta(M_T^F\rightarrow M_T^S)$: This is a two-stage heuristic (see follows);
    \item $Ada(M_C,M_P)$: Using $[M_P]$ as the surrogate model in Algorithm 2;
    \item $Ada(M_C,M_M)$: Our matheuristic algorithm, Algorithm 2.
\end{itemize}

\textbf{Two-stage Heuristic.}
In the first stage model $[M_T^F]$, we optimize the vertiport location $\bm{x}$, service route $\bm{y}$, and some essential operational variables $\bm{\alpha}$, $\bm{\psi}$, and $\bm{\varphi}$. 
Let $\overline{\bm{x}}$ and $\overline{\bm{y}}$ represent the optimal solutions for $\bm{x}$ and $\bm{y}$ obtained from solving $[M_T^F]$.
\begin{eqnarray}
    [M_T^F] &\min & \Pi (\text{Equation} (\ref{eq.m1.obj})), \nonumber\\
    &\mathrm{s.t.} & \text{(\ref{eq.m1.P})--(\ref{eq.m1.y3}), (\ref{eq.m1.rhoz}), (\ref{eq.m1.demandLinear})--(\ref{eq.m1.alphay}), (\ref{eq.m1.flow}),}  \nonumber\\
    && \bm{\alpha}, \bm{\psi}, \bm{\varphi} \text{ no-negative};
    \bm{x}, \bm{y} \text{ binary}.\nonumber
\end{eqnarray}
By fixing the locations as $\overline{\bm{x}}$ and service routes as $\overline{\bm{y}}$, we formulate the second stage model $[{M_T^S}]$. 
\begin{eqnarray}
    [{M_T^S}] &\min & \Pi (\text{Equation} (\ref{eq.m1.obj})), \nonumber\\
    &\mathrm{s.t.} & \text{(\ref{eq.m1.P})--(\ref{eq.m1.y3}), (\ref{eq.m1.rhoz}), (\ref{eq.m1.demandLinear})--(\ref{eq.m1.alphay}), (\ref{eq.m1.flow}),  (\ref{eq.m1.rhoCap1}), (\ref{eq.mc.Theta})--(\ref{eq.mc.value})}, \bm{x} = \overline{\bm{x}}, \bm{y} = \overline{\bm{y}},\nonumber\\ 
    && \bm{\alpha}, \bm{\rho}, \bm{\psi}, \bm{\varphi} \text{ no-negative};
    \Gamma \text{ integer}; \bm{x}, \bm{z}, \bm{y} \text{ binary}.\nonumber
\end{eqnarray}

\textbf{$P$-median model.} 
We develop another surrogate model, denoted as $[M_P]$, for comparison with $[M_M]$.
This model selects $P$ vertiports which minimize the total time for parcel collection and distribution.
Let $D_o^O=\sum_{d\in\mathcal{R}}D_{o,d}$ and $D_d^D=\sum_{o\in\mathcal{R}}D_{o,d}$  represent the total demand from origin $o$ and to destination $d$, respectively.
$t_{o,i}^g$ and $t_{j,d}^g$ are the collection time from origin $o$ to vertiport $i$ and distribution time from vertiport $j$ to destination $d$.
Let $x_i\in\{0,1\}$ be an integer variable indicating the construction of vertiport $i$.
The variable $y_{o,i}^O$ and $y_{j,d}^D$ is an integer variable that indicates the assignments of origin $o$ and destination $d$ to vertiports. 
\begin{eqnarray}\label{model.PMedian}
    [{M_P}]&\min & \sum_{o\in\mathcal{R}}\sum_{i\in\mathcal{N}}D_o^D t_{o,i}^g y_{o,i}^O + \sum_{d\in\mathcal{R}}\sum_{j\in\mathcal{N}}D_d^D t_{j,d}^g y_{j,d}^D, \nonumber\\
    &s.t.& \sum_{i\in\mathcal{N}} y_{o,i}^O = 1, \forall o\in\mathcal{R}, \nonumber\\
    && \sum_{j\in\mathcal{N}} y_{j,d}^D = 1, \forall d\in\mathcal{R}, \nonumber\\
    && y_{o,i}^O \leq x_i,\;\;y_{j,d}^D \leq x_j, \forall o\in\mathcal{R}, i\in\mathcal{N}, \nonumber\\ 
    && \sum_{i\in\mathcal{N}} x_i \leq P, \nonumber\\
    && \bm{x}, \bm{y}^O, \bm{y}^D \text{ binary}. \nonumber
\end{eqnarray}


\subsubsection{Results.}
The time limit for each approach was set to 7,200 seconds, with a single run for each iteration in iteration-based algorithms, i.e., $Ada(M_C,M_P)$ and $Ada(M_C,M_M)$, set to 500 seconds.
Let a string $|\mathcal{N}|$-$P$ denote the settings of an instance.

Table \ref{tab.matheuristic} reports the computational results of these methods. 
The main observation is that among the four methods compared, our matheuristic algorithm $Ada(M_C,M_{M})$ consistently achieves the best solution within a reasonable computational time, with an average CPU time of 1,760 seconds.
$Sta(M_C)$ may fail to find a high-quality solution or even a feasible solution.
The two-stage heuristic $Sta(M_T^F\rightarrow M_T^S)$ can result in an average inferiority of 15\%, and it may also produce infeasible solutions in some instances.
This result also highlights the importance of incorporating operational decisions into network planning. 
Replacing $[M_M]$ with $[M_P]$ in Algorithm 2 also results in a 10.8\% increase in average operating cost, showcasing the advantage of using $[M_M]$ to fulfill the role of $[M_R]$ in the adaptive discretization algorithm.
In conclusion, the proposed matheuristic demonstrates its ability to deliver high-quality solutions within a reasonable timeframe.

\begin{table}[H]
  \centering
  \scriptsize
  \caption{Comparison of Computational Results of Our Matheuristic Algorithm and Benchmark Methods}
    \begin{threeparttable}

    \begin{tabular}{lrrrrrrrrrrrr}
    \toprule
          & \multicolumn{3}{c}{$Sta(M_C)$} &     \multicolumn{3}{c}{$Sta(M_T^F\rightarrow M_T^S)$}  & \multicolumn{3}{c}{$Ada(M_C,M_P)$} &      \multicolumn{3}{c}{$Ada(M_C,M_{M})$}  \\
        \cmidrule(r){2-4}  \cmidrule(r){5-7} \cmidrule(r){8-10}   \cmidrule(r){11-13}  
        \multicolumn{1}{r}{$|\mathcal{N}|$-$P$} & \multicolumn{1}{r}{Obj} & \multicolumn{1}{r}{CPU} & \multicolumn{1}{r}{Gap} & \multicolumn{1}{r}{Obj} & \multicolumn{1}{r}{CPU} & \multicolumn{1}{r}{Gap} & \multicolumn{1}{r}{Obj} & \multicolumn{1}{r}{CPU} & \multicolumn{1}{r}{Gap} & \multicolumn{1}{r}{Obj} & \multicolumn{1}{r}{CPU} & \multicolumn{1}{r}{Gap} \\
        \midrule
    89-25 & {{76.65}} & 7,200  & 0.00  & {-} & 9     & 100   & {88.51} & 281   & 15.4  & {76.69} & 1,545  & 0.04 \\
    89-30 & {75.62} & 7,200  & 0.53  & {94.43} & 13    & 25.5  & {79.07} & 2,711  & 5.12  & {{75.22}} & 1,548  & 0.00 \\
    89-35 & {74.52} & 7,200  & 0.75  & {74.23} & 8     & 0.36  & {79.01} & 1,038  & 6.83  & {{73.96}} & 1,668  & 0.00 \\
    119-25 & {77.17} & 7,200  & 0.84  & {77.04} & 62    & 0.67  & {87.70} & 294   & 14.6  & {{76.53}} & 1,219  & 0.00 \\
    119-30 & {75.78} & 7,200  & 0.68  & {77.78} & 16    & 3.33  & {81.41} & 2,413  & 8.15  & {{75.27}} & 1,622  & 0.00 \\
    119-35 & {74.82} & 7,200  & 1.23  & {78.20} & 17    & 5.82  & {78.55} & 1,578  & 6.28  & {{73.90}} & 1,744  & 0.00 \\
    134-25 & {77.22} & 7,200  & 1.22  & {77.99} & 38    & 2.22  & {88.29} & 144   & 15.7  & {{76.29}} & 1,448  & 0.00 \\
    134-30 & {76.01} & 7,200  & 1.95  & {75.09} & 94    & 0.71  & {83.07} & 341   & 11.4  & {{74.56}} & 3,297  & 0.00 \\
    134-35 & {89.13} & 7,200  & 20.4  & {92.82} & 37    & 25.4  & {82.11} & 60    & 10.9  & {{74.04}} & 1,593  & 0.00 \\
    148-25 & {-} & 7,200  & 100   & {85.98} & 64    & 13.3  & {88.28} & 115   & 16.4  & {{75.86}} & 1,809  & 0.00 \\
    148-30 & {-} & 7,200  & 100   & {74.68} & 58    & 0.60  & {83.06} & 743   & 11.9  & {{74.24}} & 1,845  & 0.00 \\
    148-35 & {-} & 7,200  & 100   & {75.29} & 75    & 2.27  & {78.55} & 1,570  & 6.70  & {{73.62}} & 1,778  & 0.00 \\
    \midrule
    Total &       & 7,200  & 27.3  &       & 41    & 15.0  &       & 941   & 10.8  &       & 1,760  & 0.00 \\
    \bottomrule
    \end{tabular}%
    \begin{tablenotes}
        \item \textit{Note}. This table includes values of the upper bound (UB, in thousands CNY), optimality gap (Gap, in percentage terms), and computational time (CPU, in CUP seconds).
        If a method fails to obtain a feasible solution, the optimality gap is considered as 100\%.
    \end{tablenotes}
    \end{threeparttable}
  \label{tab.matheuristic}%
\end{table}%

\section{Details on Experimental Setting for Managerial Insights  (Section \ref{sec.insights})} 
\label{ec.insights}

\textbf{Mean Operating Cost and Lead Time.}
The mean operating cost per delivery can be obtained from dividing objective function value by the total served demand $D_{ij}^{\mathrm{T}}=\sum_{o,d\in\mathcal{R}}\alpha_{o,i,j,d}D_{od}$.
The lead time $t_{o,i,j,d}^{\mathrm{L}}$ for a service route $(o,i,j,d)$ consists of the transport time $t^{\mathrm{T}}$ on three itinerary legs and pooling time $t_{i,j}^{\mathrm{P}}$ for the parcels transported by drones on vertiport pair $(i,j)$, i.e., $t_{o,i,j,d}^{\mathrm{L}} = t_{o,i}^{\mathrm{T}} + t_{i,j}^{\mathrm{T}} + t_{j,d}^{\mathrm{T}} + t_{i,j}^{\mathrm{P}}$.
Following the probabilistic split method, the expected pooling time of a service route is formulated as $t_{i,j}^{\mathrm{P}} = \frac{Q-1}{D_{i,j}^{\mathrm{T}}}$ \citep{kai2022vertiport}. 
Finally, the mean lead time is calculated as the weighted average of lead time on each selected service route. 

\textbf{D2D--Courier.}
We formulate a model, $[M_{S}]$, which identifies the most promising O-D pairs to offer dedicated door-to-door delivery services.
$c^S$ is the basic daily wage for a dedicated courier.
$p_{o,d}^S$ is the cost for a delivery service from $o$ to $d$.
For convenient of comparison, we also allow the parcel pooling, i.e., multiple parcels within on delivery. 
Let $Q^S$ denote the pooling size. 

Let $y_{o,d}^S\in\{0,1\}$ be an integer variable that signifies whether an O-D pair $(o,d)$ is chosen for service.
Integer variable $\Gamma^S$ denote the staff size of dedicated couriers. 
Continuous variable $\rho_{o}^S$ represents the service level at origin $o\in\mathcal{R}$.
Continuous variables $\psi_{o,d}^S$ and $\varphi_{o,d}^S$ denote the average number of transit and repositioning couriers, respectively. 
The model can be formulated as follows: 
\begin{eqnarray}
    [M_{S}]&\min \Pi &= 
    \underset{\text{staff\ size}}{\underbrace{c^S\Gamma^S }} 
    + \underset{\text{courier\ transport\ cost}}{\underbrace {\sum_{o,d\in \mathcal{R}}{p_{o,d}^{S}} \left( \psi_{o,d}^S+\varphi_{o,d}^S \right) }} \label{eq.ms2s} \nonumber \\
    &\mathrm{s.t.}& \sum_{o,d\in\mathcal{R}} {\rho_o^S y_{o,d}^S} D_{od} \geq \xi \sum_{o,d\in\mathcal{R}} D_{od}, \nonumber \\
    && \psi_{o,d}^S =  \frac{D_{od}}{Q^S} {\rho_o y_{o,d}}, \;\; \forall o,d\in\mathcal{R}, \nonumber \\
    && \sum_{d\in\mathcal{R}} \left(\psi_{o,d}^S+\varphi_{o,d}^S \right) = \sum_{j\in\mathcal{N}}{\left(\psi_{c,j}^S+\varphi_{c,j}^S \right)},  \nonumber\\
    && \sum_{o\in\mathcal{R}} {\frac{\rho_o^S}{1-\rho_o^S}} +  \sum_{o,d\in\mathcal{R}} t_{o,d}^{g} \left(\psi_{o,d}^S+\varphi_{o,d}^S \right)
    \leq \Gamma^S, \nonumber\\
    && \bm{x}^S, \bm{y}^S \text{ binary}; \Gamma^S \text{ non-negative integer}; \bm{\psi}^S, \bm{\varphi}^S,\bm{\rho}^S \text{ non-negative}. \nonumber
\end{eqnarray}

\textbf{H\&S--Courier.} The model resembles $[M^*]$ but lacks the capacity constraints (\ref{eq.m1.Gamma}), (\ref{eq.m1.rhoCap1}), the service route feasibility constraints (\ref{eq.m1.y3}), and the battery charging constraints (\ref{eq.m1.charge}).
Furthermore, parameters such as drone flight costs $p_{i,j}^{\mathrm{F}}$ and flight times $t_{i,j}^f$ have to be substituted with the corresponding courier transport costs $p_{i,j}^{S}$ and travel times $t_{i,j}^g$.

\bibliographystyleannex{informs2014}
\bibliographyannex{mybibfile02}

\end{document}